\documentclass{article}

\usepackage[margin=1in]{geometry}

\usepackage[utf8]{inputenc} 
\usepackage[T1]{fontenc}    
\usepackage{hyperref}       
\usepackage{url}            
\usepackage{booktabs}       
\usepackage{amssymb, amsthm, amsfonts, amsmath, mathtools, thm-restate}       
\usepackage{nicefrac}       
\usepackage{microtype}      
\usepackage{bbm, bm}
\usepackage{subcaption}

\usepackage{algorithm}
\usepackage[noend]{algpseudocode}
\usepackage{varwidth}
\usepackage{tabularx}
\usepackage{makecell}
\usepackage{nicematrix}
\usepackage{pifont}
\usepackage{thm-restate}
\usepackage{enumitem}

\usepackage[font=small,labelfont=bf]{caption}
\usepackage{xfrac} 

\usepackage[
backref=true,
natbib=true,
style=numeric-comp,
giveninits=false, 
maxbibnames = 99, 
maxalphanames=8,
sorting=none
]{biblatex}
\usepackage{xcolor}

\definecolor{otherlightblue}{RGB}{0, 100, 200}
\definecolor{otherblue}{RGB}{0, 50, 100}
\definecolor{othergreen}{RGB}{60, 120, 0}

\hypersetup{
	colorlinks=true,
	pdfpagemode=UseNone,
	citecolor=othergreen,
	linkcolor=otherlightblue,
	urlcolor=magenta
}

\usepackage[capitalize, nameinlink]{cleveref}
\crefformat{equation}{(#2#1#3)}
\crefname{eq}{Equation}{Equations}
\crefname{fact}{Fact}{Facts}
\crefname{lemmma}{Lemma}{Lemmas}
\crefname{lemma}{Lemma}{Lemmas}
\crefname{figure}{Figure}{Figures}
\crefname{defn}{Definition}{Definitions}
\crefname{ineq}{Inequality}{Inequalities}
\crefname{prob}{Problem}{Problems}
\crefname{corollary}{Corollary}{Corollaries}
\creflabelformat{ineq}{#2{\upshape(#1)}#3} 
\crefalias{thm}{theorem}
\Crefname{listfact}{Fact}{Facts}
\newlist{assumpenum}{enumerate}{1} 
\setlist[assumpenum]{label=(\roman*), ref=\theproposition(\roman*), leftmargin = 2em, topsep=0em}
\crefalias{assumpenumi}{assumption} 
\crefname{assum}{Assumption}{Assumptions} 
\crefname{assumption}{Assumption}{Assumptions}
\crefname{section}{Section}{Sections}

\newcommand\numberthis{\addtocounter{equation}{1}\tag{\theequation}} 
\numberwithin{equation}{section}
\newcommand{\norm}[1]{\left\lVert#1\right\rVert}
\newcommand{\innerprod}[2]{\langle#1,#2\rangle}%
\global\long\def\grad{\nabla}%
\global\long\def\del{\delta}%
\DeclareMathOperator*{\argmin}{arg\,min}

\newcommand{\reals}{\mathbb{R}}
\newcommand{\R}{\mathbb{R}}
\newcommand{\NN}{\mathbb{N}}
\newcommand{\tbg}{\widetilde{{g}}}
\newcommand{\E}{\mathbb{E}}
\renewcommand{\S}{\mathbb{S}}

\theoremstyle{plain}
\newtheorem{theorem}{Theorem}[section]

\newtheorem{definition}[theorem]{Definition}

\newtheorem{lemma}[theorem]{Lemma}
\newtheorem{remark}[theorem]{Remark}
\newtheorem{corollary}[theorem]{Corollary}
\newtheorem{assumption}[theorem]{Assumption}

\newtheorem{proposition}[theorem]{Proposition}

\global\long\def\grad{\nabla}%
\newcommand{\ystar}{y^*} 
\global\long\def\ydelstar{y_{\delta}^{*}}%
\global\long\def\ystardel{y_{\delta}^{*}}%
\global\long\def\lamdeltar{\lam_{\delta}^{*}}%
\global\long\def\lam{\lambda}%
\global\long\def\lamstar{\lam^{*}}%
\newcommand{\Leqc}{\mathcal{L}_{\mathrm{eq}}}

\newcommand{\lamstarsmoothlineq}{C_{\ell}}
\newcommand{\ystarsmoothlineq}{C_y}
\newcommand{\ystarliplineq}{M_y}
\newcommand{\gtsmooth}{S_g} 
\newcommand{\gssmooth}{C_g}

\global\long\def\pq{\textrm{proj}_{\mathcal{X}}}%

\newcommand{\Unif}{\mathrm{Unif}}
\newcommand{\out}{\mathrm{out}}

\newcommand{\tF}{\widetilde{F}}

\newcommand{\inner}[1]{\left\langle#1\right\rangle}

\newcommand{\Comments}{1} 
\newcommand{\mynote}[2]{\ifnum\Comments=1\textcolor{#1}{#2}\fi}

\newcommand{\compresslist}{ 
\setlength{\itemsep}{1pt}
\setlength{\parskip}{0pt}
\setlength{\parsep}{0pt}
}

\title{First-Order Methods for \\ Linearly Constrained Bilevel Optimization\thanks{GK, SP, KW, and ZZ contributed equally; authors listed in alphabetical order.}} 

\author{%
  Guy Kornowski\thanks{Weizmann Institute of Science. \texttt{guy.kornowski@weizmann.ac.il}}
  \and
  Swati Padmanabhan\thanks{Massachusetts Institute of Technology. \texttt{pswt@mit.edu}}%
  \and
  Kai Wang\thanks{Georgia Institute of Technology. \texttt{kwang692@gatech.edu}}
  \and
  Zhe Zhang\thanks{Purdue University. \texttt{zhan5111@purdue.edu}}
  \and
  Suvrit Sra\thanks{Massachusetts Institute of Technology. \texttt{suvrit@mit.edu}}
}

\begin{document}

\maketitle

\begin{abstract}
Algorithms for bilevel optimization often encounter Hessian computations, which are prohibitive in high dimensions. While recent works offer first-order methods for {unconstrained} bilevel problems, the \textit{constrained} setting remains relatively underexplored. We present first-order linearly constrained optimization methods with finite-time hypergradient stationarity guarantees. For linear \textit{equality} constraints, we attain $\epsilon$-stationarity in $\widetilde{O}(\epsilon^{-2})$ gradient oracle calls, which is nearly-optimal. For linear \textit{inequality} constraints, we attain $(\delta,\epsilon)$-Goldstein stationarity in $\widetilde{O}(d{\delta^{-1} \epsilon^{-3}})$ gradient oracle calls, where $d$ is the upper-level dimension. Finally, we obtain for the linear inequality setting dimension-free rates of $\widetilde{O}({\delta^{-1} \epsilon^{-4}})$ oracle complexity under the additional assumption of oracle access to the optimal dual variable. Along the way, we develop new nonsmooth nonconvex optimization methods with inexact oracles. We verify these guarantees with preliminary numerical experiments.
\end{abstract}

\newpage

\section{Introduction}\label{sec:introduction}
Bilevel optimization~\cite{bracken1973mathematical, colson2007overview,bard2013practical,sinha2017review}, an important problem in optimization, is defined as follows:
\begin{align*}\numberthis\label[prob]{prob:general-constraints}
     \mbox{minimize}_{x\in \mathcal{X}} ~ F(x) \coloneqq f(x, \ystar(x)) \text{ subject to } ~\ystar(x)\in \arg\min\nolimits_{y\in S(x)} g(x, y).
\end{align*}
Here, the value of the upper-level problem at any point $x$ depends on the solution of the lower-level problem. This framework has recently found numerous applications in meta-learning~\cite{ snell2017prototypical, bertinetto2018meta, rajeswaran2019meta, ji2020convergence}, hyperparameter optimization~\cite{franceschi2018bilevel, shaban2019truncated, feurer2019hyperparameter}, and 
reinforcement learning~\cite{konda1999actor, sutton2018reinforcement, hong2020two,zhang2020bi}.
Its growing importance has spurred increasing efforts towards designing computationally efficient algorithms for it. 

As demonstrated by \citet{ghadimi2018approximation}, a key computational step in algorithms for bilevel optimization is estimating $\tfrac{dy^*(x)}{dx}$, the gradient of the lower-level solution. This gradient estimation problem has been extensively studied in differentiable optimization~\cite{amos2017optnet,agrawal2019differentiable} by applying the implicit function theorem to the KKT system of the given problem~\cite{donti2017task,wilder2019melding,kotary2021end,lee2019meta,tang2022pyepo,bai2019deep}.
However, this technique typically entails computing (or estimating) second-order derivatives, which can be prohibitive in high dimensions~\cite{mehra2021penalty, ji2021bilevel,wang2021learning}.

Recently, \citet{liu2022bome} made a big leap forward towards addressing this computational bottleneck.  
Restricting themselves to the class of unconstrained bilevel optimization, they proposed a  first-order method with finite-time stationarity guarantees. While a remarkable breakthrough,  \cite{liu2022bome} does not directly extend to the important setting of \textit{constrained} bilevel optimization. This motivates the question: 
\begin{quote}
\centering
    \emph{Can we develop a  first-order algorithm for constrained bilevel optimization?}
\end{quote}
\looseness=-1Besides being natural from the viewpoint of complexity theory, this question is well-grounded
in applications such as 
mechanism design~\cite{wang2022coordinating,dutting2021optimal},  
resource allocation~\cite{xu2013bilevel,gutjahr2016bi,zhang2010bilevel,fathollahi2022bi}, and decision-making under uncertainty~\cite{elmachtoub2022smart,munoz2022bilevel,wilder2019melding}. 
Our primary contribution is \textit{an affirmative answer to the highlighted question for bilevel programs with linear constraints}, an important problem class often arising in adversarial training, decentralized meta learning, and sensor networks~(see \cite{khanduri2023linearly} and the discussion therein). 
While there have been some other recent works~\cite{khanduri2023linearly, yao2024constrained, lu2024firstorder} on this problem, our work is first-order (as opposed to \cite{khanduri2023linearly}) and offers, in our view, a stronger guarantee on stationarity (compared to \cite{yao2024constrained, lu2024firstorder})--- cf. \cref{sec:related_work}. We hope our work paves the path to progress for first-order methods for general constrained bilevel programs.  
We now summarize our contributions.

\subsection{Our contributions}\label{sec:introduction_main_results}
We provide first-order algorithms (with associated finite-time convergence guarantees) for linearly constrained bilevel programs (\cref{prob:general-constraints}). By ``first-order'', we mean that we use only zeroth and first-order oracle access to $f$ and $g$.  
\begin{description}[style=unboxed,leftmargin=0cm, itemsep=.5em, parsep=.3em]
\item [{(1)}]
 \textbf{Linear equality constraints.}  As our first contribution, we design first-order algorithms for solving \cref{prob:general-constraints} where the lower-level constraint set $S(x):=\left\{y:Ax-By-b=0\right\}$  comprises linear equality constraints,  and $\mathcal{X}$ a convex compact set. 
With appropriate regularity assumptions on $f$ and $g$, we show in this case smoothness in $x$ of the hyperobjective $F$. Inspired by ideas from \citet{kwon2023fully}, we use implicit differentiation of the KKT matrix of a slightly perturbed version of the lower-level problem to design a first-order approximation to $\nabla F$. Constructing our first-order approximation entails solving a strongly convex optimization problem on affine constraints, which can be done efficiently. With this inexact gradient oracle in hand,  we then run projected gradient descent, which converges in $\widetilde{O}(\epsilon^{-2})$ iterations for smooth functions. 

\begin{theorem}[Informal; cf. \cref{thm:lineq-cost}]
\label{thm:lineq-cost-inf}
    Consider  \cref{prob:general-constraints} with linear equality constraints $S(x)\coloneqq\left\{y:Ax-By-b=0\right\}$ and $\mathcal{X}$ a convex compact set. Under mild regularity assumptions on  $f$ and $g$ (\cref{assumption:linEq_smoothness,assumption:eq}), there exists an algorithm, which in $\widetilde{O}(\epsilon^{-2})$ zeroth and first-order oracle calls to $f$ and $g$, converges to an $\epsilon$-stationary point of the hyperobjective $F$. 
\end{theorem}
For  linear equality constrained bilevel optimization, this is the first  first-order result attaining $\epsilon$-stationarity of $F$ with assumptions solely on the constituent functions $f$ and $g$ and none on $F$ --- cf. \cref{sec:related_work} for a discussion of the results of \citet{khanduri2023linearly} for this setting. 
 \item [{(2)}] \textbf{Linear inequality constraints.} Next, we provide first-order algorithms for solving \cref{prob:general-constraints} where the lower-level constraint set $S(x):=\left\{y:Ax-By-b\leq0\right\}$  comprises linear inequality constraints, and the upper-level variable is unconstrained. 

Our measure of convergence of  algorithms in this case is that of $(\delta,\epsilon)$-stationarity~\cite{goldstein1977optimization}: for a Lipschitz function, we say that a point $x$ is $(\delta,\epsilon)$-stationary if within a $\delta$-ball around $x$ there exists a convex combination of subgradients of the function with norm at most $\epsilon$ (cf. \cref{def:GoldsteinDeltaEpsStationary}). 

\looseness=-1
To motivate this notion of convergence, we note that the hyperobjective $F$ (in \cref{prob:general-constraints}) as a function of $x$ could be nonsmooth and nonconvex (and is Lipschitz, as  we later prove). Minimizing such a function in general is well-known to be intractable~\cite{nemirovskij1983problem}, necessitating local notions of stationarity. Indeed, not only is it impossible to attain $\epsilon$-stationarity in finite time~ \cite{zhang2020complexity}, even getting \textit{near} an
approximate stationary point of an arbitrary Lipschitz function is impossible unless
the number of queries has an exponential dependence on the dimension~\cite{kornowski2022oracle}. Consequently, for this function class, $(\delta, \epsilon)$-stationarity has recently emerged, from the work of \citet{zhang2020complexity}, to be a natural and algorithmically tractable notion of stationarity. We give the following guarantee under  regularity assumptions on $f$ and $g$.

\begin{theorem}[Informal; \cref{{{thm:izo_complete_guarantees}}}]\label{thm:informal_izo_full}
 Consider \cref{prob:general-constraints} with linear inequality constraints $S(x)\coloneqq\left\{y:Ax-By-b\leq 0\right\}$. Under mild regularity assumptions on $f$ and $g$ (\cref{assumption:linEq_smoothness}) and the lower-level primal solution $y^*$ (\cref{assumption:ineq_mild}), there exists an algorithm, which in $\widetilde{O}(d\delta^{-1}\epsilon^{-3})$  oracle calls to $f$ and $g$, converges to a $(\delta, \epsilon)$-stationary point of $F$, where $d$ is the upper-level variable dimension.  
\end{theorem}

To the best of our knowledge, this is the first result to offer a  first-order finite-time stationarity guarantee on the hyperobjective for 
linear inequality constrained bilevel optimization (cf. \cref{sec:related_work} for a discussion of related work of \citet{khanduri2023linearly, yao2024constrained, lu2024firstorder}).
We obtain our guarantee in \cref{thm:informal_izo_full} by first invoking a result by \citet{zhang2022solving} to obtain inexact hyperobjective values of $F$ using only $\widetilde{O}(1)$ oracle calls to $f$ and $g$. We also show (\cref{{lem:LipscConstrBilevel}}) that this hyperobjective $F$ is Lipschitz. We then employ our inexact zeroth-order oracle for  $F$ in an algorithm (\cref{alg: IZO}) designed to minimize Lipschitz nonsmooth nonconvex functions (in particular, $F$) with the following convergence guarantee. 

\begin{theorem}[Informal; \cref{{thm: Lipschitz-min-with-inexact-zero-oracle}}] 
Given $L$-Lipschitz $F:\reals^d\to\reals$ 
and that $|\widetilde{F}(\cdot)-F(\cdot)|\leq\epsilon$.
Then there exists an algorithm, which in 
$\widetilde{O}(d\delta^{-1}\epsilon^{-3})$
calls to $\tF(\cdot)$
outputs 
 a point $x^{\out}$ such that $\E[\mathrm{dist}(0,{\partial}_\delta F(x^{\out}))]\leq2\epsilon$. 

\end{theorem}

While such algorithms using  \emph{exact} zeroth-order access already exist~\cite{kornowski2024algorithm}, extending them to the inexact gradient  setting is non-trivial; we leverage recent ideas connecting online learning to nonsmooth nonconvex optimization by \citet{cutkosky2023optimal} (cf. \cref{{sec:nonsmooth}}). With the ubiquity of nonsmooth nonconvex  optimization problems associated with
training modern neural networks, we believe our analysis for this general task can be of independent interest to the broader optimization community.  

\item[{(3)}] \textbf{Linear inequality under assumptions on dual variable access.}
 For the {inequality} setting (i.e., \cref{prob:general-constraints} with the lower-level constraint set $S(x):=\left\{y:Ax-By-b\leq0\right\}$), we obtain dimension-free rates under an additional assumption  (\cref{item:assumption_safe_constraints}) on oracle access to the optimal dual variable $\lambda^*$ of the lower-level problem. We are not aware of a method to obtain this dual variable in a first-order fashion (though in practice, highly accurate approximations to $\lambda^*$ are readily available), hence the need for imposing this assumption. We believe that removing this assumption and obtaining dimension-free first-order rates in this setting would be a very interesting direction for future work. Our guarantee for this setting is summarized below.

\begin{theorem}[Informal; \cref{thm:Lipschitz-min-with-inexact-grad-oracle} combined with \cref{thm:cost_of_computing_ystar_gammastar_inequality}]\label{thm:linIneqInformal}
    Consider \cref{prob:general-constraints} with linear inequality constraints $S(x)\coloneqq\left\{y:Ax-By-b\leq0\right\}$ and unconstrained upper-level variable. Under mild regularity assumptions on $f$ and $g$ (\cref{assumption:linEq_smoothness}), on $y^*$ (\cref{assumption:ineq_mild}), and assuming oracle access to the optimal dual variable $\lambda^*$ (\cref{item:assumption_safe_constraints}), there exists an algorithm, which in $\widetilde{O}(\delta^{-1}\epsilon^{-4})$ oracle calls to $f$ and $g$ converges to a $(\delta, \epsilon)$-stationary point for $F$. 
\end{theorem} 

We obtain this result  by first reformulating \cref{prob:general-constraints} via the penalty method and constructing an
inexact gradient oracle for the hyperobjective $F$ (cf. \cref{sec:inequality-bilevel}). We then employ this inexact gradient oracle within {an algorithm (\cref{{alg: OIGRM}}) designed to minimize Lipschitz nonsmooth nonconvex functions}
(in particular, $F$), with the following convergence guarantee.
\begin{theorem}[Informal; \cref{{thm:Lipschitz-min-with-inexact-grad-oracle}}]
    Given Lipschitz $F:\reals^d\to\reals$ 
and $\|\widetilde{\nabla} F(\cdot)-\nabla F(\cdot)\|\leq\epsilon$,
 there exists an algorithm that, in 
$T=O(\delta^{-1}\epsilon^{-3})$
calls to $\widetilde{\nabla}F$,
outputs a $(\delta,2\epsilon)$-stationary point of $F$.
\end{theorem} 

Our  \cref{alg: OIGRM} is essentially a ``first-order'' version of  \cref{alg: IZO}. As was the case for \cref{alg: IZO}, despite the existence of algorithms with these guarantees with access to exact gradients~\cite{davis2022gradient}, their extensions to the \emph{inexact} gradient setting are not trivial and also make use of the new framework of \citet{cutkosky2023optimal}. 
Lastly, we also design a more implementation-friendly variant of \cref{alg: OIGRM} (with slightly worse theoretical guarantees) that we use in preliminary numerical experiments. 
\end{description}
\subsection{Related work}\label{sec:related_work} 
The vast body of work on asymptotic results for bilevel programming, starting with classical works such as \citet{anandalingam1990solution, ishizuka1992double, white1993penalty, vicente1994descent, ye1995optimality, ye1997exact}, typically fall into two categories: those based on approximate implicit differentiation: ~\citet{domke2012generic, pedregosa2016hyperparameter, gould2016differentiating, amos2017optnet, liao2018reviving, agrawal2019differentiable, grazzi2020iteration, lorraine2020optimizing} and those via iterative differentiation:~\citet{domke2012generic, maclaurin2015gradient, franceschi2017forward,   franceschi2018bilevel,  shaban2019truncated, grazzi2020iteration}. Another recent line of work  in this category includes \citet{liu2021value, ye2023difference, khanduri2023linearly, gao2024moreau}, which use various smoothing techniques. 

The first non-asymptotic result for bilevel programming was provided by \citet{ghadimi2018approximation}, which was followed by a flurry of work: for example, algorithms that are single-loop stochastic: ~\citet{chen2021closing, chen2022single, hong2023two}, algorithms that are projection-free: ~\citet{akhtar2022projection, jiang2023conditional, abolfazli2023inexact, cao2024projection}, those incorporating in them  variance-reduction and momentum: ~\citet{khanduri2021near, guo2021randomized, yang2021provably, dagreou2022framework}, those for single-variable bilevel programs: ~\citet{sabach2017first, amini2019iterative, amini2019iterativereg, jiang2023conditional, merchav2023convex}, and  for bilevel programs with special  constraints: ~\citet{tsaknakis2022implicit, khanduri2023linearly, xu2023efficient, abolfazli2023inexact}. 

The most direct predecessors of our work are those of \citet{liu2021value, kwon2023fully, yao2024constrained, lu2024firstorder, khanduri2023linearly}. As alluded to earlier, \citet{liu2022bome} recently made a significant contribution by providing  for bilevel programming a  first-order algorithm with finite-time stationarity guarantees. 
Their work was extended to the stochastic setting by \citet{kwon2023fully} (which we build upon for our results on the linear equality setting), simplified and improved by \citet{chen2023near}, and extended to the constrained setting by \citet{yao2024constrained, lu2024firstorder, khanduri2023linearly}.

The works of \citet{yao2024constrained, lu2024firstorder} study the more general problem of bilevel programming with {general convex} constraints. Their notion of stationarity is that of  KKT stationarity (over both upper and lower-level variables). The guarantee we provide in \cref{thm:linIneqInformal} holds for the much more  restricted setting of linear inequality constraints; however, our stationarity guarantees are expressed 
\textit{directly} in terms of the hyperobjective of interest, which we believe is more useful. Extending this notion of stationarity to bilevel programs with \emph{general} convex constraints (as studied in \citet{yao2024constrained, lu2024firstorder}) would be a very interesting direction of study. 
Moreover, the work of \citet{yao2024constrained} assumes joint convexity of the lower-level constraints in upper and lower variables to allow for efficient projections, 
whereas our guarantee requires convexity only in the lower-level variable. 

The current best result for the linearly constrained setting is that of \citet{khanduri2023linearly}. However, their algorithm is not  first-order method work 
due to requiring  Hessian computations and also imposes {some regularity assumptions}  on the hyperobjective $F$, which are, in general, impossible to verify.
In contrast, our \Cref{thm:lineq-cost-inf} imposes \textit{assumptions solely on the constituent functions $f$ and $g$} --- and  none directly on $F$ ---
thus making substantial progress on these two fronts.

\section{Preliminaries}\label{sec:preliminaries}
We use $\langle \cdot{}, {}\cdot\rangle$ to denote inner products and $\|{}\cdot{}\|$ for the Euclidean norm. Unless transposed, all vectors are column vectors. For $f:\reals^{d_2}\to\reals^{d_1}$ its Jacobian with respect to $x\in \reals^{d_2}$ is 
$\nabla f\in \reals^{d_1 \times d_2}$.  For $f:\reals^d\to\reals$, we overload  $\nabla f$ to refer to its gradient (the transposed Jacobian), a column vector. We use 
$\nabla_x$ to denote partial derivatives with respect to $x$.  

A function $f:\reals^n\to\reals^m$ is $L$-Lipschitz if for any $x,y$, we have $\|f(x) - f(y)\|\leq L \|x-y\|$.
A differentiable function $f:\reals^n\to\reals$ is convex if for any $x, y\in \reals^n$ we have $f(y)\geq f(x) + \nabla f(x)^\top (y-x)$; 
it is 
 $\mu$-strongly convex
 if $f - \tfrac{\mu}{2}\|{}\cdot{}\|^2$ is convex;
it is  $\beta$-smooth
if
$\nabla f$ is $\beta$-Lipschitz.

\begin{definition}\label{def:GoldsteinDeltaEpsStationary}
 Consider a locally Lipschitz function $f:\reals^d\to\reals$, a point $x\in\reals^d$, and a parameter $\delta>0$. The {Goldstein subdifferential}~\cite{goldstein1977optimization} of $f$ at $x$ is the set
 $\partial_{\delta} f(x):=\mathrm{conv} (\cup_{y\in \mathbb{B}_{\delta}(x)}\partial f(y)),$ where  $\partial f(x)=\mathrm{conv}\left\{\lim_{n\to \infty} \nabla f(x_n): x_n\to x, ~x_n\in \mathrm{dom}(\nabla f)\right\}$ is the Clarke subdifferential~\cite{clarke1981generalized} of $f$ and $\mathbb B_\delta(x)$ denotes the Euclidean ball of radius $\delta$ around $x$.
 A point $x$ is called $(\delta, \epsilon)$-stationary if $\mathrm{dist}(0,\partial_{\delta}f(x))\leq \epsilon$, where $\mathrm{dist}(x,S):=\inf_{y\in S}\|x-y\|$.
\end{definition}

Thus, for a Lipschitz function $f$, a point $x$ is $(\delta, \epsilon)$-stationary if within a $\delta$-ball around $x$, there exists a convex combination of  subgradients of $f$ with norm at most $\epsilon$. 
For a differentiable function $f$, we say that $x$ is $\epsilon$-stationary if $\|\nabla f(x)\|\leq \epsilon$. 

\subsection{Assumptions}\label{sec:differentiable-optimization}
We consider \cref{{prob:general-constraints}} with linear equality constraints (\cref{sec:equality-bilevel}) under \cref{assumption:linEq_smoothness,assumption:eq} and  linear inequality constraints (\cref{sec:inequality-bilevel,sec:nonsmooth}) under \cref{assumption:linEq_smoothness,assumption:ineq_mild,item:assumption_safe_constraints}. 
We assume the upper-level (UL) variable $x\in \reals^{d_x}$,  lower-level (LL) variable $y\in \reals^{d_y}$, and $A\in \reals^{d_h \times d_x}$.

\begin{assumption}\label{assumption:linEq_smoothness} 
For \cref{prob:general-constraints}, we assume the following for both settings we study: 
\begin{assumpenum}
\compresslist{
\item\label[assum]{item:assumption_upper_level} Upper-level: The objective $f$ is $C_f$-smooth and $L_f$-Lipschitz continuous in $(x,y)$. 
\item \label[assum]{item:assumption_lower_level}Lower-level: The objective $g$ is $C_g$-smooth. Fixing any $x\in\mathcal{X}$, $g(x,\cdot)$ is $\mu_g$-strongly convex.
\item \label[assum]{item:assumption_tangen_space}  We assume that the linear independence constraint qualification (LICQ) condition holds for the LL problem at every $x$ and $y$, i.e., the constraint $h(x,y)\coloneqq Ax-By-b$ has a full row rank $B$.
}
\end{assumpenum}
\end{assumption}

\begin{assumption}\label{assumption:eq}For \cref{prob:lin-eq} (with linear equality constraints), we additionally assume that the set $\mathcal{X}$ is convex and compact, and that the objective $g$ is $S_g$-Hessian smooth, that is, $\norm{\grad^2 g(x,y) -\grad^2 g(\bar x,\bar y)}\leq S_g \norm{(x,y) - (\bar x, \bar y)} \forall x,\bar x \in \mathcal{X}, \text{ and } y, \bar y \in \R^{d_y}.$ 
\end{assumption}

\begin{assumption}\label{assumption:ineq_mild}For \cref{prob:ineq} (with linear inequality constraints), we additionally assume  that $y^*$ is $L_y$-Lipschitz in $x$, where $y^*$ is the LL primal solution $y^*(x), \lambda^*(x) = \arg\max\nolimits_y \min\nolimits_{\beta \geq 0} g(x,y) + \beta^\top h(x,y)$, where $h(x,y)\coloneqq Ax-By-b$. 
\end{assumption}

\begin{assumption}\label{item:assumption_safe_constraints} 
We provide additional results for \cref{prob:ineq}  under additional stronger assumptions stated here: Denote the LL primal and dual solution $y^*(x), \lambda^*(x) = \arg\max\nolimits_y \min\nolimits_{\beta \geq 0} g(x,y) + \beta^\top h(x,y)$, where $h(x,y)\coloneqq Ax-By-b$; then, we  assume that we have exact access to $\lambda^*$ and that $\norm{\lambda^*(x)} \leq R$. 
\end{assumption}

\cref{item:assumption_upper_level,item:assumption_lower_level} are standard in bilevel optimization. 
\cref{{item:assumption_tangen_space}} is the same as the complete recourse assumption in stochastic programming \cite{shapiro2021lectures}, that is, the LL problem is feasible $y$ for every $x\in\R^{d_x}$. \cref{assumption:eq} is used only in the equality case and guarantees smoothness of $F$. \cref{assumption:ineq_mild} is used in the inequality case and implies Lipschitzness of  $F$. We need the stronger assumption detailed in  
\cref{item:assumption_safe_constraints}  for our dimension-free result for the linear inequality case. 

\section{Lower-level problem with linear equality constraint}\label{sec:equality-bilevel}
\looseness=-1We first obtain improved results for the  setting of bilevel programs with linear \textit{equality} constraints in the lower-level problem. Our formal problem statement is:   
\begin{align*}\numberthis\label[prob]{prob:lin-eq}
     & \mbox{minimize}_{x\in \mathcal{X}} ~ F(x) \coloneqq f(x, \ystar(x))~\text{ subject to } ~\ystar(x)\in \arg\min\nolimits_{y: h(x, y) = 0} g(x, y), 
\end{align*}
where $f$, $g$, $h(x,y)\coloneqq Ax-By-b$, and $\mathcal{X}$ satisfy 
\cref{assumption:linEq_smoothness,assumption:eq}.
\looseness=-1The previous best result on \cref{{prob:lin-eq}}  providing finite-time $\epsilon$-stationarity guarantees, by ~\citet{khanduri2023linearly}, required certain regularity assumptions on $F$ as well as Hessian computations. In contrast, our finite-time guarantees \emph{require assumptions only on $f$ and $g$, not on $F$};
indeed, in our work, these desirable properties of $F$ are naturally implied by our analysis. 
Specifically, our key insight is that the hypergradient $\grad F(x):=\grad_{x}f(x,y^{*})+\left(\tfrac{d y^*(x)}{dx}\right)^\top\grad_{y}f(x,y^{*})$ for  \cref{{prob:lin-eq}} is Lipschitz-continuous and  admits an easily computable --- yet highly accurate --- finite-difference approximation. Therefore, $O(\epsilon^{-2})$ iterations of gradient descent  on $F$ with this finite-difference gradient proxy  yield an $\epsilon$-stationary point. 

Specifically, for any fixed $x\in \mathcal{X}$, our proposed finite-difference gradient proxy approximating the non-trivial-to-compute component $\left(\tfrac{d y^*(x)}{dx}\right)^\top\grad_{y}f(x,y^{*})$ of the hypergradient is given
by 
 \[ v_x := {\frac{\nabla_{x}[g(x,\ystardel)+\langle\lamdeltar, h(x,\ystar)\rangle]-\nabla_{x}[g(x,y^{*})+\langle\lamstar, h(x,y^{*})\rangle]}{\delta}
 \numberthis\label{eq:part-hypergrad-approx-lin-eq},}\]
where $(\ydelstar,\lamdeltar)$ are the primal and dual solutions to the
 perturbed lower-level problem:
\begin{align}\label{eq:lower_perturb}
    & \ydelstar \coloneqq \arg\min\nolimits_{y: h(x,y) = 0} ~g(x,y)+\delta f(x,y).
\end{align} We show in \cref{{lem:lineq-finitediff-equals-gradf}} that $v$ in \cref{eq:part-hypergrad-approx-lin-eq} approximates $\left(\tfrac{d y^*(x)}{dx}\right)^\top\grad_{y}f(x,y^{*})$ up to an $O(\delta)$-additive error, implying the gradient oracle construction outlined in the pseudocode presented in 
\cref{alg:LE-inexact-gradient-oracle}. Our full implementable algorithm for solving \cref{prob:lin-eq} is displayed in \cref{{alg:LE-full-alg}}.

\begin{algorithm}[h]\caption{Inexact Gradient Oracle for Bilevel Program with Linear Equality Constraint}\label{alg:LE-inexact-gradient-oracle}
\begin{algorithmic}[1]
\State \textbf{Input:}
Current $x$, accuracy $\epsilon$, perturbation $\delta = \epsilon^2$.
\State Compute $y^*$ (as in \cref{prob:lin-eq}) and corresponding optimal dual $\lamstar$ (as in \cref{eqs:kkt-lin-eq})\label{line:LEQ-first}
\State Compute $\ystardel$ (as in \cref{eq:lower_perturb}) and and corresponding optimal dual $\lamdeltar$ (as in \cref{eqs:lineq-kkt-perturbed}) \label{line:LEQ-second} 
\State Compute $v_x$ as in \cref{eq:part-hypergrad-approx-lin-eq} \Comment{Approximates $\left({d y^{*}(x)}/{d x}\right)^\top\nabla_{y}f(x,y^{*})$}
\State \textbf{Output:} $\widetilde{\nabla} F = v_x + \nabla_x f(x, y^*)$
\end{algorithmic}
\end{algorithm}

Notice that the finite-difference term in \cref{{eq:part-hypergrad-approx-lin-eq}} avoids differentiating through the implicit function $y^*$. Instead, all we need to evaluate it are the values of $(y^*,\lambda^*,\ydelstar,\lamdeltar)$ (and gradients of $g$ and $h$). Since $(y^*,\lambda^*)$ are solutions to a smooth strongly convex linearly constrained  problem, they can be approximated
at a linear rate. Similarly, since the approximation error in \cref{eq:part-hypergrad-approx-lin-eq}  is proportional to $\delta$ (cf. \cref{{lem:lineq-finitediff-equals-gradf}}), a  small enough $\delta$ in the perturbed objective $ g+\delta f$ in \cref{eq:lower_perturb} ensures that it is dominated by the strongly convex and smooth $g$, whereby accurate approximates to $(y_\delta^*,\lambda_\delta^*)$ can also be readily obtained. Putting it all together, the proposed finite-difference hypergradient proxy in \cref{eq:part-hypergrad-approx-lin-eq} is  efficiently computable, yielding the following guarantee. 
\begin{restatable}{theorem}{linEqFullCost}\label{thm:lineq-cost}
    Consider  \cref{prob:lin-eq} under \cref{assumption:linEq_smoothness}, and let $\kappa=C_g/\mu_g$ be the condition number of $g$. Then \cref{alg:LE-full-alg}  finds an $\epsilon$-stationary point (in terms of gradient mapping, see \cref{eq:gradient-mapping}) after $T=\widetilde{O}(C_F (F(x_0)-\inf F)\sqrt{\kappa}\epsilon^{-2})$  oracle calls to $f$ and $g$, where $C_F:= 2(L_f +C_f+C_g)C_H^3 S_g (L_g +\|A\|)^2$ is the smoothness constant of the hyperobjective $F$.
\end{restatable}  
We now sketch the proofs of the key components that together imply \cref{thm:lineq-cost}. The complete proofs may be found in \cref{{sec:appendix_linear_equality}}.   
\subsection{Main technical ideas}\label{sec:warmup-lin-eq} 
We briefly outline the two key technical building blocks alluded to above, that together give us  \cref{{thm:lineq-cost}}: the approximation guarantee of our finite-difference gradient proxy (\cref{{eq:part-hypergrad-approx-lin-eq}}) and the smoothness of hyperobjective $F$ (for \cref{{prob:lin-eq}}). The starting point for both these results is the following simple observation obtained by implicitly differentiating, with respect to $x$, the KKT system associated with $y^*=\arg\min_{h(x,y)=0} g(x,y)$ and optimal dual variable $\lamstar$: 
 \[ \begin{bmatrix}\frac{d\ystar(x)}{dx}\\
\frac{d\lamstar(x)}{dx}
\end{bmatrix}= {\begin{bmatrix}\grad^2_{yy}g(x,\ystar) & \nabla_y h(x,y^{*})^{\top}\\
\nabla_y h(x,y^{*}) & 0
\end{bmatrix}^{-1}} \begin{bmatrix}-\nabla^2_{yx}g(x,y^{*})\\
-\nabla_x h(x,\ystar)\end{bmatrix} \numberthis\label{eq:def_matrix_H}\]
The invertibility of the matrix in the preceding equation is proved in  \cref{{cor:nonsingularH}}:
essentially, this invertibility is implied by strong convexity of $g$ and $\grad_y h(x,y^*)=B$ having full row rank. This in conjunction with the compactness of $\mathcal{X}$ implies that the inverse of the matrix is bounded by some constant $C_H$ (cf. \cref{cor:nonsingularH} for details). Our hypergradient approximation guarantee follows: 

\begin{restatable}{lemma}{lemLineqFiniteDiffEqualsGradF}\label{lem:lineq-finitediff-equals-gradf} 
For  \cref{prob:lin-eq}
under \cref{assumption:linEq_smoothness}, with $v_x$ as in \cref{eq:part-hypergrad-approx-lin-eq}, 
the
following holds:
\[
\textstyle{\left\|v_x-\left(\frac{d y^{*}(x)}{d x}\right)^\top\nabla_{y}f(x,y^{*})\right\|\leq O(C_F \delta).}
\]
\end{restatable}
\begin{proof}[Proof sketch; see \cref{{sec:appendix_linear_equality}} for the complete proof]
    The main idea for this proof is that the two terms being compared are essentially the same by an application of the implicit function theorem. First, we can use the expression for $\tfrac{dy^*(x)}{dx}$ from \cref{eq:def_matrix_H} as follows: 
\begin{align*}
\left(\frac{d y^{*}(x)}{d x}\right)^\top\nabla_{y}f(x,y^{*})
&=\begin{bmatrix}-\nabla^2_{yx} g(x,y^{*})\\
-\nabla_x h(x,\ystar)
\end{bmatrix}^\top\begin{bmatrix}\nabla^2_{yy} g(x,y^{*}) & \nabla_y h(x,y^{*})^{\top}\\
\nabla_y h(x,y^{*}) & 0
\end{bmatrix}^{-1}\begin{bmatrix}\nabla_{y}f(x,y^{*}) \\ 0\end{bmatrix}.
\end{align*} We now examine our finite-difference gradient proxy. For simplicity of exposition, we instead consider $\lim_{\delta\rightarrow0}\frac{\nabla_{x}[g(x,\ystardel)+\langle\lamdeltar, h(x,y^*)\rangle]-\nabla_{x}[g(x,y^{*})+\langle\lamstar, h(x,y^{*})\rangle]}{\delta},$ which, by an application of the fundamental theorem of calculus and \cref{assumption:eq}
, equals our finite-difference gradient proxy up to an $O(\delta)$-additive error. Note that this expression  is exactly \[\nabla^2_{xy} g(x,y^*) \frac{d\ydelstar}{d\delta} + \nabla_x h(x,y^*)^\top \frac{d\lamdeltar}{d\delta}.\numberthis\label{eq:LINEQ-lim-findiff}\] Since  $y^*_{\delta}$ is also a minimizer of a strongly convex function over a linear equality constraint \cref{{eq:lower_perturb}}, the same reasoning that yields \cref{eq:def_matrix_H} can be also used to obtain 
\[
\begin{bmatrix}\frac{d\ydelstar(x)}{d\delta}\\
\frac{d\lamdeltar(x)}{d\delta}
\end{bmatrix}\Bigg|_{\delta=0}=\begin{bmatrix}\nabla^2_{yy} g(x,\ystar) & \nabla_y h(x,\ystar)^{\top}\\
\nabla_y h(x,\ystar) & 0
\end{bmatrix}^{-1}\begin{bmatrix}-\nabla_y f(x,\ystar)\\
0
\end{bmatrix}. \numberthis\label{eq:LINEQ-dystardelta}
\] Combining \cref{eq:LINEQ-lim-findiff} and \cref{eq:LINEQ-dystardelta}  immediately yields the claim. 
\end{proof}

Having shown our hypergradient approximation, we now turn to showing smoothness of the hyperobjective $F$, crucial to getting our claimed rate. 

\begin{restatable}{lemma}{lemYdelstarLamdelstarSmooth}\label{lem:smoothness_of_ydelstar_lamdelstar}
The solution $y^*$ (as defined in \cref{{prob:lin-eq}}) is $O(C_H \cdot (\gssmooth + \|A\|))$-Lipschitz continuous and $O(C_H^3\cdot\gtsmooth\cdot(\gssmooth+\|A\|)^2)$-smooth as a function
of $x$. Thus the hyper-objective $F$ is gradient-Lipschitz  with a  smoothness constant of $C_F:=O\{ (L_f +C_f+C_g)C_H^3 S_g (L_g +\|A\|)^2\}.$ 
\end{restatable}
\begin{proof}[Proof sketch; see \cref{{sec:appendix_linear_equality}}]
 The Lipschitz bound follows by combining \cref{{eq:def_matrix_H}} with the smoothness of $g$ and bound $C_H$ on the matrix inverse therein. Differentiating \cref{{eq:def_matrix_H}} by $x$ yields a linear system with the same matrix;  we repeat this approach to get the  smoothness bound. 
\end{proof}

\section{Lower-level problem with linear inequality constraint: nonsmooth nonconvex optimization with inexact  oracles}\label{sec:nonsmooth}
We now shift gears from the case of linear \emph{equality} constraints to that of linear \emph{inequality} constraints. Specifically, defining $h(x,y)=Ax-By-b$,  the problem we now consider is 
\begin{align}\label[prob]{prob:ineq}
     \mathop{\text{minimize}}\nolimits_{x} ~ F(x) \coloneqq f(x, \ystar(x)) \quad \text{ subject to} ~\ystar(x)\in \argmin\nolimits_{y: h(x, y) \leq 0} g(x, y). 
\end{align}
As noted earlier, for this larger class of problems, the hyperobjective $F$ can be nonsmooth nonconvex, thereby necessitating  our measure of convergence to be the now popular notion of Goldstein stationarity~\cite{zhang2020complexity}. 

Our first algorithm for solving \cref{prob:ineq}  is  presented in \cref{alg: IZO}, with its convergence guarantee in \cref{thm:izo_complete_guarantees}.
At a high level, this algorithm first obtains access to an inexact zeroth-order oracle to $F$ (we shortly explain how this is done) and uses this oracle to construct a (biased) gradient estimate of $F$. It then uses this gradient estimate to update the iterates with an update rule motivated by  recent works reducing nonconvex optimization to online optimization~\cite{cutkosky2023optimal}. We explain these components in \cref{sec:ncns_inexact_ZO}.

\begin{algorithm}[h]
\begin{algorithmic}[1]\caption{Nonsmooth Nonconvex Algorithm with Inexact Zero-Order oracle}\label{alg: IZO}
\State \textbf{Input:}
Initialization $x_0\in\reals^d$, clipping parameter $D>0$,
step size $\eta>0$, smoothing parameter $\rho>0$, accuracy parameter $\nu>0$,
iteration budget $T\in\NN$, inexact zero-order oracle $\tF:\reals^d\to\reals$.
\State \textbf{Initialize:} ${\Delta}_1=\mathbf{0}$
\For{$t=1,\dots,T$}
\State Sample $s_t\sim\Unif[0,1]$,~~$w_t\sim\Unif(\S^{d-1})$
\State $x_t=x_{t-1}+{\Delta}_t$, ~~$z_t=x_{t-1}+s_t{\Delta}_t$
\State $\tbg_t=\tfrac{d}{2\rho}(\tF(z_t+\rho w_t)-\tF( z_t-\rho w_t))w_t$
\State ${\Delta}_{t+1}
=\mathrm{clip}_{D}\left({\Delta}_t-\eta\tbg_t\right) 
$
\Comment{$\mathrm{clip}_D(z):=\min\{1,\tfrac{D}{\norm{z}}\}\cdot z$}
\EndFor
\State $M=\lfloor\frac{\nu}{D}\rfloor,~K=\lfloor\frac{T}{M}\rfloor$
\For{$k=1,\dots,K$}
\State $\overline{x}_{k}=\frac{1}{M}\sum_{m=1}^{M} z_{(k-1)M+m}$
\EndFor
\State Sample $ x^{\out}\sim\Unif\{\overline{ x}_1,\dots,\overline{ x}_{K}\}$
\State \textbf{Output:} $ x^{\out}$. 
\end{algorithmic}
\end{algorithm}

\begin{theorem}\label{thm:izo_complete_guarantees}
    Consider  \cref{{prob:ineq}} under \cref{assumption:linEq_smoothness} and \cref{assumption:ineq_mild}.  Let $\kappa=C_g/\mu_g$ be the condition number of $g$. Then combining  the procedure for \cref{{lem:ZeroOrderApprox}}
with \cref{alg: IZO} run with parameters 
$\rho=\min\left\{\tfrac{\delta}{2},\tfrac{F(x_0)-\inf F}{L_fL_y}\right\},\nu=\delta-\rho,~D=\Theta\left(\frac{\nu\epsilon^2\rho^2}{d\rho^2 L_f^2L_y^2+\alpha^2 d^2}\right)$, and $\eta=\Theta\left(\frac{\nu\epsilon^3\rho^4}{(d\rho^2L_f^2L_y^2+\alpha^2d^2)^2}\right)$ 
outputs a point $x^{\out}$ such that $\E[\mathrm{dist}(0,{\partial}_\delta F(x^{\out}))]\leq\epsilon+\alpha$ with
\[T=
O\left(\frac{\sqrt{\kappa}d(F(x_0)-\inf F)}{\delta\epsilon^3}\cdot \left(L_f^2L_y^2+\alpha^2 \left(\frac{d}{\delta^{2}}+\frac{dL_f^2L_y^2}{(F(x_0)-\inf F)^2}\right)\right)\cdot\log(L_f/\alpha)\right)
\text{ oracle calls to } f \text{ and } g.\] 
\end{theorem}

\subsection{Nonsmooth nonconvex optimization with inexact zeroth-order oracle}\label{sec:ncns_inexact_ZO}
We can obtain inexact zeroth-order oracle access to $F$  because (as formalized in \cref{lem:ZeroOrderApprox}) despite potential nonsmoothness and nonconvexity of the hyperobjective $F$ in \cref{prob:ineq},  estimating its \emph{value} $F(x)$  at any given point $x$ amounts to solving a single smooth and strongly convex optimization problem, which can be done can be done using $\widetilde{O}(1)$ oracle calls to $f$ and $g$ by appealing to a result by \citet{zhang2022solving}.

\begin{lemma}\label{lem:ZeroOrderApprox}
    Given any $x$, we can return $\widetilde{F}(x)$ such that $|F(x)-\widetilde{F}(x)|\leq\alpha$ using $O(\sqrt{C_g/\mu_g}\log(L_f/\alpha))$ first-order oracle calls to $f$ and $g$.
\end{lemma}

\begin{proof}[Proof of \cref{lem:ZeroOrderApprox}]
    We note that it suffices to find $\tilde{y}^*$ such that $\norm{\tilde{y}^*-y^*(x)}\leq \alpha/L_f$, since setting $\tF(x):=f(x,\tilde{y}^*)$ will then satisfy $|\tF(x)-F(x)|=|f(x,\tilde{y}^*)-f(x,y^*(x))|\leq L_f\cdot \tfrac{\alpha}{L_f}=\alpha$ by Lispchitzness of $f$, as required. Noting that $y^*(x)=\arg\min_{h(x,y)\leq 0}g(x,y)$ is the solution to a constrained smooth, strongly-convex problem with condition number $C_g/\mu_g$, it is possible to approximate it up to $\alpha/L_f$ with $O(\sqrt{C_g/\mu_g}\log(L_f/\alpha))$ first-order oracle calls using the result of \citet{zhang2022solving}.
\end{proof}

Having computed an inexact value of the  hyperobjective $F$ (in \cref{thm:cost_of_computing_ystar_gammastar_inequality}), we now show how to use it to develop an algorithm for solving \cref{prob:general-constraints}. To this end, we first note that  $F$, despite being 
 possibly nonsmooth and nonconvex, is Lipschitz and hence amenable to the use of recent algorithmic developments  in nonsmooth nonconvex optimization pertaining to Goldstein stationarity. 
\begin{restatable}{lemma}{lemLipscConstrBilevel}\label{lem:LipscConstrBilevel}Under \cref{assumption:linEq_smoothness} and \ref{item:assumption_safe_constraints}, $F$ in \cref{prob:general-constraints} is $O(L_fL_y)$-Lipschitz in $x$.
\end{restatable}

\looseness=-1With this guarantee on the Lipschitzness of $F$, we prove \cref{{thm: Lipschitz-min-with-inexact-zero-oracle}} for attaining Goldstein stationarity using the inexact zeroth-order oracle of a Lipschitz function. Our proof of \cref{{thm: Lipschitz-min-with-inexact-zero-oracle}} crucially uses the recent online-to-nonconvex framework of \citet{cutkosky2023optimal}. Combining \cref{lem:ZeroOrderApprox} and \cref{thm: Lipschitz-min-with-inexact-zero-oracle} then immediately implies \cref{thm:izo_complete_guarantees}. 

\begin{theorem}\label{thm: Lipschitz-min-with-inexact-zero-oracle}
Suppose $F:\reals^d\to\reals$ is $L$-Lipschitz, 
and that $|\widetilde{F}(\cdot)-F(\cdot)|\leq\alpha$.
Then running \cref{alg: IZO} with
$\rho=\min\left\{\tfrac{\delta}{2},\tfrac{F(x_0)-\inf F}{L}\right\},\nu=\delta-\rho,~D=\Theta\left(\frac{\nu\epsilon^2\rho^2}{d\rho^2 L^2+\alpha^2 d^2}\right),\eta=\Theta\left(\frac{\nu\epsilon^3\rho^4}{(d\rho^2L^2+\alpha^2d^2)^2}\right)$,
outputs a point $x^{\out}$ such that $\E[\mathrm{dist}(0,{\partial}_\delta F(x^{\out}))]\leq\epsilon+\alpha$ with
\[T=
O\left(\frac{d(F(x_0)-\inf F)}{\delta\epsilon^3}\cdot \left(L^2+\alpha^2 (\frac{d}{\delta^{2}}+\frac{dL^2}{(F(x_0)-\inf F)^2})\right)\right)
\text{ calls to } \tF(\cdot).\] 
\end{theorem}

\subsection{Nonsmooth nonconvex optimization with inexact \emph{gradient}  oracle}
In the next section (\cref{{sec:inequality-bilevel}}), we provide a way to generate approximate gradients of $F$. Here, we present an algorithm that attains Goldstein stationarity of \cref{prob:ineq} using this inexact gradient oracle.  While there has been a long line of recent work on algorithms for nonsmooth nonconvex optimization with convergence to Goldstein stationarity~\cite{zhang2020complexity, davis2022gradient, jordan2023deterministic, kong2023cost, grimmer2023goldstein}, these results necessarily require \textit{exact} gradients. 
This brittleness to any error in gradients renders them ineffective in our setting, where our computed (hyper)gradient necessarily suffers from an
additive error.
While inexact oracles are known to be effective for smooth or convex objectives \citep{devolder2014first}, utilizing inexact gradients in the nonsmooth nonconvex regime presents a nontrivial challenge.
Indeed, without any local bound on gradient variation due to smoothness, or convexity that ensures that gradients are everywhere correlated with the direction pointing at the minimum, it is not clear a priori how to control the accumulating price of inexactness throughout the run of an algorithm.
To derive such results, we again use the recently proposed connection between online learning and nonsmooth nonconvex optimization by \citet{cutkosky2023optimal}. In particular, by controlling the accumulated error suffered by online gradient descent for \emph{linear} losses (cf. \cref{lem: inexact OGD}), 
we derive guarantees for our setting of interest, providing Lipschitz optimization algorithms that converge to Goldstein stationary points even with inexact gradients.

This algorithm matches the best known complexity in first-order nonsmooth nonconvex optimization \cite{zhang2020complexity,davis2022gradient,cutkosky2023optimal}, merely replacing the convergence to a $(\delta,\epsilon)$-stationary point by $(\delta,\epsilon+\alpha)$-stationarity, where $\alpha$ is the inexactness of the gradient oracle.

\begin{algorithm}[h]
\begin{algorithmic}[1]\caption{Nonsmooth Nonconvex Algorithm with Inexact Gradient Oracle}\label{alg: OIGRM}
\State \textbf{Input:}
Initialization $x_0\in\reals^d$, clipping parameter $D>0$,
step size $\eta>0$,
accuracy parameter $\delta>0$,
iteration budget $T\in\NN$, inexact gradient oracle $\widetilde{\nabla}F:\reals^d\to\reals^d$.
\State \textbf{Initialize:} ${\Delta}_1=\mathbf{0}$
\For{$t=1,\dots,T$}
\State Sample $s_t\sim\Unif[0,1]$ 
\State $x_t=x_{t-1}+{\Delta}_t$, ~~$z_t=x_{t-1}+s_t{\Delta}_t$
\State $\tbg_t=\widetilde{\nabla}F(z_t)$
\State ${\Delta}_{t+1}
=\mathrm{clip}_{D}\left({\Delta}_t-\eta\tbg_t\right) 
$
\Comment{$\mathrm{clip}_D(z):=\min\{1,\tfrac{D}{\norm{z}}\}\cdot z$}
\EndFor
\State $M=\lfloor\frac{\delta}{D}\rfloor,~K=\lfloor\frac{T}{M}\rfloor$
\For{$k=1,\dots,K$}
\State $\overline{x}_{k}=\frac{1}{M}\sum_{m=1}^{M} z_{(k-1)M+m}$
\EndFor
\State Sample $ x^{\out}\sim\Unif\{\overline{ x}_1,\dots,\overline{ x}_{K}\}$
\State \textbf{Output:} $ x^{\out}$. 
\end{algorithmic}
\end{algorithm}

\begin{restatable}{theorem}{thmLipscMinWithInexactGradOracle}\label{thm:Lipschitz-min-with-inexact-grad-oracle}
Suppose $F:\reals^d\to\reals$ is $L$-Lipschitz  
and that $\|\widetilde{\nabla} F(\cdot)-\nabla F(\cdot)\|\leq\alpha$. 
Then running \cref{alg: OIGRM} with
$D=\Theta(\frac{\delta\epsilon^2}{L^2}),\eta=\Theta(\frac{\delta\epsilon^3}{L^4})$,
outputs a point $x^{\out}$ such that $\E[\mathrm{dist}(0,{\partial}_\delta F(x^{\out}))]\leq\epsilon+\alpha$, with
$T=O\left(\frac{(F(x_0)-\inf F) L^2}{\delta\epsilon^3}\right)$
calls to $\widetilde{\nabla}F(\cdot)$.
\end{restatable}

We defer the proof of \cref{thm:Lipschitz-min-with-inexact-grad-oracle} to \cref{{sec:zeroth-order-algs}}.
Plugging the complexity of computing inexact gradients,
as given by \cref{thm:cost_of_computing_ystar_gammastar_inequality},
into the result above,
we immediately obtain convergence to a $(\delta,\epsilon)$-stationary point of \cref{prob:general-constraints} with $\widetilde{O}(\delta^{-1}\epsilon^{-4})$ gradient calls overall.

\paragraph{Implementation-friendly algorithm.}
While \cref{{alg: OIGRM}} matches the best-known results in nonsmooth nonconvex optimization, it could be  impractical  due to several hyperparameters which need tuning to guarantee optimal results. 
Arguably, a more natural application of the hypergradient estimates would be simply plugging them into gradient descent, which requires tuning only the stepsize. Since the objective $F$ is neither smooth nor convex, perturbations are required to guarantee differentiability along the trajectory.
We therefore complement \cref{thm:Lipschitz-min-with-inexact-grad-oracle} by analyzing perturbed (inexact) gradient descent  in the nonsmooth nonconvex setting (\cref{{alg: PIGD}}). Despite its suboptimal worst-case theoretical guarantees, we find this algorithm  easier to implement in practice. We defer the proof of \cref{{thm:practical_Lipschitz-min-with-inexact-grad-oracle}} to \cref{{sec:zeroth-order-algs}}.

\begin{algorithm}[h]
\begin{algorithmic}[1]\caption{Perturbed Inexact GD}\label{alg: PIGD}
\State \textbf{Input:}
Inexact gradient oracle $\widetilde{\nabla}F:\reals^d\to\reals^d$, initialization $x_0\in\reals^d$, spatial parameter $\delta>0$, step size $\eta>0$, iteration budget $T\in\NN$.
\For{$t=0,\dots,T-1$}
\State Sample $w_t\sim\Unif(\S^{d-1})$
\State $\tbg_t=\widetilde{\nabla}F(x_{t}+\delta\cdot w_t)$
\State $x_{t+1}
=x_t-\eta\tbg_t
$
\EndFor
\State \textbf{Output:} $x^{\out}\sim\mathrm{Unif}\{x_0,\dots,x_{T-1}\}$. 
\end{algorithmic}
\end{algorithm}

\begin{theorem}\label{thm:practical_Lipschitz-min-with-inexact-grad-oracle}
Suppose $F:\reals^d\to\reals$ is $L$-Lipschitz, 
and that $\|\widetilde{\nabla} F(\cdot)-\nabla F(\cdot)\|\leq\alpha$. 
Then running \cref{alg: PIGD} with
$\eta=\Theta\left(\frac{{((F(x_0)-\inf F)+\delta L)^{1/2}\delta^{1/2}}}{{T^{1/2} L^{1/2}{d}^{1/4}(\alpha+L)}}\right)$
outputs a point $x^{\out}$ such that $\E[\mathrm{dist}(0,{\partial}_\delta F(x^{\out}))]\leq\epsilon+\sqrt{\alpha L}$, with \[T=O\left(\frac{(F(x_0)-\inf F+\delta L) L^3 \sqrt{d}}{\delta\epsilon^4}\right) \text{ 
calls to } \widetilde{\nabla}F(\cdot).\] 

\end{theorem}

\section{Lower-level problem with linear inequality constraint: constructing the inexact gradient oracle}\label{sec:inequality-bilevel}
In this section, we present a method to approximate the hypergradient $\nabla F(x)$ at any given point $x$, where: 
\begin{align}\label{eqn:second-order-method}
    \nabla F(x) = \nabla_x f(x,y^*) + \left(\tfrac{ dy^*(x)}{d x}\right)^\top \nabla_y f(x,y^*).
\end{align}
The key computational challenge here lies in computing  $\tfrac{d y^*(x)}{d x}$ since it  requires differentiating through an argmin operator, 
which typically requires second-order derivatives.  
Instead, here we 
differentiate (using the implicit function theorem) through the KKT conditions describing $y^*(x)$
and get: 
\begin{align}\label{eqn:kkt-system}
\begin{bmatrix}
\nabla^2_{yy} g + (\lambda^*)^\top \nabla_{yy}^2 h & \nabla_y h_\mathcal{I}^\top \\
\text{diag}(\lambda^*_\mathcal{I}) \nabla_y h_\mathcal{I} & 0
\end{bmatrix}
\begin{bmatrix}
    \frac{d y^*(x)}{d x} \\
    \frac{d \lambda^*_\mathcal{I}(x)}{d x}
\end{bmatrix}
= 
-
\begin{bmatrix}
    \nabla^2_{yx} g + (\lambda^*)^\top \nabla_{yx}^2 h \\
    \text{diag}(\lambda^*_\mathcal{I}) \nabla_x h_\mathcal{I}
\end{bmatrix}
\end{align}
where given $x$, we assume efficient access to the optimal dual solution $\lambda^*(x)$ of the LL problem in \cref{prob:ineq}. In \cref{eqn:kkt-system}, we use $\mathcal{I} \coloneqq \{i \in [d_h]: h_i(x,y) = 0, \lambda^*_i > 0 \}$ to denote the set of active constraints with non-zero dual solution, with
$h_\mathcal{I} \coloneqq [h_i]_{i \in \mathcal{I}}$ and $\lambda^*_{\mathcal{I}} \coloneqq [\lambda^*_i]_{i \in \mathcal{I}}$ being the constraints and dual variables, respectively,  corresponding to $\mathcal{I}$.  

Observe that as currently stated, \cref{eqn:kkt-system} leads to a second-order computation of ${dy^*(x)}/{dx}$.
In the rest of the section, we provide a fully first-order {approximate} hypergradient oracle by constructing an equivalent reformulation of \cref{prob:ineq} using a penalty function.

\subsection{Reformulation via the penalty method}
We begin by reformulating \cref{prob:ineq} into a single level constrained optimization problem:
\begin{align}\label{eqn:inequality_reformulation}
    \text{minimize}_{x,y} ~ f(x,y)  \text{ subject to}
    ~\begin{cases}
        g(x,y) + (\lambda^*(x))^\top h(x,y) \leq g^*(x) \\ 
        h(x,y) \leq 0
    \end{cases},
\end{align}
where $g^*(x)  \coloneqq \min_{y: h(x,y) \leq 0} g(x,y) = g(x,y^*(x))$ and $\lambda^*(x)$ is the optimal dual solution. The equivalence of this reformulation to  \cref{prob:ineq} is spelled out in \cref{appendix:reformulation-equivalence}. 
From \cref{eqn:inequality_reformulation}, we  define the following penalty function, crucial to our analysis:
\begin{align}
    \mathcal{L}_{\lambda^*, \boldsymbol{\alpha}}(x,y) = f(x,y) + \alpha_1 \left( g(x,y) + (\lambda^*)^\top h(x,y) - g^*(x)  \right) + \frac{\alpha_2}{2} \norm{h_\mathcal{I}(x,y)}^2, \numberthis\label{eqn:penalty-lagrangian}
\end{align}
where $\boldsymbol{\alpha} = [\alpha_1, \alpha_2] \geq 0$ are the penalty parameters.
Notably, we can compute its derivative with respect to $x$ of \cref{eqn:penalty-lagrangian} in a fully first-order fashion by the following expression: 
\begin{align*}
    \nabla_x \mathcal{L}_{\lambda^*, \boldsymbol{\alpha}}(x, y) = \nabla_x f(x,y) \! + \! \alpha_1 (\nabla_x g(x,y) \! + \! \nabla_x h(x,y)^\top \lambda^* \! - \! \nabla_x g^*(x) ) \! + \! \alpha_2 \nabla_x h_\mathcal{I}(x,y)^\top h_\mathcal{I}(x,y).
\end{align*}
To give some intuition for our choice of penalties in \cref{{eqn:penalty-lagrangian}}, we note that  the two constraints in \cref{eqn:inequality_reformulation} behave quite differently.
The first constraint $g(x,y) + \lambda^*(x)^\top h(x,y)  \leq g^*(x)$ is one-sided, i.e., can only be violated or met, and hence just needs a penalty parameter $\alpha_1$ to weight the ``violation''.
As to the second constraint $h(x,y) \leq 0$, it can be arbitrary. To allow for such a ``two-sided'' constraint, we penalize only the active constraints $\mathcal{I}$, i.e., we use $\norm{h_{\mathcal{I}}(x,y)}^2$ to penalize deviation.

Next, we define the optimal solutions to the penalty function optimization by:
\begin{align*}
    y_{\lambda^*, \boldsymbol{\alpha}}^*(x) := \arg\min\nolimits_{y} \mathcal{L}_{\lambda^*, \boldsymbol{\alpha}}(x,y).\numberthis\label{eq:def_y_lambda_star} 
\end{align*}
We now show that this minimizer is close to the optimal solution of the LL problem, while suffering only a small constraint violation. 
\begin{restatable}[]{lemma}{solutionApproximation}\label{thm:solution-bound}
Given any $x$, the corresponding dual solution $\lambda^*(x)$,  primal solution $y^*(x)$ of the lower optimization problem in \cref{prob:ineq}, and $y_{\lambda^*, \boldsymbol{\alpha}}^*(x)$ as in \cref{{eq:def_y_lambda_star}}, satisfy:
\begin{align}\label{eqn:solution-bound}
    \norm{y_{\lambda^*, \boldsymbol{\alpha}}^*(x) - y^*(x)} \leq O(\alpha_1^{-1}) \text{~~and~~} \norm{h_\mathcal{I}(x,y_{\lambda^*, \boldsymbol{\alpha}}^*(x))} \leq O(\alpha_1^{-1/2} \alpha_2^{-1/2}).
\end{align}
\end{restatable}
The proof of \cref{thm:solution-bound} is based on the Lipschitzness of $f$ and strong convexity of $g$ for sufficiently large $\alpha_1$.
The aforementioned constraint violation bound on $h_\mathcal{I}(x,y)$ is later used in \cref{thm:diff_in_hypergrad_and_gradLagr} to bound the inexactness of our proposed gradient oracle.

\subsection{Main result: approximating the hypergradient}
The main export of this section is the following bound on the approximation of  the hypergradient. This, together with the bounds in 
\cref{thm:solution-bound}, validate our use of the penalty function in \cref{eqn:penalty-lagrangian}. 
\begin{restatable}[]{lemma}{gradientApproximation}\label{thm:diff_in_hypergrad_and_gradLagr}
Consider $F$ as in \cref{{prob:ineq}},   $\mathcal{L}$ as in  \cref{{eqn:penalty-lagrangian}},  a fixed $x$, and $y_{\lambda^*, \boldsymbol{\alpha}}^*$ as in \cref{{eq:def_y_lambda_star}}. Then under \cref{assumption:linEq_smoothness,item:assumption_safe_constraints}, we have: 
\begin{align*}
   & \norm{\nabla F(x) - \nabla_x \mathcal{L}_{\lambda^*, \boldsymbol{\alpha}}(x, y_{\lambda^*, \boldsymbol{\alpha}}^*)} \leq O({\alpha_1^{-1}}) + O({\alpha_1^{-1/2}\alpha_2^{-1/2}}) + O({\alpha_1^{1/2}}{\alpha_2^{-1/2}}) + O({\alpha_1^{-3/2}}{\alpha_2^{1/2}}).
\end{align*}
\end{restatable}
The proof can be found in \cref{appendix:proof-of-inexact-gradient}. 
With this hypergradient approximation guarantee, we design \cref{alg:inexact-gradient-oracle} to compute an inexact gradient oracle for the hyperobjective $F$.
\begin{algorithm}[h]\caption{Inexact Gradient Oracle for General Inequality Constraints}\label{alg:inexact-gradient-oracle}
\begin{algorithmic}[1]
\State \textbf{Input:}
Upper level variable $x$, accuracy $\alpha$, penalty parameters $\alpha_1 = {\alpha^{-2}}, \alpha_2 = {\alpha^{-4}}$.
\State Compute $y^*$, $\lambda^*$, and active constraints $\mathcal{I}$ of the constrained LL problem $\min_{y: h(x,y) \leq 0} g(x,y)$. \label{line:inner-optimization-problem} 
\State Define penalty function $\mathcal{L}_{\lambda^*, \boldsymbol{\alpha}}(x,y)$ by ~\cref{eqn:penalty-lagrangian}
\State Compute the minimizer $y^*_{\lambda^*, \boldsymbol{\alpha}} = \arg\min\nolimits_{y} \mathcal{L}_{\lambda^*, \boldsymbol{\alpha}}(x,y)$ (as in ~\cref{eq:def_y_lambda_star}).\label{line:lagrangian-optimization}
\State \textbf{Output:} $\widetilde{\nabla} F \coloneqq \nabla_x \mathcal{L}_{\lambda^*, \boldsymbol{\alpha}}(x,y^*_{\lambda^*, \boldsymbol{\alpha}})$.
\end{algorithmic}
\end{algorithm}
\begin{restatable}[]{theorem}{computationCostInequality}\label{thm:cost_of_computing_ystar_gammastar_inequality}
    Given any accuracy parameter $\alpha > 0$, \cref{alg:inexact-gradient-oracle} outputs $\widetilde{\nabla}_x F(x)$ such that $\|\widetilde{\nabla} F(x) - \nabla F(x)\| \leq \alpha$ within $\widetilde{O}({\alpha^{-1}})$ gradient oracle evaluations. 
\end{restatable}
\begin{proof}[Proof sketch (full proof in \cref{appendix:cost_of_computing_ystar_gammastar_inequality})]
    By choosing 
    the penalty parameters $\alpha_1 = \alpha^{-2}$ and $\alpha_2 = \alpha^{-4}$, ~\cref{thm:diff_in_hypergrad_and_gradLagr} guarantees the inexactness of the gradient oracle is upper bounded by $O(\alpha)$.
    
    For the penalty minimization, the objective $\mathcal{L}$ is strongly convex but with a  condition number $\kappa_\mathcal{L} = O(\alpha^{-2})$ due to strong convexity $O(\alpha_1 \mu_g) = O( {\alpha^{-2} \mu_g})$ and smoothness $O(\alpha_2 C_h) = O(\alpha^{-4})$. This leads to $O(\sqrt{\kappa_\mathcal{L}} \log \| \alpha^{-1} \|) = \widetilde{O}({\alpha^{-1}})$ iterations to converge to $\alpha$ accuracy using standard accelerated methods.
    
    Using these error bounds, we can bound the error propagation of inexact solution $y^*_{\lambda,\alpha}$ and inexact gradient oracle using \cref{thm:diff_in_hypergrad_and_gradLagr} by total error $O(\alpha)$.
    Therefore, putting all together, we  need only  $O(\alpha^{-1})$ oracle calls and can maintain $O(\alpha)$ inexactness.
\end{proof}

\cref{thm:cost_of_computing_ystar_gammastar_inequality} shows that we can compute an inexact gradient from ~\cref{alg:inexact-gradient-oracle} with $\Tilde{O}(\alpha^{-1})$ gradient oracles. This concludes the proof of finding $(\delta,\epsilon)$-stationary point of ~\cref{prob:general-constraints} with linear inequality constraints using $\Tilde{O}(\delta^{-1} \epsilon^{-4})$ gradient oracle calls.

\section{Experiments}\label{sec:exp}

\begin{figure}
    \centering
    \begin{subfigure}[b]{0.32\textwidth}
        \centering
        \includegraphics[width=\textwidth]{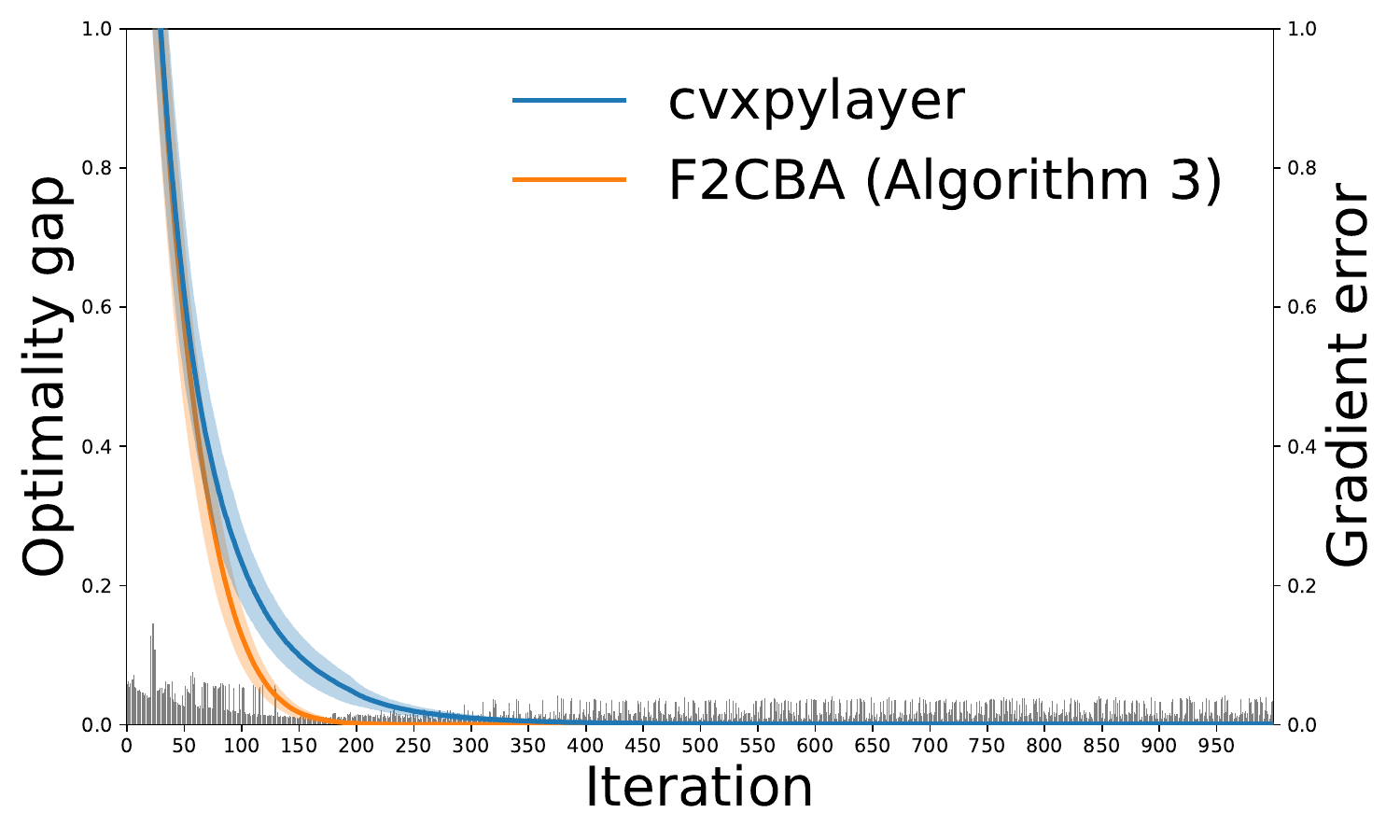}
        \caption{Convergence comparison and gradient error plot.}
        \label{fig:convergence-comparison}
    \end{subfigure}
    \hfill
    \begin{subfigure}[b]{0.32\textwidth}
        \centering
        \includegraphics[width=\textwidth]{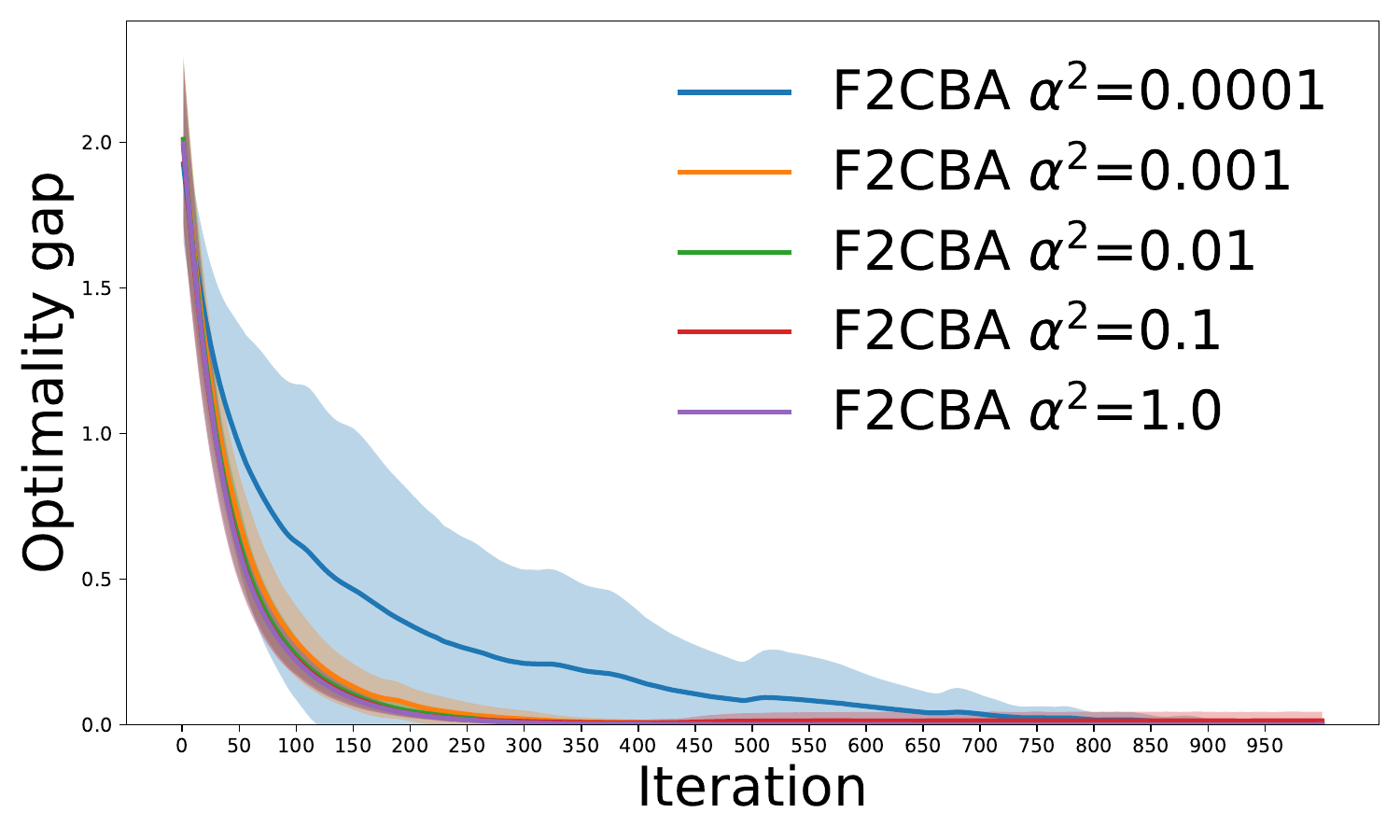}
        \caption{Convergence analysis with varying gradient inexactness $\alpha$.}
        \label{fig:convergence-comparison-2}
    \end{subfigure}
    \hfill
    \begin{subfigure}[b]{0.32\textwidth}
        \centering
        \includegraphics[width=\textwidth]{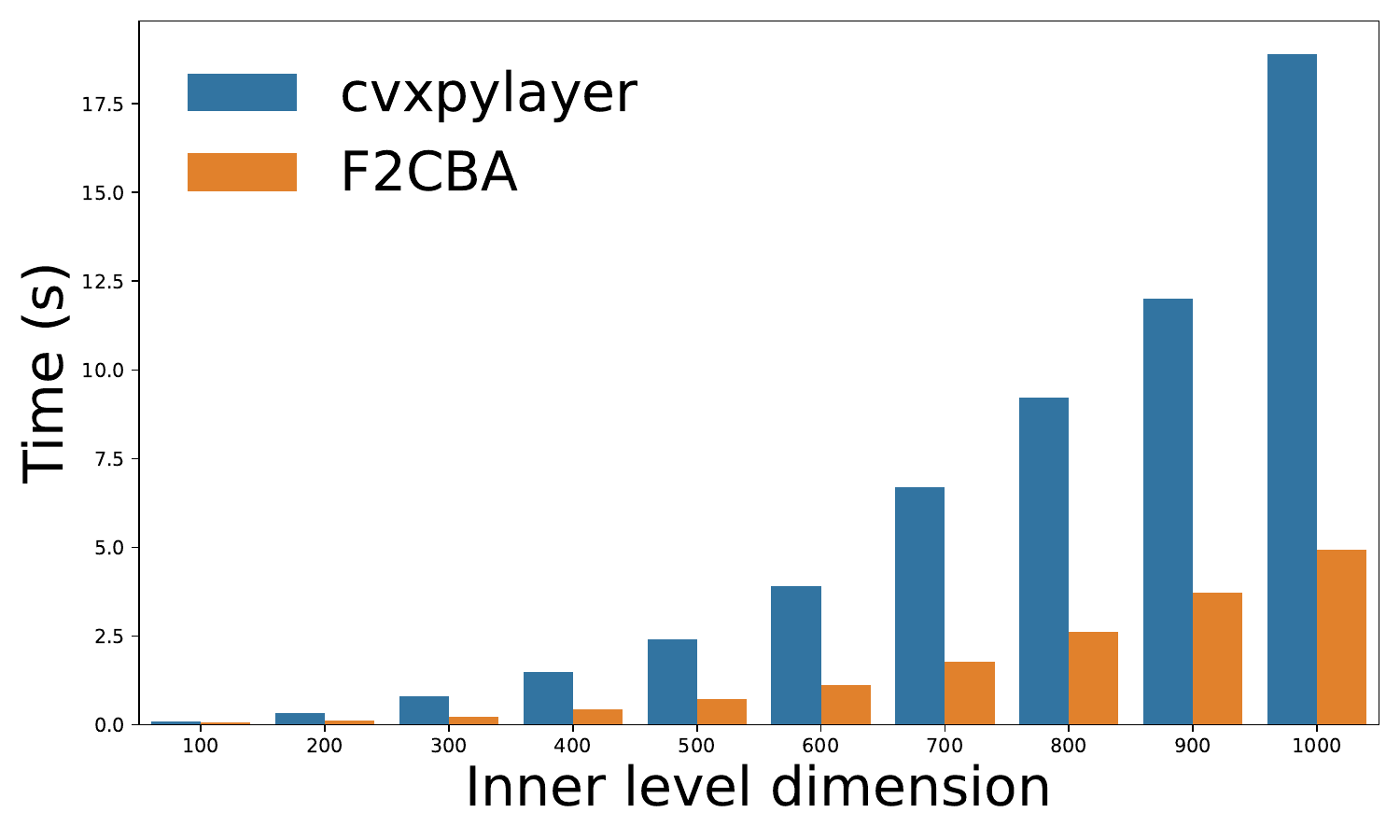}
        \caption{Computation cost per gradient step of varying problem size $d_y$.}
        \label{fig:computation-comparison}
    \end{subfigure}
    \caption{We run \cref{alg: PIGD} using \cref{alg:inexact-gradient-oracle} on the bilevel optimization in the toy example in \cref{eqn:experiment-bilevel-optimization} with $d_x = 100$, $d_y = 200$, $n_{\text{const}} = d_y / 5$, and accuracy $\alpha = 0.1$. \cref{fig:convergence-comparison}, \cref{fig:convergence-comparison-2}, \cref{fig:computation-comparison} vary \# of iterations, gradient exactness $\alpha$, and $d_y$, respectively, to compare the performance under different settings.}
    \label{fig:three graphs}
\end{figure}

We generate instances of the following constrained bilevel optimization problem with inequality constraints: 
\begin{align*}\numberthis\label[prob]{eqn:experiment-bilevel-optimization}
     & \mbox{minimize}_{x} ~  c^\top y^* + 0.01 \norm{x}^2 + 0.01 \norm{y^*}^2 ~\text{ subject to } ~ y^* = \argmin_{y: h(x,y) \leq 0} \frac{1}{2} y^\top Q y + x^\top P y, 
\end{align*}
where $h(x,y) = A y - b$ is a $d_h$-dim linear constraint. The PSD matrix $Q \in \reals^{d_y \times d_y}$, $c \in \reals^{d_y}$, $P \in \reals^{d_x \times d_y}$, and constraints $A \in \reals^{d_h \times d_y}$, $b \in \reals^{d_h}$ are randomly generated from normal distributions (cf. \cref{appendix:experiment-setup}).
We compare our \cref{alg: PIGD} with a non-fully first-order method implemented using \texttt{cvxpyLayer}~\cite{agrawal2019differentiable}. Both algorithms use Adam~\cite{kingma2014adam} to control the learning rate in gradient descent. All the experiments are averaged over ten random seeds.

\looseness=-1\cref{fig:convergence-comparison} shows that  
both the algorithms converge to the same optimal solution at the same rate. Simultaneously, the colorful bars represent the gradient differences between two methods, showing the inexactness of our gradients.
\cref{fig:convergence-comparison-2} additionally varies this inexactness to demonstrate its impact on convergence with standard deviation plotted. 
\cref{fig:computation-comparison} compares the computation costs for different lower-level problem sizes. Our fully first-order method significantly outperforms, in computation cost, the non-fully first-order method implemented using differentiable optimization method.

The implementation can be found in \url{https://github.com/guaguakai/constrained-bilevel-optimization}.

\section{Limitations and future directions}\label{sec:limitation}
One  limitation to our approach is that the inexact gradient oracle we constructed in \cref{sec:inequality-bilevel} requires access to the exact dual multiplier $\lambda^*$.
For a first-order algorithm, the closest proxy one could get to this would be a highly accurate approximation (which could be computed up to $\epsilon$ error within $O(\log(1/\epsilon))$ gradient oracle evaluations). Quantifying the effect of using such an approximate $\lambda^*$ would be an interesting direction of future research. 

Another important direction for future work would be to extend the hypergradient stationarity guarantee of \cref{thm:linIneqInformal} to bilevel programs with \emph{general convex} constraints. 

To this end, we conjecture that the use of a 
primal-only gradient approximation oracle could be potentially effective. 
Finally, our current rate of $\widetilde{O}(\delta^{-1} \epsilon^{-4})$ oracle calls for reaching $(\delta,\epsilon)$-Goldstein stationarity is not necessarily inherent to the problem; indeed, it might be the case that an alternate approach could improve it to the best known rate of ${O}(\delta^{-1} \epsilon^{-3})$ for nonsmooth nonconvex optimization.

\section*{Acknowledgments} 
GK is supported by an Azrieli Foundation graduate fellowship. 
KW is supported by Schmidt Sciences AI2050 Fellowship. SS acknowledges generous support from the Alexander von Humboldt Foundation. SP is supported by NSF CCF-2112665 (TILOS AI Research Institute). 
ZZ is supported by NSF FODSI Fellowship. 

\clearpage
\printbibliography

\newpage
\appendix
\section*{Appendix}

\section{Proofs from \cref{sec:equality-bilevel}} \label{sec:appendix_linear_equality}
In this section, we provide the full proofs of 
claims for bilevel programs with linear equality constraints, as stated in \cref{sec:equality-bilevel}. We first state a few technical results using the implicit function theorem that we repeatedly invoke in our results for this setting. 

\begin{restatable}{lemma}{lemDystarDxLineq}\label{lem:dystarDxLinEq}
    Fix a point $x$. Given $y^*= \arg\min_{y: h(x,y)=0} g(x,y)$ where $g$ is strongly convex in $y$ and $\lamstar$ is the dual optimal variable for this problem, define $\Leqc(x,y,\lam)=g(x,y)+\langle\lam, h(x,y)\rangle$. Then, we have 
\[\underbrace{\begin{bmatrix}\grad^2_{yy}\Leqc(x,\ystar,\lamstar) & \nabla_y h(x,y^{*})^{\top}\\
\nabla_y h(x,y^{*}) & 0
\end{bmatrix}}_{H ~\text{for linear equality constraints}}\begin{bmatrix}\frac{d\ystar}{dx}\\
\frac{d\lamstar}{dx}
\end{bmatrix}=\begin{bmatrix}-\nabla^2_{yx}g(x,y^{*})-\nabla^2_{yx}\langle\lamstar, h(x,y^{*})\rangle\\
-\nabla_x h(x,\ystar)\end{bmatrix}.
\]
\end{restatable}
\begin{proof} 
Since $g$ is strongly convex, by the linear
constraint qualification condition, the KKT condition
is both sufficient and necessary  for optimality.
Hence, consider the
following KKT system obtained via first order optimality of $y^*$, with dual optimal variable $\lamstar$:
\begin{align*}
    \nabla_y g(x,y^{*})+\nabla_y\langle\lamstar,  h(x,y^{*})\rangle  =0, \text{ and } h(x, \ystar)=0.\numberthis\label{eqs:kkt-lin-eq}
\end{align*}
Differentiating 
the system of equations in \cref{eqs:kkt-lin-eq} with respect to $x$ and rearranging terms yields:
\begin{equation}\label{eqs:differentiated-kkt-lin-eq}
    \begin{aligned}
\begin{bmatrix}\nabla^2_{yy} g(x,\ystar) + \nabla^2_{yy}\langle\lamstar,   h(x,\ystar)\rangle  & \nabla_y h(x,y^{*})^{\top}\\
\nabla_y h(x,y^{*}) & 0
\end{bmatrix}\begin{bmatrix}\frac{d\ystar}{dx}\\
\frac{d\lamstar}{dx}
\end{bmatrix}=\begin{bmatrix}-\nabla^2_{yx}g(x,y^{*})-\nabla^2_{yx}\langle\lamstar, h(x,y^{*})\rangle\\
-\nabla_x h(x,\ystar)\end{bmatrix}
    \end{aligned}
\end{equation} 
Finally, noting that $\nabla_{yy}^{2}\Leqc(x,y,\lam)=\nabla^2_{yy} g(x,y)+\nabla^2_{yy}\langle\lam,  h(x,y)\rangle$, we can write \cref{eqs:differentiated-kkt-lin-eq} as claimed. 
\end{proof}

\begin{restatable}{lemma}{lemLineqHinvertibility}\label{lem:non-singular-req} 
Consider the setup in \cref{{lem:dystarDxLinEq}}. The matrix $H$ defined in \cref{eq:def_matrix_H}  is invertible if
the Hessian $\grad_{yy}^{2}\Leqc(x,\ystar,\lamstar):=\nabla^2_{yy} g(x, y^*) +\nabla^2_{yy}\langle\lamstar,  h(x, y^*)\rangle$
satisfies $\nabla_{yy}^{2}\Leqc(x,\ystar,\lamstar)\succ0$ over
the tangent plane $T:=\{y:\nabla_y h(x,y^{*})y=0\}$ and $\nabla_y h$ has full
rank.
\end{restatable}
\begin{proof}
Let $u=[y,\lam].$ We show that $Hu=0$ implies $u=0$, which in turn implies invertibility of $H$. If $\nabla_y h(x,\ystar)y\neq0,$ then by construction of $u$ and $H$, we must also have $Hu\neq0$. Otherwise if $\nabla_y h(x,\ystar)y=0$ and $y\neq0$,
the quadratic form $u^{\top}Hu$ is positive, as seen by 
\[
u^{\top}Hu=y^{\top}\grad^2_{yy}\Leqc(x,\ystar,\lamstar)y>0,
\] where the final step is by the assumption of $\Leqc$ being positive definite over the defined tangent plane $T=\{y:\nabla_y h(x,y^{*})y=0\}$. 
If $y=0$ while $Hu=0$, then $\nabla_y h$ having full rank implies $\lam=0$. Combined with $y=0$, this means $u=0$, as required when $Hu=0$. This concludes the proof. 
\end{proof}

\begin{restatable}{corollary}{corNonSingularH}\label{cor:nonsingularH}
For \cref{prob:lin-eq} under \cref{assumption:linEq_smoothness} and \cref{assumption:eq},
 the matrix $H$  (as defined in \cref{eq:def_matrix_H}) is non-singular. Further, there exists a finite $C_H$ such that $\|H^{-1}\|\leq C_H$. 
\end{restatable}
\begin{proof} Since we are assuming strong convexity of $g$, \cref{{lem:non-singular-req}} applies, yielding the claimed invertibility of $H$. Combined with the boundedness of variables $x$ (per \cref{assumption:eq}) 
and continuity of the inverse implies a bound on $\|H^{-1}\|$. 
\end{proof}

\subsection{Construction of the inexact gradient oracle}
We now show how to construct the inexact gradient oracle for the objective $F$ in \cref{prob:lin-eq}. As sketched in \cref{sec:equality-bilevel}, we then use this oracle in a projected gradient descent algorithm to get the claimed guarantee. 
\begin{restatable}{lemma}{lemLEydelstarCloseToystar}
\label{lem:y-delstar-Lip} Consider  \cref{prob:lin-eq} under \cref{assumption:linEq_smoothness} and \cref{assumption:eq}.
Let  $\ystardel$ be as defined in \cref{eq:lower_perturb}. 
Then, for any $\delta\in[0,\Delta]$ with $\Delta\leq\mu_{g}/2C_{f}$, the following
relation is valid: \[
\|y_{\delta}^{*}-y^{*}\|\leq M(x)\delta,   \textrm{ with } M(x):=\frac{2}{\mu_{g}}\|\grad_{y} f(x,\ystar)\|\leq \frac{2L_f}{\mu_g}.
\]
\end{restatable}
\begin{proof}
    The first-order optimality condition applied to $g(x,y)+\delta f(x,y)$ at $\ystar$ and $\ydelstar$ gives 
\[ \innerprod{\grad_y g(x,\ydelstar)+\delta\grad_y f(x,\ydelstar)}{y^{*}-y_{\del}^{*}} \geq0,\] which upon adding and subtracting $\grad_y f(x,y^{*})$ transforms into
\[ \innerprod{\grad_y g(x,\ydelstar)+\delta[\grad_y f(x,\ydelstar)-\grad_y f(x,y^{*})]+\delta\grad_y f(x,y^{*})}{y^{*}-y_{\del}^{*}}\geq0.\numberthis\label[ineq]{eq:add-sub-fo-gpf}\] Similarly, the first-order optimality condition applied to $g$ at $\ystar$ and $y_{\delta}^{*}$ gives 
\[\innerprod{\grad_y g(x,y^{*})}{y_{\delta}^{*}-\ystar} \geq0.\numberthis\label[ineq]{eq:fo-g}\]
Adding \cref{eq:add-sub-fo-gpf} and \cref{eq:fo-g} and rearranging yields 
\begin{align*}
\innerprod{\grad_y g(x,y_{\delta}^{*})-\grad_y g(x,y^{*})+\delta[\grad_y f(x,\ydelstar)-\grad_y f(x,y^{*})]}{y_{\delta}^{*}-\ystar} & \leq\innerprod{\delta\grad_y f(x,\ystar)}{y^{*}-\ydelstar}.
\end{align*} Applying to the left side above a lower bound via strong convexity of $g+\delta f$ and to the right hand side an upper bound via Cauchy-Schwarz inequality, we have
\[s\|\ydelstar-\ystar\|\leq \delta \|\nabla_y f(x, y^*)\|, \numberthis\label[ineq]{eq:combined-fo-g-gpf}\] where $s$ is the strong convexity of $g+\delta f$. Since $f$ is $C_f$-smooth, the worst case
value of this is $s=\mu_{g}-\delta C_{f}=\mu_{g}-\frac{\mu_{g}}{2C_{f}}C_{f}=\mu_{g}/2$, which when plugged in \cref{eq:combined-fo-g-gpf} then gives the claimed bound. 
\end{proof}

\begin{restatable}{lemma}{lemLinEqLimitFiniteDiffEqualsGradF}\label{lem:lineq-in-limit-finitediff-equals-gradf}
Consider  \cref{prob:lin-eq} under \cref{assumption:linEq_smoothness} and \cref{assumption:eq}.
Then the
following relation is valid.
\[
\lim_{\delta\rightarrow0}\frac{\nabla_{x}[g(x,\ystardel(x))+\lamdeltar h(x,\ystar)]-\nabla_{x}[g(x,y^{*}(x))+\lamstar h(x,y^{*})]}{\delta}=\left(\frac{d y^{*}(x)}{d x}\right)^\top\nabla_{y}f(x,y^{*}(x)).
\]
\end{restatable}
\begin{proof}
Recall that by definition, $g$ is strongly convex and $y^* = \arg\min_{y: h(x,y)=0} g(x,y)$. Hence, we can apply \cref{lem:dystarDxLinEq}. Combining this with \cref{lem:non-singular-req} and further applying that linearity of $h$ implies $\nabla^2_{yy}h = 0$ and $\nabla^2_{xy}h=0$, we obtain the following: 
\[ \begin{bmatrix}\frac{d\ystar}{dx}\\
\frac{d\lamstar}{dx}
\end{bmatrix}=\begin{bmatrix}\nabla^2_{yy} g(x,y^{*}) & \nabla_y h(x,y^{*})^{\top}\\
\nabla_y h(x,y^{*}) & 0
\end{bmatrix}^{-1}\begin{bmatrix}-\nabla^2_{yx}g(x,y^{*})\\
-\nabla_x h(x,\ystar)
\end{bmatrix}.\]
So we can express the right-hand side of the claimed equation in the lemma statement by 
\begin{align*}
\left(\frac{d y^{*}(x)}{d x}\right)^\top\nabla_{y}f(x,y^{*}(x))&=\begin{bmatrix} \left(\frac{dy^*}{dx}\right)^\top & \left(\frac{d\lamstar}{dx}\right)^\top\end{bmatrix}\begin{bmatrix}\nabla_{y}f(x,y^{*}(x))\\ 
0\end{bmatrix},
\end{align*} which can be further simplified to \[\begin{bmatrix}-\nabla^2_{yx} g(x,y^{*})^\top & -\nabla_x h(x,\ystar)^\top\end{bmatrix}\begin{bmatrix}\nabla^2_{yy} g(x,y^{*}) & \nabla_y h(x,y^{*})^{\top}\\
\nabla_y h(x,y^{*}) & 0
\end{bmatrix}^{-1}\begin{bmatrix}\nabla_{y}f(x,y^{*}(x))\\ 
0\end{bmatrix}.\numberthis\label{prop3eq:RHS}\]
We now apply \cref{lem:dystarDxLinEq} to the perturbed problem defined in \cref{eq:lower_perturb}.
We know from \cref{lem:y-delstar-Lip} that $\lim_{\delta\rightarrow0}\ydelstar=\ystar$. 
The associated KKT system is given by
\begin{align*}
\delta \nabla_y f (x,\ydelstar)+\nabla_y g(x,\ydelstar)+\nabla_y \langle\lamdeltar,  h(x,\ydelstar)\rangle  =0 \text{ and }
h(x,\ydelstar) =0. \numberthis\label{eqs:lineq-kkt-perturbed}
\end{align*}
Taking the derivative with respect of \cref{eqs:lineq-kkt-perturbed} gives the following  implicit  system, where we used the fact that $h$ is linear and hence $\nabla^2_{yy} h=0$:  
\begin{equation}
\underbrace{\begin{bmatrix}\delta \nabla^2_{yy} f(x,\ydelstar)+\nabla^2_{yy} g(x,\ydelstar) & \nabla_y h(x,\ydelstar)^{\top}\\
\nabla_y h(x,\ydelstar) & 0
\end{bmatrix}}_{H_{\delta}}\begin{bmatrix}\frac{d\ydelstar}{d\delta}\\
\frac{d\lamdeltar}{d\delta}
\end{bmatrix}=\begin{bmatrix}-\nabla_y f(x,\ydelstar)^\top\\
0
\end{bmatrix}.\label{eq:implicit_function}
\end{equation}
For a sufficiently small  $\delta$, we have $\nabla^2_{yy} g(x,\ystardel)+\delta  \nabla^2_{yy} f(x,\ydelstar)\succeq \tfrac{\mu_{g}}{2}I$, which implies  invertibility of $H_{\delta}$ by an application of \cref{{lem:non-singular-req}}. 
Since \cref{lem:y-delstar-Lip} implies  $\lim_{\delta\rightarrow0}\ydelstar=\ystar$, we get 
\[
\begin{bmatrix}\frac{d\ydelstar}{d\delta}\\
\frac{d\lamdeltar}{d\delta}
\end{bmatrix}|_{\delta=0}=\begin{bmatrix}\nabla^2_{yy} g(x,\ystar) & \nabla_y h(x,\ystar)^{\top}\\
\nabla_y h(x,\ystar) & 0
\end{bmatrix}^{-1}\begin{bmatrix}-\nabla_y f(x,\ystar)\\
0
\end{bmatrix}.
\]
So we can express the left-hand side of the expression in the lemma statement by 
\begin{align*}
&\lim_{\delta\rightarrow0}\frac{\nabla_{x}[g(x,\ystardel(x))+\langle\lamdeltar, h(x,y^*)\rangle]-\nabla_{x}[g(x,y^{*}(x))+\langle\lamstar, h(x,y^{*})\rangle]}{\delta}\\
&= \nabla^2_{xy} g(x,y^*) \frac{d\ydelstar}{d\delta} + \nabla_x h(x,y^*)^\top \frac{d\lamdeltar}{d\delta}\\ 
&=\begin{bmatrix}\nabla^2_{xy} g(x,\ystar) & \nabla_x h(x,\ystar)^\top\end{bmatrix}\begin{bmatrix}\nabla^2_{yy} g(x,\ystar) & \nabla_y h(x,\ystar)^{\top}\\
\nabla_y h(x,\ystar) & 0
\end{bmatrix}^{-1}\begin{bmatrix}-\nabla_y f(x,\ystar)\\
0
\end{bmatrix},
\end{align*}
which matches \cref{prop3eq:RHS} (since $(\nabla^2_{yx} g)^\top=\nabla^2_{xy}g$), thus concluding the proof.
\end{proof}


\lemYdelstarLamdelstarSmooth*
\begin{proof} 
Rearranging \cref{{{eq:def_matrix_H}}} and applying \cref{cor:nonsingularH}, we have \[
\begin{bmatrix}\frac{d\ystar}{dx}\\
\frac{d\lamstar}{dx}
\end{bmatrix}=\begin{bmatrix}\nabla^2_{yy} g(x,y^{*}) & B^{\top}\\
B & 0
\end{bmatrix}^{-1}\begin{bmatrix}-\nabla^2_{yx} g(x,y^{*})\\
-\nabla_x h(x,\ystar)
\end{bmatrix}. 
\] This implies a Lipschitz bound of $C_H \cdot (\gssmooth + \|A\|)$. 
    Next, note that in the case with linear equality constraints, the terms in  \cref{{{eqs:differentiated-kkt-lin-eq}}}  involving second-order derivatives of $h$ are all zero; differentiating \cref{{{eqs:differentiated-kkt-lin-eq}}}    with respect to $x$, we notice that the linear system we get again has the same matrix $H$ from before. We can therefore again perform the same inversion and apply the bound on $\|H^{-1}\|$ and on the third-order derivatives of $g$ (\cref{assumption:eq}) 
to observe that $\|\frac{d^2y^*}{dx^2}\|\leq O(C_H \cdot \gtsmooth \|\frac{dy^*}{dx}\|^2)= O(C_H^3\cdot\gtsmooth\cdot(\gssmooth+\|A\|)^2)$, where we are hiding numerical constants in the Big-Oh notation. 

As a result, we can calculate the Lipschitz smoothness constant associated with the hyper-objective $F$ by
\begin{align*}
&\|\nabla F(x) - \nabla F(\bar x)\| \\
&\leq \|\frac{dy^*(x)}{dx}\nabla_y f(x,y^*(x))-\frac{dy^*(\bar x)}{dx}\nabla_y f(\bar x,y^*(\bar x))\|+ \|\nabla_x f(x,y^*(x)) - \nabla_x f(\bar x, y^*(\bar x))\|\\
&\leq  [C_fC_H (L_g + \|A\|) + C_f C^2_H (L_g + \|A\|)^2+ L_f C_H^3 S_g (L_g+\|A\|)^2]\|x-\bar x\| \\
&\ \ +[C_f + C_f C_H (L_g +\|A\|)] \|x-\bar x\|\\
&\leq \underbrace{ 2(L_f +C_f+C_g)C_H^3 S_g (L_g +\|A\|)^2}_{C_F}\|x-\bar x\|.
\end{align*}

\end{proof}

\lemLineqFiniteDiffEqualsGradF*
\begin{proof}
For simplicity, we adopt the following notation throughout this proof:
$g_{xy}(x,y) = \grad^2_{xy} g,$ and $g_{xyy}$ denotes the tensor such that its $ijk$ entry is given by $\frac{\partial^3 g}{\partial x_i \partial y_j \partial y_k}$. 
We first consider the terms involving $g$. 
By the fundamental theorem of calculus, we have  \[ \nabla_{x}g(x,\ystardel(x))-\nabla_{x}g(x,y^{*}(x)) =\int_{t=0}^{\delta}g_{xy}(x,y_{t}^{*}(x))\frac{dy_{t}^{*}(x)}{dt} dt. \] As a result, we have 
\begin{align*}
&\frac{\nabla_{x}g(x,\ystardel(x))-\nabla_{x}g(x,y^{*}(x))}{\delta}-g_{xy}(x,y^{*}(x))\frac{dy_{t}^{*}(x)}{dt}|_{t=0}\\
& =\frac{1}{\delta}\int_{t=0}^{\delta}\left(g_{xy}(x,y_{t}^{*}(x))\frac{dy_{t}^{*}(x)}{dt}-g_{xy}(x,y^{*}(x))\frac{dy_{t}^{*}(x)}{dt}|_{t=0} \right)dt\\
 & =\frac{1}{\delta}\int_{t=0}^{\delta}\left(g_{xy}(x,y_{t}^{*}(x))\frac{dy_{t}^{*}(x)}{dt}-g_{xy}(x,y^{*}(x))\frac{dy_t^{*}(x)}{dt}|_{t=0}\right) dt\\
 & =\frac{1}{\delta}\int_{t=0}^{\delta}\left(g_{xy}(x,y_{t}^{*}(x))-g_{xy}(x,y^{*}(x))\right) \frac{dy_{t}^{*}(x)}{dt} dt +\frac{1}{\delta}\int_{t=0}^{\delta}g_{xy}(x,y^{*}(x))\cdot\left(\frac{dy_{t}^{*}(x)}{dt}-\frac{dy_t^{*}(x)} {dt}|_{t=0}\right)dt.\numberthis\label{eq:LE-findiff-to-ftimespartial} \end{align*} We now bound each of the terms on the right-hand side of \cref{eq:LE-findiff-to-ftimespartial}. 
 For the first term, we have  
 \begin{align*}
&\|\frac{1}{\delta}\int_{t=0}^{\delta}\left(g_{xy}(x,y_{t}^{*}(x))-g_{xy}(x,y^{*}(x))dt\right) \frac{dy_{t}^{*}(x)}{dt}\|\\
&\leq\frac{1}{\delta}\int_{t=0}^{\delta}\|\frac{dy_{t}^{*}(x)}{dt}\|\cdot\int_{s=0}^{t}\|g_{xyy}(x,y_{s}^{*}(x))\|\|\frac{dy_{s}^{*}(x)}{ds}\|ds\cdot dt\\
&\leq\frac{1}{\delta}\int_{t=0}^{\delta}\|\frac{dy_{t}^{*}(x)}{dt}\|\cdot\max_{s\in[0,\delta]}\|g_{xyy}(x,y_{s}^{*}(x))\|\cdot\|\frac{dy_{s}^{*}(x)}{ds}\|tdt\\
&\leq\frac{1}{\delta}\cdot\max_{u\in[0,\delta]}\|g_{xyy}(x,y_{u}^{*}(x))\|\cdot\delta^{2}\cdot\max_{t\in[0,\delta]}\|\frac{dy_{t}^{*}(x)}{dt}\|^{2}\\&\leq\del\cdot\max_{u\in[0,\delta]}\|g_{xyy}(x,y_{u}^{*}(x))\|\cdot\max_{t\in[0,\delta]}\|\frac{dy_{t}^{*}(x)}{dt}\|^{2}
 \\
 &= \del\cdot\gtsmooth\cdot \ystarliplineq^2,\numberthis\label[ineq]{eq:finite-diff-grad-first_term}
 \end{align*} where $\ystarliplineq$ is the Lipschitz bound on $y^*$ as shown in \cref{lem:smoothness_of_ydelstar_lamdelstar}, and $\gtsmooth$ is the smoothness of $g$ from
 \cref{assumption:eq}.
 For the second term on the right-hand side of \cref{eq:LE-findiff-to-ftimespartial}, we have  
 \begin{align*}
\|\frac{1}{\delta}\int_{t=0}^{\delta}g_{xy}(x,y^{*}(x))\cdot\left(\frac{dy_{t}^{*}(x)}{dt}-\frac{dy^{*}(x)}{dt}\right)\| & \leq\frac{1}{\delta}\cdot\|g_{xy}(x,y^{*}(x))\|\cdot\int_{t=0}^{\del}\left(\int_{s=0}^{t}\|\frac{d^{2}}{ds^{2}}y_{s}^{*}(x)\|ds\right)dt \\
 & \leq\frac{1}{\del}\cdot\|g_{xy}(x,y^{*}(x))\|\cdot\max_{s\in[0,\delta]}\|\frac{d^{2}}{ds^{2}}y_{s}^{*}(x)\|\cdot\delta^{2} \\
 & \leq\delta\cdot\|g_{xy}(x,y^{*}(x))\|\cdot\max_{s\in[0,\delta]}\|\frac{d^{2}}{ds^{2}}y_{s}^{*}(x)\|\\
 &= \delta\cdot \gssmooth\cdot \ystarsmoothlineq, \numberthis\label[ineq]{eq:finite-diff-grad-second_term} 
\end{align*} where $\gssmooth$ is the bound on smoothness of $g$ as in
\cref{assumption:eq},
and $\ystarsmoothlineq$ is the bound on $\|\frac{d^2y^*}{dx^2}\|$ from \cref{lem:smoothness_of_ydelstar_lamdelstar}. 
For the terms involving the function $h$, we have \begin{align*}
\|\frac{\lamdeltar-\lam^{*}}{\delta}-\frac{d\lam_{\delta}^{*}}{d\del}|_{\delta=0}\| & =\frac{1}{\delta}\int_{t=0}^{\del}\|\frac{d\lam_{t}^{*}}{dt}-\frac{d\lam_{\delta}^{*}}{d\del}|_{\delta=0}\|dt\\
 & =\frac{1}{\delta}\int_{t=0}^{\delta}\int_{s=0}^{t}\|\frac{d^{2}}{ds^{2}}\lam_{s}^{*} \| ds\cdot dt\\
 & \leq\frac{1}{\delta}\max_{s\in[0,\delta]}\|\frac{d^{2}}{ds^{2}}\lam_{s}^{*}\|\cdot\delta^{2}\leq\delta\cdot\max_{s\in[0,\delta]}\|\frac{d^{2}}{ds^{2}}\lam_{s}^{*}\|\\
 &= \delta\cdot \lamstarsmoothlineq,\numberthis\label[ineq]{eq:finite-diff-grad-third-term}
\end{align*} where $\lamstarsmoothlineq$ is the bound on $\|\frac{d^2 \lam^*}{ds^2}\|$ from \cref{lem:smoothness_of_ydelstar_lamdelstar}. 
Combining  \cref{eq:LE-findiff-to-ftimespartial}, \cref{eq:finite-diff-grad-first_term}, \cref{eq:finite-diff-grad-second_term}, and \cref{eq:finite-diff-grad-third-term}, along with \cref{{lem:lineq-in-limit-finitediff-equals-gradf}}, \cref{cor:nonsingularH}, and  \cref{{lem:smoothness_of_ydelstar_lamdelstar}}, we have that overall bound is \begin{align*}\delta\cdot( \gtsmooth \ystarliplineq^2 + \gssmooth \ystarsmoothlineq + \lamstarsmoothlineq) &\leq O(\delta\cdot(\gtsmooth \cdot C_H^3 \cdot (\gssmooth + \|A\|)^2\cdot(\gssmooth+C_f + L_f))).\end{align*} 
\end{proof}

\subsection{The cost of linear equality constrained bilevel program}\label{sec:LEQ-main-thm-full-proof}

We present here our full (implementable) algorithm for \cref{{prob:lin-eq}} and restate and prove its convergence guarantee. 
\begin{algorithm}[h]\caption{The Fully First-Order Method for Bilevel Equality Constrained Problem}\label{alg:LE-full-alg}
\begin{algorithmic}[1]
\State \textbf{Input:}
Current $x_0$, accuracy $\epsilon$, perturbation $\delta = \epsilon^2/8C^2_F R_\mathcal{X}$ with $C_F= 2(L_f +C_f+C_g)C_H^3 S_g (L_g +\|A\|)^2$,  accuracy for the lower level problem $\tilde\delta = 2(C_g +\|A\|)\delta^2$.
\For{t=0,1,2,...}

\State Run \cref{alg:LE-approximate-prima-dual-solution}   to generate $\tilde \delta$-accurate primal and dual solutions $(\hat{y}^*, \hat{\lambda}^*)$ for $$\min_{y: Ax_t+By=b} g(x_t,y)$$
\State  Run \cref{alg:LE-approximate-prima-dual-solution} to generate $\tilde \delta$-accurate primal and dual solutions $(\hat{y}_\delta^*, \hat{\lambda}_\delta^*)$ for 
$$\min_{y: Ax_t+By=b} g(x_t,y)+\delta f(x_t, y)$$ 

\State Compute $\hat{v}_t:= \frac{\nabla_{x}[g(x_t,\hat{y}_\delta^*)+\hat{\lambda}_\delta^* h(x,\hat{y}^*)]-\nabla_{x}[g(x_t,\hat{y}^*)+\hat{\lambda}^* h(x,\hat{y}^*)]}{\delta}$, set $$\widetilde{\nabla} F(x_t) := \hat{v}^t + \nabla_x f(x, \hat{y}^*(x)).$$
\State Set $x_{t+1} \leftarrow \arg\min_{z\in \mathcal{X}} \| z- (x_{t}-\frac{1}{C_F}\widetilde{\nabla} F(x_{t}))\|^2.$
\EndFor

\end{algorithmic}
\end{algorithm}

\linEqFullCost*
\begin{proof}

We first show the inexact gradient $\widetilde{\nabla}F(x_t)$ generated in \cref{alg:LE-full-alg} is an $\delta$-accurate approximation to the hyper-gradient $\nabla F(x_t).$ Consider the inexact gradient defined in \eqref{eq:part-hypergrad-approx-lin-eq}
\begin{align*}
    \|v_{t}-\hat{v}_{t}\|&\leq\frac{1}{\delta}\{\|[\nabla_{x}g(x_{t},\hat{y}_{\delta}^{*})-\nabla_{x}[g(x_{t},\hat{y}^{*})]-[\nabla_{x}g(x_{t},y_{\delta}^{*})-\nabla_{x}[g(x_{t},y^{*})\|\\&\,\,\,+\ensuremath{\|\hat{\lambda}_{\delta}^{*}-\hat{\lambda}^{*}-[\lambda_{\delta}^{*}-\lambda^{*}\|\|A\|\}}\\&\leq\frac{2}{\delta}[C_{g}+\| A\|]\tilde{\delta}.
\end{align*}
Thus we get 
\begin{align*}
    \|\widetilde{\grad}F(x_{t})-\nabla F(x_{t})\|&\leq\|\grad_{x}f(x_{t},y^{*})-\grad_{x}f(x_{t},\hat{y}^{*})\|+\norm{\hat{v}^{t}-v^{t}}+\|v^{t}-\frac{dy^{*}(x^{t})}{dx}\grad_{y}f(x_{t},y^{*}(x_{t}))\|\\&\leq C_{f}\tilde{\delta}+\frac{2}{\delta}[C_{g}+\norm A]\tilde{\delta}+C_{F}\delta\\&\leq\frac{2\tilde{\delta}}{\delta}[C_{f}+C_{g}+\|A\|]+C_{F}\delta\\&\leq\frac{\epsilon^{2}}{4C_{F}R_{\mathcal{X}}}.
\end{align*}

Applied to the $C_F$-smooth hyper-objective $F$, such an inexact gradient oracle satisfies the requirement for  \cref{pr:inexact-pgd}. Thus an $\epsilon$-stationary point with $\|\mathcal{G}_F(x^t)\|\leq \epsilon$ (see Eq. \eqref{eq:gradient-mapping}  for the definition of gradient mapping) must be found in  $N=O(\frac{C_F (F(x^0)-F^*)}{\epsilon^2})$ iterations. Noting the evaluation of inexact solutions $(\hat y^*, \hat \lambda^*, \hat y_\delta^*, \hat \lambda_\delta^* )$ requires $\tilde{O}(\sqrt{C_g/\mu_g})$ first order oracle evaluations, we arrive at the total oracle complexity of $\tilde{O}(\sqrt{C_g/\mu_g}\frac{C_F (F(x^0)-F^*)}{\epsilon^2})$ for finding an $\epsilon$-stationary point. 
\end{proof}
\subsection{The cost of inexact projected gradient descent method}
In this subsection, we state the number of iterations required by  projected gradient descent method to find an $\epsilon$-stationary point using inexact gradient oracles. Specifically, we consider the following non-convex smooth problem  where the objective $F$ is assumed to be $C_F$-Lipschitz smooth: 
\begin{equation}\label{eq:prob-smooth-constrained}
    \mbox{minimize}_{x\in \mathcal{X}} F(x).
\end{equation}
Since the feasible region $\mathcal{X}$ is compact, we use the norm of the following gradient mapping $\mathcal{G}_F(x)$ as the stationarity criterion 
\begin{equation}\label{eq:gradient-mapping}
    \mathcal{G}_F(x):= {C_F}(x - x^+) \text{ where } x^+ = \arg\min_{z\in \mathcal{X}} \left\| z-\left(x -\frac{1}{C_F} \nabla F(x)\right)\right\|^2.
\end{equation}
Initialized to some $x_0$ and the inexact gradient oracle $\widetilde{\nabla} F$, the updates of the inexact projected gradient descent method is given by 
\begin{equation}\label{alg:inexact-projected-gd-algorithm}
    \begin{split}
        \text{\textbf{For}}&\text{ t=1,2,..., N \textbf{do}:}\\
        &\text{Set } x_t \leftarrow \arg\min_{z\in \mathcal{X}} \left\| z- \left(x_{t-1}-\frac{1}{C_F}\widetilde{\nabla} F(x_{t-1})\right)\right\|^2.
    \end{split}
\end{equation}
The next proposition calculates the complexity result. 
\begin{proposition}\label{pr:inexact-pgd}
    Consider the constrained optimization problem in \eqref{eq:prob-smooth-constrained} with $F$ being $C_F$-Lipschitz smooth and $\mathcal{X}$ having a radius of $R$. When supplied with a $\delta=\epsilon^2/4C_F R$ -inexact gradient oracle $\widetilde{\nabla} F$, that is, $\|\nabla F(x)-\widetilde{\nabla}F(x)\|\leq \delta$, the solution generated by the projected gradient descent method \eqref{alg:inexact-projected-gd-algorithm} satisfies 
    $$\min_{t \in [N]} \| \mathcal{G}_F(x_t)\|^2 \leq \frac{C_F(F(x_0)-F^*)}{N}+\delta C_F R,$$
    that is, it takes at most $O(\frac{C_F (F(x^0)-F^*)}{\epsilon^2})$ iterations to generate some $\bar x$ with $\|\mathcal{G}_F(x)\|\leq \epsilon$.
\end{proposition}
\begin{proof}
    By $C_F$-smoothness of $F$, we have 
\begin{align*}
f(x_{t+1})  =f(x_t-\frac{1}{C_F}\widetilde{\mathcal{G}_F}(x_t))
 & \leq f(x_t)-\frac{1}{C_F}\widetilde{\mathcal{G}_F}(x_t)^{\top}\nabla f(x_t)+\frac{1}{2C_F}\|\widetilde{\mathcal{G}_F}(x_t)\|^{2}\\
 & =f(x_t)-\frac{1}{2C_F}\|\widetilde{\mathcal{G}_F}(x_t)(x_t)\|^{2}+\frac{1}{C_F}\widetilde{\mathcal{G}_F}(x_t)^{\top}(\widetilde{\mathcal{G}_F}(x_t)-\nabla f(x_t)).\numberthis\label[ineq]{eq:fn_decrease_lineq_pgd}
\end{align*}
We now show that $\frac{1}{\beta}\widetilde{\mathcal{G}_F}(x_t)^{\top}(\widetilde{\mathcal{G}_F}(x_t)-\nabla f(x_t))\leq0.$ Let $\widetilde{y}_t = x_t-\frac{1}{C_F}\widetilde{\nabla}F(x_t)$, and let $y_t = x_t-\frac{1}{C_F}{\nabla}f(x_t)$. 
Then have that 
\begin{align*}
\frac{1}{C_F}\widetilde{\mathcal{G}_F}(x_t)^{\top}(\frac{1}{C_F}\widetilde{\mathcal{G}_F}(x_t)-\nabla f(x_t)) & =C_F(x_t-\pq(\widetilde{y}_t))^{\top}(y_t - \pq(\widetilde{y}_t))\\
 & = C_F(x_t-\pq(\widetilde{y}_t))^\top (\widetilde{y}_t - \pq(\widetilde{y}_t)) \\
 &\quad +C_F(x_t-\pq(\widetilde{y}_t))^\top(y_t - \widetilde{y}_t)\\
&\leq C_F(x_t-\pq(\widetilde{y}_t))^\top(y_t - \widetilde{y}_t)\\
&\leq \delta C_F  R,
\end{align*}
where the penultimate inequality uses the fact that $\mathcal{X}$ is a convex set, and $R$ is the diameter of the set $X$. Combining this with \cref{eq:fn_decrease_lineq_pgd}, we have that the function decrease per iteration is 
\[ F(x_{t+1})\leq F(x_t)-\frac{1}{2C_F}\|\widetilde{\mathcal{G}_F}(x_t)\|^{2} + \delta C_F R. \] 

Summing over $N$ iterations telescopes the terms, we get 
$$\min_{t\in[N]}\|\widetilde{\mathcal{G}_F}(x_t)\|^{2}\leq \frac{1}{N} C_F (F(x^0)-F^*) + \delta C_F R.$$
Substituting in $N=\frac{4}{\epsilon^2}C_F (F(x^0)-F^*)$  and the choice of $\delta = \epsilon^2/4C_F R$, we get 
$$\min_{t\in[N]}\|\widetilde{\mathcal{G}_F}(x_t)\|^{2}\leq \frac{\epsilon^2}{2}.$$
Taking into account the fact that $\|\widetilde{\mathcal{G}_F}(x_t)- \mathcal{G}_F (x_t) \|\leq \|\nabla F(x^t)-\widetilde{\nabla} F(x^t)\|\leq \delta$,  we obtain the desired result. 

\end{proof}

\subsection{The cost of generating approximate solutions of \cref{{prob:lin-eq}}}\label{sec:LEQ-cost-computing-ystar-lamstar}
In this subsection, we address the issue of generating approximations to the primal and dual solutions  $(y^*,\lambda^*)$ associated with the lower-level problem in \cref{{{prob:lin-eq}}}. These approximations are required for computing the approximate hypergradient in \cref{alg:LE-inexact-gradient-oracle}.  For notational simplicity, we are going to consider the following constrained strongly convex problem: 
\[ 
\begin{array}{ll}
    \mbox{minimize}_{y\in\R^d} &g(y)\\
    \mbox{subject to } & By=b.
\end{array}\numberthis\label[prob]{eq:simple-linear-cosntrained-problem}
\] 

We propose the following simple scheme to generate approximate solutions to \cref{{eq:simple-linear-cosntrained-problem}}. 
\begin{center}
\fbox{\begin{varwidth}{\dimexpr\textwidth-2\fboxsep-2\fboxrule\relax}
  Compute a feasible $\hat{y}$  such that $\|\hat{y}-y^*\|\leq\delta$.
  Then solve
  
     \begin{equation}\label{eq:approximate-lam-hat}
        \hat\lambda = \arg\min_{\lambda \in \R^m} \|\nabla_y g(\hat y)-B^\top\lambda\|^2.
    \end{equation}
    
\end{varwidth}}
\end{center}

The following lemma tells us that $\hat \lambda$ is close to $\lambda^*$ if $B$ has full row rank. 

\begin{lemma}\label{lm:generating-lamhat}
    Suppose $g$ in \cref{eq:simple-linear-cosntrained-problem} is a $C_g$-Lipschitz smooth, and the matrix $B$ has full row rank such that the following matrix $M_{B}$ is invertible
    $$M_B =\begin{bmatrix}
        I & B^\top\\
        B & 0
    \end{bmatrix}.$$
    Then the approximate solution $(\hat \lambda, \hat y)$ from \cref{eq:approximate-lam-hat} satisfies $\|\hat\lambda -\lambda^*\|\leq \|M_{B}^{-1}\| (1+C_g)\delta$.
\end{lemma}

\begin{proof}
    Since $(\lambda^*,y^*)$ satisfy the KKT conditions, they are the solution to the following linear system 
    \begin{equation}\label{lmeq:ystarlamstar-equation}
        \underbrace{\begin{bmatrix}I & B^\top\\ B & 0 \end{bmatrix}}_{=M_B} 
    \begin{bmatrix} y^* \\ \lambda^*\end{bmatrix}
    =\begin{bmatrix}-\nabla_y g(y^*) + I y^* \\ b\end{bmatrix}.
    \end{equation}
    That is 
    $$ \begin{bmatrix} y^* \\ \lambda^*\end{bmatrix}
    =M_B^{-1}\begin{bmatrix}-\nabla_y g(y^*) + I y^* \\ b\end{bmatrix}.
    $$
    On the other hand, the approximate solutions $(\hat{y}, \hat{\lambda})$ in \cref{eq:approximate-lam-hat} satisfies 
        $$\begin{bmatrix}I & B^\top\\ B & 0 \end{bmatrix}
    \begin{bmatrix} \hat{y} \\ \hat{\lambda}\end{bmatrix}
    =\begin{bmatrix} B^\top \hat{\lambda} + I \hat{y} \\ b\end{bmatrix}.$$
    We show the right hand side (r.h.s) of the above equation to be close to the r.h.s of \cref{lmeq:ystarlamstar-equation}. Let $S:=\{B^\top \lambda: \lambda\in \R^m\}$ denote the subspace spanned by the rows of $B$. We can rewrite $B^\top \hat\lambda$ as the projection of $\nabla g(\hat y)$ onto $S$, that is, 
    \begin{align*}
        &B^\top \hat \lambda = \arg\min_{s\in S}\|\nabla_y g(\hat y) -s\|^2\\
     -\nabla_y g( y^*)=&B^\top \lambda^* = \arg\min_{s\in S}\|\nabla_y g( y^*) -s\|^2,
    \end{align*}
    where the second relation follows from the KKT conditon associated with $(\lambda^*,y^*)$. Since the projection is an non-expansive operation, we have 
    $$\|B^\top \hat\lambda -(- \nabla_y g(y^*))\|=\|B^\top \hat \lambda - B^\top \lambda^*\|\leq \|\nabla_y g(\hat y) -\nabla g(y^*)\|\leq C_g \|\hat y - y^*\|\leq C_g \delta.$$
    We can rewrite $(\hat y, \hat \lambda)$ as solutions to the following linear system with some $\|\tau\|\leq (1+C_g)\delta$, 
       $$ \begin{bmatrix} \hat{y} \\ \hat \lambda\end{bmatrix}
    =M_B^{-1}\begin{bmatrix}-\nabla_y g(y^*) + I y^* +\tau \\ b\end{bmatrix}.
    $$
    Thus we get 
     $$ \|\begin{bmatrix} \hat{y} \\ \hat \lambda\end{bmatrix}-\begin{bmatrix} y^* \\ \lambda^*\end{bmatrix}\|
    =\|M_B^{-1}\|\|\begin{bmatrix}\tau \\ 0\end{bmatrix}\leq \|M_B^{-1}\|(1+C_g)\delta.
    $$    
\end{proof}

Now we can just use the AGD method to generate a close enough approximate solution $\hat y$ and call up the Subroutine in \cref{eq:approximate-lam-hat} to generate the approximate dual solution $\hat \lambda.$

\begin{algorithm}[h]\caption{The Projected Gradient Method to Generate Primal and Dual Solutions for a Linearly Constrained Problem}\label{alg:LE-approximate-prima-dual-solution}
\begin{algorithmic}[1]
\State \textbf{Input}: accuracy requirement $\epsilon>0$ and linearly constrained problem $\min_{y: By =b} g(y)$.
 \State Starting from $y^0=0$ and using $Y:=\{y\in\R^d: By=b\}$ as the simple feasible region. 
 \State Run the Accelerated Gradient Descent (AGD) Method (Section 3.3 in \cite{lan2020first}) for $N=\lceil 4\sqrt{C_g/\mu_g}\log\frac{\|y^*\| \|M^{-1}_B\|(C_g+1)}{\mu_g\epsilon} \rceil$ iterations.
  \State Use the $y^N$ as the approximate solution $\hat y$ to generate $\hat \lambda$ according to \cref{eq:approximate-lam-hat}.
  \State \Return $(\hat y, \hat \lambda)$
\end{algorithmic}
\end{algorithm}

\begin{proposition}
    Suppose the objective function $g$ is both $L_g$-smooth and $\mu_g$-strongly convex, and that the constraint satisfies the assumption in \cref{lm:generating-lamhat}. Fix an $\epsilon>0$, the solution  $(\hat y, \hat \lambda)$ returned by the above procedure satisfies $\|y^* - \hat y\|\leq \epsilon$ and $\|\hat \lambda - \lambda^*\|\leq \epsilon$. In another words, the cost of generating $\epsilon$-close primal and dual solutions are bounded by $O(\sqrt{\frac{C_g}{\mu_g}}\log{\frac{1}{\epsilon}}).$
\end{proposition}
\begin{proof}
With $N:=\lceil 4\sqrt{C_g/\mu_g}\log\frac{\|y^*\| \|M^{-1}_B\|(L_g+1)}{\mu_g\epsilon} \rceil$, Theorem 3.7 in \cite{lan2020first} shows that $\|y^N - \hat y\|\leq \epsilon/ \|M_{B}^{-1}\| (1+L_g)$. Then we can apply \cref{lm:generating-lamhat} to obtain the desired bound.
\end{proof}
\section{Proofs from \cref{{sec:nonsmooth}}}\label{sec:zeroth-order-algs}
In this section we present a variant of \cref{alg: OIGRM}, which as stated in \cref{{sec:nonsmooth}},  requires  evaluating only the hyperobjective $F$ (as opposed to estimating the hypergradient in \cref{alg: OIGRM}).
The motivation for designing such an algorithm, is that
while evaluating $\nabla F$ up to $\alpha$ accuracy requires $O(\alpha^{-1})$ gradient evaluations, the hyperobjective value can be estimated at a linear rate, as stated in \cref{{lem:ZeroOrderApprox}}. 
Our algorithms are based on the Lipschitzness of $F$, which we prove below. 
\lemLipscConstrBilevel*
\begin{proof}
    By Lemma $2.1$ of \cite{ghadimi2018approximation}, the hypergradient of $F$ computed with respect to the variable $x$ may be expressed as $\nabla_x F(x) = \nabla_x f(x, y^\ast(x)) + \left(\frac{d y^\ast(x)}{dx}\right)^\top \cdot \nabla_y f(x, y^\ast(x))$. Since we impose Lipschitzness on $f$ and $y^*$, we can bound each of the terms of $\nabla_x F(x)$ by the claimed bound. 
\end{proof}

\subsection{Proof of \cref{thm: Lipschitz-min-with-inexact-zero-oracle}}

Denoting the uniform randomized smoothing $F_\rho(x):=\E_{\norm{z}\leq1}[F(x+\rho\cdot z)]$ where the expectation, here and in what follows, is taken with respect to the uniform measure, it is well known \citep[Lemma 10]{shamir2017optimal} that
\begin{align}
\E_{\norm{w}=1}\left[\tfrac{d}{2\rho}(F(x+\rho w)-F(x-\rho w)) w\right]
&=\nabla F_\rho( x)~,
\nonumber
\\
\E_{ \norm{w}=1}\norm{\nabla F_\rho( x)-\tfrac{d}{2\rho}( F( x+\rho w)- F( x-\rho w)) w}^2
&\lesssim dL^2~.\label{eq:grad_var_bound}
\end{align}
We first show that replacing the gradient estimator with the inexact evaluations $\tF(\cdot)$ leads to a biased gradient estimator of $ F$.

\begin{lemma}\label{lem: inexact gradient}
Suppose $| F(\cdot)-\tF(\cdot)|\leq\alpha$. Denoting
\begin{align*}
 g_x&=\tfrac{d}{2\rho}( F( x+\rho w)- F( x-\rho w)) w~,
\\
\tbg_x&=\tfrac{d}{2\rho}(\tF(x+\rho w)-\tF( x-\rho w)) w~,
\end{align*}
it holds that
\[
\E_{ \norm{w}=1}\norm{ g_x-\tbg_x}\leq\frac{\alpha d}{\rho}~,
~~~~~\text{and}~~~~~
\E_{\norm{w}=1}\norm{\tbg_x}^2\lesssim \frac{\alpha^2 d^2}{\rho^2}+dL^2~.
\]
\end{lemma}

\begin{proof}
For the first bound, we have
\[
\E_{\norm{w}=1}\norm{g_x-\tbg_x}
\leq \frac{d}{2\rho}(2\alpha)\E_{\norm{w}=1}\norm{w}
=\frac{\alpha d}{\rho}~,
\]
while for the second bound
\[
\E_{\norm{w}=1}\norm{\tbg_x}^2
=\E_{\norm{w}=1}\norm{\tbg_x-g_x+g_x}^2
\leq 2\E_{\norm{w}=1}\norm{\tbg_x-g_x}^2+2\E_{\norm{w}=1}\norm{g_x}^2
\lesssim \frac{d^2}{\rho^2}\cdot\alpha^2+dL^2~,
\]
where the last step invoked \cref{eq:grad_var_bound}.
\end{proof}

\begin{proof}[Proof of \cref{thm: Lipschitz-min-with-inexact-zero-oracle}]
We denote $\alpha'=\frac{\alpha d}{\rho},~\widetilde{G}=\sqrt{\frac{\alpha^2 d^2}{\rho^2}+dL^2}$. Since $x_t=x_{t-1}+\Delta_{t}$, we have
\begin{align*}
 F_\rho(x_t)- F_\rho(x_{t-1})
&=\int_{0}^{1}\inner{\nabla F_{\rho}(x_{t-1}+s\Delta_t),\Delta_t}ds
\\
&=\E_{s_t\sim\Unif[0,1]}\left[\nabla F_{\rho}(x_{t-1}+s_t{\Delta}_t),\Delta_t\right]
\\&=\E\left[\inner{\nabla F_{\rho}(z_{t}),\Delta_t}\right]~.
\end{align*}
By summing over $t\in[T]=[K\times M]$, we get for any fixed sequence $u_1,\dots,u_K\in\reals^d:$
\begin{align*}
    \inf  F_\rho\leq  F_\rho(x_T)
    &\leq  F_\rho(x_0)+\sum_{t=1}^{T}\E\left[\inner{\nabla F_{\rho}(z_{t}),\Delta_t}\right]
    \\&= F_\rho(x_0)+\sum_{k=1}^{K}\sum_{m=1}^{M}\E\left[\inner{\nabla F_{\rho}(z_{(k-1)M+m}),\Delta_{(k-1)M+m}-u_k}\right]\\
    &~~~~+\sum_{k=1}^{K}\sum_{m=1}^{M}\E\left[\inner{\nabla F_{\rho}(z_{(k-1)M+m}),u_k}\right]
    \\
    &\leq  F_\rho(x_0)+\sum_{k=1}^{K}\mathrm{Reg}_M(u_k)+\sum_{k=1}^{K}\sum_{m=1}^{M}\E\left[\inner{\nabla F_{\rho}(z_{(k-1)M+m}),u_k}\right]
    \\
    &\leq F_\rho(x_0)+KD\widetilde{G}\sqrt{M}+K\alpha'DM+\sum_{k=1}^{K}\sum_{m=1}^{M}\E\left[\inner{\nabla F_{\rho}(z_{(k-1)M+m}),u_k}\right]
\end{align*}
where the last inequality follows by combining \cref{lem: inexact gradient} and \cref{lem: inexact OGD}.
By setting $u_k:=-D\frac{\sum_{m=1}^{M}\nabla F_{\rho}(z_{(k-1)M+m})}{\norm{\sum_{m=1}^{M}\nabla F_{\rho}(z_{(k-1)M+m})}}$, rearranging and dividing by $DT=DKM$ we obtain
\begin{align}
\frac{1}{K}\sum_{k=1}^{K}\E\norm{\frac{1}{M}\sum_{m=1}^{M}\nabla F_{\rho}(z_{(k-1)M+m})}
&\leq \frac{ F_\rho(x_0)-\inf F_\rho}{DT}+\frac{\widetilde{G}}{\sqrt{M}}+\alpha'
\nonumber\\
&=\frac{ F_\rho(x_0)-\inf F_\rho}{K\nu}+\frac{\sqrt{\frac{\alpha^2 d^2}{\rho^2}+L^2d}}{\sqrt{M}}+\frac{\alpha d}{\rho}  
\nonumber\\
&\leq
\frac{ F_\rho(x_0)-\inf F_\rho}{K\nu}
+\frac{{\frac{\alpha d}{\rho}}}{\sqrt{M}}
+\frac{L\sqrt{d}}{\sqrt{M}}
+\frac{\alpha d}{\rho}~. \label{eq: bound eps}
\end{align}
Finally, note that for all $m\in[M]:\norm{z_{(k-1)M+m}-\overline{x}_{k}}\leq M D\leq\nu$, therefore 
$\nabla F_{\rho}(z_{(k-1)M+m})\in\partial_\nu F_\rho(\overline{x}_{k})\subset \partial_{\delta}F(\overline{x}_{k})$, where the last containment is due to \citep[Lemma 4]{kornowski2024algorithm} by using our assignment $\rho+\nu= \delta$.
Invoking the convexity of the Goldstein subdifferential, this implies that
\[
\frac{1}{M}\sum_{m=1}^{M}\nabla F_{\rho}(z_{(k-1)M+m})
\in\partial_{\delta}F(\overline{x}_{k})
~,
\]
thus it suffices to bound the first three summands in \eqref{eq: bound eps} by $\epsilon$ in order to finish the proof.
This happens as long as $
\frac{ F_\rho(x_0)-\inf F_\rho}{K\nu}\leq\frac{\epsilon}{3}$, $\frac{{\frac{\alpha d}{\rho}}}{\sqrt{M}}\leq\frac{\epsilon}{3}$, and $\frac{L\sqrt{d}}{\sqrt{M}}\leq\frac{\epsilon}{3}$, which imply $K\gtrsim \frac{ F_\rho(x_0)-\inf F_\rho}{\nu\epsilon}$, $ M\gtrsim \frac{\alpha^2 d^2}{\rho^2\epsilon^2}$, and $M\gtrsim \frac{L^2 d}{\epsilon^2}~$. 
By our assignments of $\rho$ and $\nu$, these  result in
\begin{align*}
T=KM
&=O\left(\frac{ F_\rho(x_0)-\inf F_\rho}{\nu\epsilon}\cdot \left(\frac{\alpha^2 d^2}{\rho^2\epsilon^2}+\frac{L^2 d}{\epsilon^2}\right)\right)
\\
&=O\left(\frac{(F(x_0)-\inf F)d}{\delta\epsilon^3}\cdot \left(\frac{\alpha^2 d}{\rho^2}+L^2\right)\right)
\\
&=O\left(\frac{(F(x_0)-\inf F)d}{\delta\epsilon^3}\cdot \left(\alpha^2 d\cdot\max\left\{\frac{1}{\delta^2},\frac{L^2}{(F(x_0)-\inf F)^2}\right\}+L^2\right)\right)
~,
\end{align*}
completing the proof of convergence guarantee for \cref{alg: IZO}.

\end{proof}

\subsection{Proof of \cref{thm:Lipschitz-min-with-inexact-grad-oracle}}

We recall  \cref{thm:Lipschitz-min-with-inexact-grad-oracle} below to keep this section self-contained. 

\thmLipscMinWithInexactGradOracle*

Our analysis is inspired by the reduction from online learning to nonconvex optimization given by \cite{cutkosky2023optimal}.
To that end, we start by proving a seemingly unrelated result, asserting that online gradient descent minimizes the regret with respect to inexact evaluations. Recalling standard definitions from online learning, given a sequence of linear losses $\ell_m(\cdot)=\inner{g_m,\cdot}$, if an algorithm chooses ${\Delta}_1,\dots,{\Delta}_M$ we denote the regret with respect to $u$ as
\[
\mathrm{Reg}_M(u):=\sum_{m=1}^{M}\inner{g_m,{\Delta}_m-u}.
\]
Consider an update rule according to online projected \emph{inexact} gradient descent:
\[
\Delta_{m+1}:=\mathrm{clip}_{D}(\Delta_{m}-\eta_m\tbg_m).
\]

\begin{lemma}[Inexact Online Gradient Descent] \label{lem: inexact OGD}
In the setting above, suppose that $(\tbg_m)_{m=1}^{M}$ are possibly randomized vectors, such that
$\E\norm{\tbg_m-g_m}\leq\alpha$ and $\E\norm{\tbg_m}^2\leq\widetilde{G}^2$ for all $m\in[M]$.
Then for any $\norm{u}\leq D$ it holds that
\[
\E\left[\mathrm{Reg}_M(u)\right]\leq \frac{D^2}{\eta_M}+\widetilde{G}^2\sum_{m=1}^{M}\eta_m+ \alpha DM~,
\]
where the expectation is with respect to the (possible) randomness of $(\tbg_m)_{m=1}^{M}$.
In particular, setting $\eta_m\equiv\frac{D}{\widetilde{G}\sqrt{M}}$ yields
\[
\E\left[\mathrm{Reg}_M(u)\right]\lesssim D\widetilde{G}\sqrt{M}+\alpha DM~.
\]
\end{lemma}

\begin{proof}
For any $m\in[M]:$
\begin{align*}
\norm{\Delta_{m+1}-u}^2
&=\norm{\mathrm{clip}_{D}(\Delta_{m}-\eta_m\tbg_m)-u}^2
\\&\leq \norm{\Delta_{m}-\eta_m\tbg_m-u}^2
=\norm{\Delta_{m}-u}^2+\eta_m^2\norm{\tbg_m}^2-2\eta_m\inner{\Delta_{m}-u,\tbg_m}~,
\end{align*}
thus
\[
\inner{\tbg_m,\Delta_m-u}
\leq \frac{\norm{\Delta_m-u}^2-\norm{\Delta_{m+1}-u}^2}{2\eta_m}+\frac{\eta_m}{2}\norm{\tbg_m}^2~,
\]
from which we get that
\begin{align*}
\E\inner{g_m,\Delta_m-u}
&=\E\inner{\tbg_m,\Delta_m-u}+
\E\inner{g_m-\tbg_m,\Delta_m-u}
\\
&\leq\frac{\norm{\Delta_m-u}^2-\norm{\Delta_{m+1}-u}^2}{2\eta_m}+\frac{\eta_m}{2}\E\norm{\tbg_m}^2+\E\norm{g_m-\tbg_m}\cdot \norm{\Delta_m-u}
\\&\leq\frac{\norm{\Delta_m-u}^2-\norm{\Delta_{m+1}-u}^2}{2\eta_m}+\frac{\eta_m}{2}\widetilde{G}^2+\alpha D~.
\end{align*}
Summing over $m\in[M]$, we see that
\begin{align*}
\E\left[\mathrm{Reg}_M(u)\right]
&\leq \sum_{m=1}^{M}\norm{\Delta_m-u}^2\left(\frac{1}{\eta_m}-\frac{1}{\eta_{m-1}}\right)+\frac{\widetilde{G}^2}{2}\sum_{m=1}^{M}\eta_m+M\alpha D
\\
&\leq \frac{D^2}{\eta_M}+\widetilde{G}^2\sum_{m=1}^{M}\eta_m+ \alpha DM~.
\end{align*}
The simplification for $\eta_m\equiv\frac{D}{\widetilde{G}\sqrt{M}}$ readily follows.
\end{proof}

We are now ready to analyze \cref{alg: OIGRM} in the inexact gradient setting.
\begin{proof}[Proof of \cref{thm:Lipschitz-min-with-inexact-grad-oracle}]
Since \cref{alg: OIGRM}  has $x_t=x_{t-1}+\Delta_{t}$, we have
\begin{align*}
F(x_t)-F(x_{t-1})
&=\int_{0}^{1}\inner{\nabla F(x_{t-1}+s\Delta_t),\Delta_t}ds
\\
&=\E_{s_t\sim\Unif[0,1]}\left[\langle\nabla F(x_{t-1}+s_t{\Delta}_t),\Delta_t\rangle\right]
\\&=\E\left[\inner{\nabla F(z_{t}),\Delta_t}\right]~. 
\end{align*}
By summing over $t\in[T]=[K\times M]$, we get for any fixed sequence $u_1,\dots,u_K\in\reals^d:$
\begin{align*}
    \inf F\leq F(x_T)
    &\leq F(x_0)+\sum_{t=1}^{T}\E\left[\inner{\nabla F(z_{t}),\Delta_t}\right]
    \\&=F(x_0)+\sum_{k=1}^{K}\sum_{m=1}^{M}\E\left[\inner{\nabla F(z_{(k-1)M+m}),\Delta_{(k-1)M+m}-u_k}\right]\\
    &~~~~+\sum_{k=1}^{K}\sum_{m=1}^{M}\E\left[\inner{\nabla F(z_{(k-1)M+m}),u_k}\right]
    \\
    &\leq  F(x_0)+\sum_{k=1}^{K}\mathrm{Reg}_M(u_k)+\sum_{k=1}^{K}\sum_{m=1}^{M}\E\left[\inner{\nabla F(z_{(k-1)M+m}),u_k}\right]
    \\
    &\leq F(x_0)+KD\widetilde{G}\sqrt{M}+K\alpha DM+\sum_{k=1}^{K}\sum_{m=1}^{M}\E\left[\inner{\nabla F(z_{(k-1)M+m}),u_k}\right]
\end{align*}
where the last inequality follows from \cref{lem: inexact OGD}
for $\widetilde{G}=\sqrt{L^2+\alpha^2},~\eta=\frac{D}{\widetilde{G}\sqrt{M}}$, since 
$\norm{\tbg_t-\nabla F(z_t)}\leq\alpha$ (deterministically) for all $t\in[T]$ by assumption.
Letting $u_k:=-D\frac{\sum_{m=1}^{M}\nabla F(z_{(k-1)M+m})}{\norm{\sum_{m=1}^{M}\nabla F(z_{(k-1)M+m})}}$, rearranging and dividing by $DT=DKM$, we obtain
\begin{align}
\frac{1}{K}\sum_{k=1}^{K}\E\norm{\frac{1}{M}\sum_{m=1}^{M}\nabla F(z_{(k-1)M+m})}
&\leq \frac{F(x_0)-\inf F}{DT}+\frac{\widetilde{G}}{\sqrt{M}}+\alpha
\nonumber\\
&=\frac{F(x_0)-\inf F}{K\delta}
+\frac{\widetilde{G}}{\sqrt{M}}+\alpha~.
\label{eq: bound eps_first}
\end{align}
Finally, note that for all $k\in[K],m\in[M]:\norm{z_{(k-1)M+m}-\overline{x}_{k}}\leq M D\leq\delta$, therefore
$\nabla F(z_{(k-1)M+m})\in\partial_\delta F(\overline{x}_{k})$. Invoking the convexity of the Goldstein subdifferential, we see that
\[
\frac{1}{M}\sum_{m=1}^{M}\nabla F(z_{(k-1)M+m})\in\partial_\delta  F(\overline{x}_{k})~,
\]
thus it suffices to bound the first two summands on the right-hand side in \cref{eq: bound eps_first} by $\epsilon$ in order to finish the proof. This happens as long as $
\frac{F(x_0)-\inf F}{K\delta}\leq\frac{\epsilon}{2}$ and $\frac{\widetilde{G}}{\sqrt{M}}\leq\frac{\epsilon}{2}$. These are  equivalent to $ K\geq \frac{2(F(x_0)-\inf F)}{\delta\epsilon}$ and $M\geq\frac{4\widetilde{G}^2}{\epsilon^2}$, 
which results in \[T=KM=O\left(\frac{F(x_0)-\inf F}{\delta\epsilon}\cdot \frac{L^2+\alpha^2}{\epsilon^2}\right)=O\left(\frac{(F(x_0)-\inf F)L^2}{\delta\epsilon^3}\right),\] completing the proof.
\end{proof}

\subsection{Proof of \cref{thm:practical_Lipschitz-min-with-inexact-grad-oracle}}

Throughout the proof we denote $z_t=x_{t}+\delta\cdot w_t$.
Since $F$ is $L$-Lipschitz, $F_\delta(x):=\E_{w\sim\mathrm{Unif}(\S^{d-1})}[F(x+\delta\cdot w)]$ is $L$-Lipschitz and $O(L\sqrt{d}/\delta)$-smooth. By smoothness we get
\begin{align*}
F_{\delta}(x_{t+1})-F_{\delta}(x_{t})
&\leq \inner{\nabla F_{\delta}(x_t),x_{t+1}-x_{t}}+O\left(\frac{L\sqrt{d}}{\delta}\right)\cdot\norm{x_{t+1}-x_{t}}^2
\\
&=-\eta\inner{\nabla F_{\delta}(x_t),\tbg_t}+O\left(\frac{\eta^2 L\sqrt{d}}{\delta}\right)\cdot\norm{\tbg_t}^2
\\
&=-\eta\inner{\nabla F_{\delta}(x_t),\nabla F(z_t)}-\eta\inner{\nabla F_{\delta}(x_t),\tbg_t-\nabla F(z_t)}+O\left(\frac{\eta^2 L\sqrt{d}}{\delta}\right)\cdot\norm{\tbg_t}^2~.
\end{align*}

Noting that $\E[\nabla F(z_t)]=\nabla F_\delta(x_t)$
and that
$\norm{\tbg_t}\leq
\norm{\tbg_t-\nabla F(z_t)}+\norm{\nabla F(z_t)}
\leq \alpha+L$, we see that
\begin{align*}
    \E[F_{\delta}(x_{t+1})-F_{\delta}(x_{t})]
    \leq -\eta \E\norm{\nabla F_\delta (x_t)}^2+\eta L\alpha+O\left(\frac{\eta^2 L\sqrt{d}}{\delta}(\alpha+L)^2\right)~,
\end{align*}
which implies
\[
\E\norm{\nabla F_\delta (x_t)}^2 \leq \frac{\E[F_\delta(x_t)]-\E[F_\delta(x_{t+1})]}{\eta}+L\alpha+O\left(\frac{\eta L\sqrt{d}(\alpha+L)^2}{\delta}\right)~.
\]
Averaging over $t=0,\dots,T-1$ 
and noting that $F_\delta(x_0)-\inf F_\delta \leq (F(x_0)-\inf F) +\delta L$ results in
\[
\E\norm{\nabla F_\delta (x^{\out})}^2=\frac{1}{T}\sum_{t=0}^{T-1}E\norm{\nabla F_\delta (x_t)}^2
\leq \frac{(F(x_0)-\inf F)+\delta L}{\eta T}+L\alpha+O\left(\frac{\eta L\sqrt{d}(\alpha+L)^2}{\delta}\right)~.
\]
By Jensen's inequality and the sub-additivity of the square root,
\[
\E\norm{\nabla F_\delta (x^{\out})} \leq \sqrt{\frac{(F(x_0)-\inf F)+\delta L}{\eta T}}+\sqrt{L\alpha}+O\left(\sqrt{\frac{\eta L\sqrt{d}(\alpha+L)^2}{\delta}}\right)~.
\]
Setting $\eta=\frac{\sqrt{((F(x_0)-\inf F)+\delta L)\delta}}{\sqrt{T L\sqrt{d}(\alpha+L)^2}}$ yields the final bound
\[
\E\norm{\nabla F_\delta (x^{\out})} \lesssim
\frac{((F(x_0)-\inf F)+\delta L)^{1/4}L^{1/4}d^{1/8}(\alpha+L)^{1/2}}{\delta^{1/4}T^{1/4}}+\sqrt{L\alpha}~,
\]
and the first summand is bounded by $\epsilon$ for $T=O\left(\frac{((F(x_0)-\inf F)+\delta L) L \sqrt{d}(L+\alpha)^2}{\delta\epsilon^4}\right)$.

\section{Proofs of \cref{sec:inequality-bilevel}}

\subsection{Reformulation equivalence}\label{appendix:reformulation-equivalence}
\begin{restatable}[Reformulation equivalence]{theorem}{reformulation}\label{thm:reformulation_equivalence}
When $\lambda^*$ matches to an optimal dual solution to the lower level problem $y^* = \arg\min_y g(x,y) ~\text{s.t.} ~h(x,y) \leq 0$, we show that for each $x$, the reformulation has the same feasible region of $y$. 
\end{restatable}
\begin{proof}
 We first show that lower-level feasibility implies feasibility of the reformulated problem.     
    Let $y^*, \lambda^* = \min\limits_y \max\limits_{\beta \geq 0} g(x,y) + \beta^\top h(x,y)$ be the primal and the dual solution to the lower level problem with parameter $x$.
    We can verify that $y^*$ satisfies all the constraints in the reformulation problem. The feasibility condition $h(x,y^*)$ is automatically satisfied.
    We just need to check:
    \begin{align*}
        g^*(x) & \coloneqq  \min\limits_\theta g(x,\theta) + (\lambda^*)^\top h(x,\theta) \\
        & = g(x,y^*) + (\lambda^*)^\top h(x,y^*).\numberthis\label[eq]{eq:g_gamma_star_equals_g_plus_gamma_h}
    \end{align*}
    Therefore, $x, y^*$ is a feasible point to the reformulation problem.
    
    We now show the other direction, i.e., that feasibility of the reformulaed problem implies that of the lower-level problem. 
    Given $\lambda^*$, let us assume $y$ satisfies $g(x,y) \leq g^*_{\lambda^*}(x)$ and $h(x,y) \leq 0$.
    On the other hand, assume $y^*, \lambda^* = \min\limits_y \max\limits_{\beta \geq 0} g(x,y) + \beta^\top h(x,y)$ be the primal and the dual solution.
    We can show that:
    \begin{align}
        g(x,y) + (\lambda^*)^\top h(x,y) \leq g^*(x) \coloneqq \min_\theta g(x,\theta) + (\lambda^*)^\top h(x,\theta). 
    \end{align}
    By the strong convexity of $g + (\lambda^*)^\top h$, we know that $y$ matches to the unique minimum $y^*$, which implies that $y = y^*$ is also a feasible point to the original bilevel problem.
\end{proof}

\subsection{Active constraints in differentiable optimization}\label{sec:inactive-constraints-in-differentiable-optimization}
By computing the derivative of the KKT conditions in \cref{sec:differentiable-optimization}, we get:
\begin{align}
    (\nabla^2_{yx} g + (\lambda^*)^\top \nabla_{yx}^2 h) + (\nabla^2_{yy} g + (\lambda^*)^\top \nabla_{yy}^2 h) \frac{dy^*}{dx} + (\nabla_y h)^\top \frac{d\lambda^*}{dx} &= 0 \label{eqn:appendix-KKT1} \\
    \text{diag}(\lambda^*) \nabla_x h + \text{diag}(\lambda^*) \nabla_y h \frac{dy^*}{dx} + \text{diag}(h) \frac{d\lambda^*}{dx} &= 0. \label{eqn:appendix-KKT2}
\end{align}

Let $\mathcal{I} = \{ i \in [d_h] | h(x,y^*)_i = 0, \lambda^*_i > 0 \}$ be the set of active constraints with positive dual solution, and $\mathcal{I}_1 = \{ i | h(x,y^*)_i \neq 0 \}$ be the set of inactive constraints and $\mathcal{I}_2 = \{ i | h(x,y^*)_i = 0, \lambda^*_i = 0 \}$. We know that $\bar{\mathcal{I}} = \mathcal{I}_1 \cup \mathcal{I}_2$. For each $i \in \mathcal{I}_1$, due to complementary slackness, we know that $\lambda^*_i = 0$. 

For $i \in \mathcal{I}_1$ in \cref{eqn:appendix-KKT1}, we have $\lambda^*_i \nabla_x h(x,y^*)_i + \lambda^*_i \nabla_y h(x,y^*)_i \frac{dy^*}{dx} + h(x,y^*)_i \frac{d \lambda^*_i}{dx} = 0$, which implies $h(x,y^*)_i \frac{d \lambda^*_i}{dx} = 0$ because  $\lambda^*_i = 0$. This in turn implies $\frac{d \lambda^*_i}{dx} = 0 $ because $h(x,y^*)_i < 0$. 
That means the dual variable $\lambda^*_i = 0$ and has zero gradient $\frac{d \lambda^*_i}{dx} = 0$ for any index $i \in \mathcal{I}_1$.
Therefore, we can remove row $i \in \mathcal{I}_1$ in \cref{eqn:appendix-KKT2} and obtain $\lambda^*_i = 0$ and $\frac{d \lambda^*_i}{dx} = 0$.

For $i \in \mathcal{I}_2$, the KKT condition in \cref{eqn:appendix-KKT2} is  degenerate. Therefore, $\frac{d \lambda^*_i}{d x}$ can be arbitrary, i.e., non-differentiable.
As a subgradient choice, we can set $\frac{d \lambda^*_i}{d x} = 0$ for such $i$. 
This choice will also eliminate its impact on the KKT condition in \cref{eqn:appendix-KKT1} because $\frac{d \lambda^*_i}{d x}$ is set to be $0$.
By this choice of subgradient, we can also remove row $i \in \mathcal{I}_2$ \cref{eqn:appendix-KKT2}.

Thus \cref{eqn:appendix-KKT2} can be written as the following set of equations, for  $h_\mathcal{I} = [h_i]_{i \in \mathcal{I}}$ and $\lambda^*_\mathcal{I} = [\lambda^*_i]_{i\in \mathcal{I}}$:
\begin{align}
    & \text{diag}(\lambda^*) \nabla_x h_\mathcal{I} + \text{diag}(\lambda^*_\mathcal{I}) \nabla_y h_\mathcal{I} \frac{dy^*}{dx} + \text{diag}(h_\mathcal{I}) \frac{d\lambda^*_\mathcal{I}}{dx} = 0 \nonumber \\
    \Longrightarrow \quad & \text{diag}(\lambda^*) \nabla_x h_\mathcal{I} + \text{diag}(\lambda^*_\mathcal{I}) \nabla_y h_\mathcal{I} \frac{dy^*}{dx} = 0 \quad \text{(due to $h_\mathcal{I}(x,y^*) = 0$)}. \label{eqn:appendix-new-KKT2}
\end{align}

In \cref{eqn:appendix-KKT1}, due to $\frac{d \lambda^*_i}{dx} = 0$ for all $i \in \bar{\mathcal{I}}$, we can remove $\frac{d \lambda^*_i}{dx} ~\forall i \in \bar{\mathcal{I}}$ in \cref{eqn:appendix-KKT1} by:
\begin{align}
    0 = & ~ (\nabla^2_{yx} g + (\lambda^*)^\top \nabla_{yx}^2 h) + (\nabla^2_{yy} g + (\lambda^*)^\top \nabla_{yy}^2 h) \frac{dy^*}{dx} + (\nabla_y h)^\top \frac{d\lambda^*}{dx} \nonumber \\
    = & ~ (\nabla^2_{yx} g + (\lambda^*)^\top \nabla_{yx}^2 h) + (\nabla^2_{yy} g + (\lambda^*)^\top \nabla_{yy}^2 h) \frac{dy^*}{dx} + (\nabla_y h_\mathcal{I})^\top \frac{d\lambda^*_\mathcal{I}}{dx}. 
    \label{eqn:appendix-new-KKT1} 
\end{align}

Combining \cref{eqn:appendix-new-KKT1} and \cref{eqn:appendix-new-KKT2}, we get:
\begin{align*}
    (\nabla^2_{yx} g + (\lambda^*)^\top \nabla_{yx}^2 h) + (\nabla^2_{yy} g + (\lambda^*)^\top \nabla_{yy}^2 h) \frac{dy^*}{dx} + (\nabla_y h_\mathcal{I})^\top \frac{d\lambda^*_\mathcal{I}}{dx} &= 0 \\
    \text{diag}(\lambda^*) \nabla_x h_\mathcal{I} + \text{diag}(\lambda^*_\mathcal{I}) \nabla_y h_\mathcal{I} \frac{dy^*}{dx} &= 0, 
\end{align*}
which can be written in its matrix form:
\begin{align}\label{eqn:kkt-system_appendix}
\begin{bmatrix}
\nabla^2_{yy} g + (\lambda^*)^\top \nabla_{yy}^2 h & \nabla_y h_\mathcal{I}^\top \\
\text{diag}(\lambda^*_\mathcal{I}) \nabla_y h_\mathcal{I} & 0
\end{bmatrix}
\begin{bmatrix}
    \frac{dy^*}{dx} \\
    \frac{d\lambda^*_\mathcal{I}}{dx}
\end{bmatrix}
= 
-
\begin{bmatrix}
    \nabla^2_{yx} g + (\lambda^*)^\top \nabla_{yx}^2 h \\
    \text{diag}(\lambda^*_\mathcal{I}) \nabla_x h_\mathcal{I}
\end{bmatrix}
\end{align}
This concludes the derivation of the derivative of constrained optimization in \cref{eqn:kkt-system}.

\subsection{Inequality case: bounds on primal solution error and constraint violation}
\solutionApproximation*

\begin{proof}
We first provide the claimed bound on $\|y^*_{\alpha_1, \alpha_2} - y^*(x)\|$. 

\noindent\textbf{Part 1: Bound on the convergence of $y$.}

Since $y_{\lambda^*,\boldsymbol{\alpha}}^*$ minimizes $\mathcal{L}_{\boldsymbol{\alpha}, \lambda^*}(x,y)$, the first-order condition gives us:
\begin{align*}
    0 = \nabla_y \mathcal{L}_{\boldsymbol{\alpha}, \lambda^*}(x,y_{\lambda^*,\boldsymbol{\alpha}}^*). 
\end{align*}
Similarly, we can compute the gradient of $\mathcal{L}_{\boldsymbol{\alpha}, \lambda^*}(x,y)$ at $y^*$: 
\begin{align*}
    \nabla_y \mathcal{L}_{\alpha}(x,y^*) & = \nabla_y f(x,y^*) + \alpha_1(\nabla_y g(x,y^*) + (\lambda^*)^\top \nabla_y h(x,y^*)) + \alpha_2 \nabla_y h_\mathcal{I}(x,y^*)^\top h_\mathcal{I}(x,y^*)  \\
    & = \nabla_y f(x,y^*),
\end{align*} 
where the second step is due to the property of the primal and dual solution: $\nabla_y g(x,y^*) + (\lambda^*)^\top \nabla_y h(x,y^*) = 0$ by the stationarity condition in the KKT conditions, and by definition of the active constraints $h_\mathcal{I}$ where the optimal $y^*$ must have $h_\mathcal{I}(x,y^*) = 0$.

Since, for a sufficiently large $\alpha_1$, the penalty function is $\alpha_1 \mu_g - L_f \geq \frac{\alpha_1 \mu_g}{2}$  strongly convex in $y$, we have:
\begin{align*}
    \frac{\alpha_1 \mu_g}{2} \norm{y^* - y_{\lambda^*,\boldsymbol{\alpha}}^*} \leq \norm{\nabla_y\mathcal{L}_{\boldsymbol{\alpha}, \lambda^*}(x,y^*) - \nabla_y\mathcal{L}_{\boldsymbol{\alpha}, \lambda^*}(x,y_{\lambda^*,\boldsymbol{\alpha}}^*)} = \norm{\nabla_y f(x,y^*)} \leq L_f.
\end{align*}
Therefore, upon rearranging the terms,  we obtain the claimed bound:
\begin{align*}
    \norm{y^* - y^*_{\boldsymbol{\alpha}, \lambda^*}} \leq \frac{2 L_f}{\alpha_1 \mu_g}.
\end{align*}

\noindent\textbf{Part 2: bound on the constraint violation.}

When we plug $y^*$ into \cref{eqn:penalty-lagrangian}, we get:
\begin{align*}
    \mathcal{L}_{\boldsymbol{\alpha}, \lambda^*}(x,y^*) & = f(x,y^*) + \alpha_1 (g(x,y^*) + (\lambda^*)^\top h(x,y^*) - g^*_{\lambda^*}(x)) + \frac{\alpha_2}{2} \norm{h_\mathcal{I}(x,y^*)}^2 = f(x,y^*).
\end{align*}
Plugging in $y^*_{\boldsymbol{\alpha},\lambda^*}$, we may obtain: 
\begin{align*}
    \mathcal{L}_{\boldsymbol{\alpha}, \lambda^*}(x,y_{\lambda^*,\boldsymbol{\alpha}}^*) & = f(x,y_{\lambda^*,\boldsymbol{\alpha}}^*) + \alpha_1 (g(x,y_{\lambda^*,\boldsymbol{\alpha}}^*) + (\lambda^*)^\top h(x,y_{\lambda^*,\boldsymbol{\alpha}}^*) - g^*(x)) + \frac{\alpha_2}{2} \norm{h_\mathcal{I}(x,y_{\lambda^*,\boldsymbol{\alpha}}^*)}^2 \\
    & = f(x,y_{\lambda^*,\boldsymbol{\alpha}}^*) + \alpha_1 (  g(x,y_{\lambda^*,\boldsymbol{\alpha}}^*) + (\lambda^*)^\top h(x,y_{\lambda^*,\boldsymbol{\alpha}}^*) - g(x,y^*) - (\lambda^*)^\top h(x,y^*) ) \\
    & \qquad + \frac{\alpha_2}{2} \norm{h_\mathcal{I}(x,y_{\lambda^*,\boldsymbol{\alpha}}^*)}^2 \\
    & \geq f(x,y_{\lambda^*,\boldsymbol{\alpha}}^*) + \alpha_1 \mu_g \norm{y^* - y_{\lambda^*,\boldsymbol{\alpha}}^*}^2 + \frac{\alpha_2}{2} \norm{h_\mathcal{I}(x,y_{\lambda^*,\boldsymbol{\alpha}}^*)}^2,
\end{align*} where we used the strong convexity (with respect to $y$) of $g(x,y)+(\lambda^*)^\top h(x,y)$ and the optimality of $y^*$ for $g(x,y)+(\lambda^*)^\top h(x,y)$. 
By the optimality of $y_{\lambda^*,\boldsymbol{\alpha}}^*$ for $\mathcal{L}_{\boldsymbol{\alpha}, \lambda^*}$, we know that 
\begin{align*}
    f(x,y^*) = \mathcal{L}_{\boldsymbol{\alpha}, \lambda^*}(x,y^*) \geq \mathcal{L}_{\boldsymbol{\alpha}, \lambda^*}(x,y_{\lambda^*,\boldsymbol{\alpha}}^*) \geq f(x,y_{\lambda^*,\boldsymbol{\alpha}}^*) + \alpha_1 \mu_g \norm{y^* - y_{\lambda^*,\boldsymbol{\alpha}}^*}^2 + \frac{\alpha_2}{2} \norm{h_\mathcal{I}(x,y_{\lambda^*,\boldsymbol{\alpha}}^*)}^2.
\end{align*}
Therefore, by the Lipschitzness of the function $f$ in terms of $y$, and the bound $\|y^* - y_{\lambda^*,\boldsymbol{\alpha}}^*\| \leq \frac{2L_f}{\alpha_1 \mu_g}$, we know that:\begin{align*}
    \frac{\alpha_2}{2} \norm{h_\mathcal{I}(x,y_{\lambda^*,\boldsymbol{\alpha}}^*)}^2 & \leq f(x,y^*) - f(x,y_{\lambda^*,\boldsymbol{\alpha}}^*) - \alpha_1 \mu_g \norm{y^* - y_{\lambda^*,\boldsymbol{\alpha}}^*}^2 \\
    & \leq L_f \norm{y^* - y_{\lambda^*,\boldsymbol{\alpha}}^*} - \alpha_1 \mu_g \norm{y^* - y_{\lambda^*,\boldsymbol{\alpha}}^*}^2 \\
    & \leq L_f \norm{y^* - y_{\lambda^*,\boldsymbol{\alpha}}^*} \\
    & = O({\alpha_1^{-1}}). 
\end{align*}
Rearranging terms then gives the claimed bound. 
\end{proof}
The bound on the constraint violation in \cref{{thm:solution-bound}} is an important step in the following theorem.

\subsection{Proof of \cref{thm:diff_in_hypergrad_and_gradLagr}: gradient approximation for inequality constraints}\label{appendix:proof-of-inexact-gradient}
\gradientApproximation*
\begin{proof}
First, we recall \cref{eqn:penalty-lagrangian} here: \[\mathcal{L}_{\lambda^*,\boldsymbol{\alpha}}(x,y) = f(x,y) + \alpha_1 \left( g(x,y) + (\lambda^*)^\top h(x,y) - g^*(x)  \right) + \frac{\alpha_2}{2} \norm{h_\mathcal{I}(x,y)}^2.\] Next, recall from \cref{eq:g_gamma_star_equals_g_plus_gamma_h}, we can express $g^*(x)=g(x,y^*)+(\lambda^*)^\top h(x,y^*)$, which we use in the first step below:  
\begin{align}
    \footnotesize
    & \nabla_x F(x) - \frac{d}{dx} \mathcal{L}_{\lambda^*,\boldsymbol{\alpha}}(x,y_{\lambda^*,\alpha}^*) \nonumber \\
    = & \left( \nabla_x f(x,y^*) + \frac{d y^*}{d x}^\top \nabla_y f(x,y^*) \right) - \Biggl( \nabla_x f(x,y_{\lambda^*,\boldsymbol{\alpha}}^*) + \alpha_1 (\nabla_x g(x,y_{\lambda^*,\boldsymbol{\alpha}}^*) + \nabla_x h(x,y_{\lambda^*,\boldsymbol{\alpha}}^*)^\top \lambda^* \nonumber \\
    &  - \alpha_1(\nabla_x g(x,y^*) + \nabla_x h(x,y^*)^\top \lambda^*) 
    + \alpha_2
    \nabla_x h_\mathcal{I}(x,y_{\lambda^*,\boldsymbol{\alpha}}^*)^\top h_\mathcal{I}(x,y_{\lambda^*,\boldsymbol{\alpha}}^*)
\Biggl) \nonumber \\
    = & \nabla_x f(x,y^*) - \nabla_x f(x,y_{\lambda^*,\boldsymbol{\alpha}}^*) \label{eqn:f-difference} \\
    & + \frac{d y^*}{d x}^\top \nabla_y f(x,y^*) - \frac{d y^*}{d x}^\top \nabla_y f(x,y_{\lambda^*,\boldsymbol{\alpha}}^*) \label{eqn:df-difference} \\ 
    & + \frac{d y^*}{d x}^\top \nabla_y f(x,y_{\lambda^*,\boldsymbol{\alpha}}^*) - \underbrace{\alpha_1 \begin{bmatrix}
        \nabla^2_{yx} g + (\lambda^*)^\top \nabla_{yx}^2 h \\
        \text{diag}(\lambda^*_\mathcal{I}) \nabla_x h_\mathcal{I}
    \end{bmatrix}^\top  \begin{bmatrix}
        y_{\lambda^*,\boldsymbol{\alpha}}^* - y^*  \\
        0
    \end{bmatrix} }_{\text{added term 1}} \nonumber\\
    &\qquad- \underbrace{\alpha_2 \begin{bmatrix}
        \nabla^2_{yx} g + (\lambda^*)^\top \nabla_{yx}^2 h \lambda^* \\
        \text{diag}(\lambda^*_\mathcal{I}) \nabla_x h_\mathcal{I}
    \end{bmatrix}^\top \begin{bmatrix}
        0 \\
        \text{diag}(1/\lambda^*_\mathcal{I}) h_\mathcal{I}(x,y_{\lambda^*,\boldsymbol{\alpha}}^*) 
    \end{bmatrix}}_{\text{added term 2}} \label{eqn:df-and-added-term} \\
    & + \alpha_1 \Biggl( \nabla_x g(x,y^*) - \nabla_x g(x,y_{\lambda^*,\boldsymbol{\alpha}}^*) + \nabla_x h(x, y^*)^\top \lambda^* - \nabla_x h(x, y_{\lambda^*,\boldsymbol{\alpha}}^*)^\top \lambda^* \nonumber\\
    &\qquad+ \underbrace{\begin{bmatrix}
        \nabla^2_{yx} g + (\lambda^*)^\top \nabla_{yx}^2 h \\
        \text{diag}(\lambda^*_\mathcal{I}) \nabla_x h_\mathcal{I}
    \end{bmatrix}^\top \begin{bmatrix}
        y_{\lambda^*,\boldsymbol{\alpha}}^* - y^*  \\
        0 
    \end{bmatrix}}_{\text{added term 1}}
    \Biggl)
    \label{eqn:dgdh-and-added-term} \\
    & - \alpha_2 \nabla_x h_\mathcal{I}(x,y_{\lambda^*,\boldsymbol{\alpha}}^*)^\top h_\mathcal{I}(x,y_{\lambda^*,\boldsymbol{\alpha}}^*) + \underbrace{\alpha_2 \begin{bmatrix}
        \nabla^2_{yx} g + (\lambda^*)^\top \nabla_{yx}^2 h \lambda^* \\
        \text{diag}(\lambda^*_\mathcal{I}) \nabla_x h_\mathcal{I}
    \end{bmatrix}^\top \begin{bmatrix}
        0 \\
        \text{diag}(1/\lambda^*_\mathcal{I}) h_\mathcal{I}(x,y_{\lambda^*,\boldsymbol{\alpha}}^*) 
    \end{bmatrix}}_{\text{added term 2}}. \label{eqn:dh2-difference} 
\end{align}
According to \cref{eqn:kkt-system} and \cref{eqn:kkt-system_appendix}, we let 
\[H = \begin{bmatrix}
    \nabla^2_{yy} g + (\lambda^*)^\top \nabla_{yy}^2 h & \nabla_y h_\mathcal{I}^\top \\
    \text{diag}((\lambda^*_\mathcal{I}) \nabla_y h_\mathcal{I} & 0
\end{bmatrix},\] which is invertible by \cref{item:assumption_tangen_space} and by the fact that we remove all the inactive constraints.
We now bound the terms in \cref{eqn:f-difference}, \cref{eqn:df-difference}, \cref{eqn:df-and-added-term}, \cref{eqn:dgdh-and-added-term}, and \cref{eqn:dh2-difference}.

\noindent\textbf{Bounding \cref{eqn:f-difference} and \cref{eqn:df-difference}:}
\cref{eqn:f-difference} can be easily bounded by the smoothness of $f$ 
in terms of $x$ and $y$, and the bound on $\|y^*- y_{(\lambda^*,\boldsymbol{\alpha}}^*\| \leq O({\alpha_1^{-1}})$ from \cref{thm:solution-bound}. Therefore, we know:
\begin{align*}
    \norm{\nabla_x f(x,y^*) - \nabla_x f(x,y_{(\lambda^*,\boldsymbol{\alpha}}^*)} \leq C_f \norm{y^* - y_{(\lambda^*,\boldsymbol{\alpha}}^*} \leq C_f\cdot O({\alpha_1^{-1}}). 
\end{align*}
Similarly, given \cref{item:assumption_safe_constraints}  by which $y^*(x)$ is $L_y$-Lipschitz in $x$, we have the bound $\norm{\frac{d y^*}{d x}} \leq L_y$. Therefore, \cref{eqn:df-difference} can be bounded by:
\begin{align*}
    \norm{\frac{dy^*}{dx}^\top \nabla_y f(x,y^*) - \frac{dy^*}{dx}^\top \nabla_y f(x,y_{(\lambda^*,\boldsymbol{\alpha}}^*)} \leq C_f \norm{\frac{dy^*}{dx}} \norm{y^* - y_{(\lambda^*,\boldsymbol{\alpha}}^*} \leq C_f L_y\cdot O({\alpha_1^{-1}}). 
\end{align*}

\noindent\textbf{Bounding \cref{eqn:df-and-added-term}:}

Using~\cref{eqn:kkt-system} to solve $\begin{bmatrix}
    \frac{d y^*}{d x} \\
    \frac{d (\lambda^*}{d x}
\end{bmatrix} = -H^{-1} \begin{bmatrix}
    \nabla^2_{yx} g + ((\lambda^*)^\top \nabla_{yx}^2 h \\
    \text{diag}((\lambda^*_\mathcal{I}) \nabla_x h_\mathcal{I}
\end{bmatrix}$, 
we can write:
\begin{align*}
    & \frac{d y^*}{d x}^\top \nabla_y f(x,y_{(\lambda^*,\boldsymbol{\alpha}}^*) =  \begin{bmatrix}
        \nabla^2_{yx} g + ((\lambda^*)^\top \nabla_{yx}^2 h \\
        \text{diag}((\lambda^*_\mathcal{I}) \nabla_x h_\mathcal{I}
    \end{bmatrix}^\top (H^{-1})^\top \begin{bmatrix}
        - \nabla_y f(x,y_{(\lambda^*,\boldsymbol{\alpha}}^*) \\ 0
    \end{bmatrix}   \\
     &= -\frac{d y^*}{d x}^\top  \Biggl( \alpha_1 \begin{bmatrix}
        \nabla_y g(x,y_{(\lambda^*,\boldsymbol{\alpha}}^*) + \nabla_y h(x,y_{(\lambda^*,\boldsymbol{\alpha}}^*)^\top (\lambda^* \\
        0
    \end{bmatrix} \\ 
    &\quad + 
    \alpha_2 \begin{bmatrix}
         \nabla_y h_\mathcal{I}(x,y_{(\lambda^*,\boldsymbol{\alpha}}^*)^\top h_\mathcal{I}(x,y_{(\lambda^*,\boldsymbol{\alpha}}^*) \\
        0
    \end{bmatrix}
    \Biggl), \numberthis\label{eqn:dfdy_dydx_expansion}
\end{align*}
where we use the optimality of $y_{(\lambda^*,\boldsymbol{\alpha}}^*$ from \cref{{eq:def_y_lambda_star}}:
\begin{align*}\numberthis\label{eqn:y_lambda_star_optimality}
    & \nabla_y f(x,y_{(\lambda^*,\boldsymbol{\alpha}}^*) + \alpha_1 \left( \nabla_y g(x,y_{(\lambda^*,\boldsymbol{\alpha}}^*) + \nabla_y h(x,y_{(\lambda^*,\boldsymbol{\alpha}}^*)^\top (\lambda^* \right) \\
    &\quad+ \alpha_2 \nabla_y h_\mathcal{I}(x,y_{(\lambda^*,\boldsymbol{\alpha}}^*)^\top h_\mathcal{I}(x,y_{(\lambda^*,\boldsymbol{\alpha}}^*) = 0. 
\end{align*}
Further, recall that $H$ is non-degenerate by \cref{assumption:linEq_smoothness}, as a result of which, the added term 1 in \cref{eqn:df-and-added-term} can be modified as follows: 
\begin{align}
    \footnotesize
    & \begin{bmatrix}
        \nabla^2_{yx} g + ((\lambda^*)^\top \nabla_{yx}^2 h \\
        \text{diag}((\lambda^*_\mathcal{I}) \nabla_x h_\mathcal{I}
    \end{bmatrix}^\top \begin{bmatrix}
        \alpha_1 (y_{(\lambda^*,\boldsymbol{\alpha}}^* - y^*) \\
        0
    \end{bmatrix} \nonumber \\
    = & \begin{bmatrix}
        \nabla^2_{yx} g + ((\lambda^*)^\top \nabla_{yx}^2 h \\
        \text{diag}((\lambda^*_\mathcal{I}) \nabla_x h_\mathcal{I}
    \end{bmatrix}^\top (H^{-1})^\top  H^\top  \begin{bmatrix}
        \alpha_1 (y_{(\lambda^*,\boldsymbol{\alpha}}^* - y^*) \\
        0
    \end{bmatrix} \nonumber \\
    = & \alpha_1 
    \begin{bmatrix}
        \nabla^2_{yx} g + ((\lambda^*)^\top \nabla_{yx}^2 h \\
        \text{diag}((\lambda^*_\mathcal{I}) \nabla_x h_\mathcal{I}
    \end{bmatrix}^\top (H^{-1})^\top \begin{bmatrix}
        (\nabla^2_{yy} g + ((\lambda^*)^\top \nabla_{yy}^2 h)^\top (y_{(\lambda^*,\boldsymbol{\alpha}}^* - y^*) \\
        \nabla_y h_\mathcal{I}(x,y^*)(y_{(\lambda^*,\boldsymbol{\alpha}}^* - y^*)
    \end{bmatrix}. \label{eqn:added_term}
\end{align}

The added term 2 in \cref{eqn:df-and-added-term} can be expanded to: 
\begin{align}
    & \alpha_2 \begin{bmatrix}
        \nabla^2_{yx} g + ((\lambda^*)^\top \nabla_{yx}^2 h (\lambda^* \\
        \text{diag}((\lambda^*_\mathcal{I}) \nabla_x h_\mathcal{I}
    \end{bmatrix}^\top \begin{bmatrix}
        0 \\
        \text{diag}(1/(\lambda^*_\mathcal{I}) h_\mathcal{I}(x,y_{(\lambda^*,\boldsymbol{\alpha}}^*) 
    \end{bmatrix} \nonumber \\
    = & \alpha_2 \begin{bmatrix}
        \nabla^2_{yx} g + ((\lambda^*)^\top \nabla_{yx}^2 h (\lambda^* \\
        \text{diag}((\lambda^*_\mathcal{I}) \nabla_x h_\mathcal{I}
    \end{bmatrix}^\top (H^{-1})^\top H^\top  \begin{bmatrix}
        0 \\
        \text{diag}(1/(\lambda^*_\mathcal{I}) h_\mathcal{I}(x,y_{(\lambda^*,\boldsymbol{\alpha}}^*) 
    \end{bmatrix} \nonumber \\
    = & \alpha_2  \begin{bmatrix}
        \nabla^2_{yx} g + ((\lambda^*)^\top \nabla_{yx}^2 h (\lambda^* \\
        \text{diag}((\lambda^*_\mathcal{I}) \nabla_x h_\mathcal{I}
    \end{bmatrix}^\top (H^{-1})^\top \begin{bmatrix}
        \nabla_y h_\mathcal{I}(x,y^*)^\top h_\mathcal{I}(x,y_{(\lambda^*,\boldsymbol{\alpha}}^*) \\
        0
    \end{bmatrix} \label{eqn:added-term-2}
\end{align}    
Therefore, we can compute the difference between \cref{eqn:dfdy_dydx_expansion}, \cref{eqn:added_term}, and \cref{eqn:added-term-2} to bound \cref{eqn:df-and-added-term}, and use the fact that $\nabla_y g(x,y^*) + (\lambda^*)^\top \nabla_y h(x,y^*) = 0$: 
\begin{align}
    \footnotesize
    & \frac{d y^*}{d x}^\top \nabla_y 
 f(x,y_{\lambda^*,\boldsymbol{\alpha}}^*) - \text{added term 1 } - \text{added term 2}  \nonumber \\
    = & \begin{bmatrix}
    \nabla^2_{yx} g + (\lambda^*)^\top \nabla_{yx}^2 h \\
        \text{diag}(\lambda^*_\mathcal{I}) \nabla_x h_\mathcal{I}
    \end{bmatrix}^\top (H^{-1})^\top \biggr( \alpha_1 \begin{bmatrix}
        \nabla_y g(x,y_{\lambda^*,\boldsymbol{\alpha}}^*) - \nabla_y g(x,y^*) - \nabla^2_{yy} g(x,y^*) (y_{\lambda^*,\boldsymbol{\alpha}}^* - y^*) \\ 
        0
        \end{bmatrix}\label{eqn:g_second_order_difference} \\
     & +
    \alpha_1 \begin{bmatrix}
        \nabla_y h(x,y_{\lambda^*,\boldsymbol{\alpha}}^*)^\top \lambda^* - \nabla_y h(x,y^*)^\top \lambda^* - \nabla^2_{yy} h(x,y^*)^\top \lambda^* (y_{\lambda^*,\boldsymbol{\alpha}}^* - y^*) \\
        0
    \end{bmatrix}\label{eqn:h_second_order_difference}
    \\
    & \quad -
    \alpha_1 \begin{bmatrix}
        0 \\ 
        \nabla_y h_\mathcal{I}(x,y^*) (y_{\lambda^*,\boldsymbol{\alpha}}^* - y^*)
    \end{bmatrix} \label{eqn:constraint_difference}
    \\
    & \quad \quad + \alpha_2  \begin{bmatrix}
        \nabla_y h_\mathcal{I}(x,y_{\lambda^*,\boldsymbol{\alpha}}^*)^\top h_\mathcal{I}(x,y_{\lambda^*,\boldsymbol{\alpha}}^*)  \\
        0 
    \end{bmatrix} - \begin{bmatrix}
        \nabla_y h_\mathcal{I}(x,y^*)^\top h_\mathcal{I}(x,y_{\lambda^*,\boldsymbol{\alpha}}^*) \\
        0
    \end{bmatrix} \biggr)
    \label{eqn:h2-difference}
\end{align}
The terms in \cref{eqn:g_second_order_difference} and \cref{eqn:h_second_order_difference} can both be bounded by $\alpha_1 C_{g} L_y \|y_{\lambda^*,\boldsymbol{\alpha}}^* - y^*\|^2$ and $\alpha_1 R C_{h} L_y \|y_{\lambda^*,\boldsymbol{\alpha}}^* - y^*\|^2$ by the smoothness of $g$ and $h^\top \lambda^*$. Further, plugging in   $\|y^*- y_{\lambda^*,\boldsymbol{\alpha}}^*\| \leq O({\alpha_1^{-1}})$ from \cref{thm:solution-bound} bounds both these terms by $O({\alpha_1^{-1}})$.

To bound the term in \cref{eqn:constraint_difference}, we use:
\begin{align*}
    \norm{h_\mathcal{I}(x,y_{\lambda^*,\boldsymbol{\alpha}}^*) - h_\mathcal{I}(x,y^*) - \nabla_y h_\mathcal{I}(x,y^*) (y_{\lambda^*,\boldsymbol{\alpha}}^* - y^*)} \leq C_h \norm{y_{\lambda^*,\boldsymbol{\alpha}}^* - y^*}^2. 
\end{align*}
Therefore, we have: 
\begin{align*}
     \norm{\nabla_y h_\mathcal{I}(x,y^*) (y_{\lambda^*,\boldsymbol{\alpha}}^* - y^*)} & \leq  \norm{h_\mathcal{I}(x,y_{\lambda^*,\boldsymbol{\alpha}}^*)} +\norm{h_\mathcal{I}(x,y^*)} + C_h O(\norm{y_{\lambda^*,\boldsymbol{\alpha}}^* - y^*}^2) \\
    & \leq O({{\alpha_1^{-1/2} \alpha_2^{-1/2}}}) + 0 + O({\alpha_1^{-2}}) \\
    & = O({{\alpha_1^{-1/2} \alpha_2^{-1/2}}} + {\alpha_1^{-2}}),
\end{align*}
which upon scaling by $\alpha_1$ gives us the following bound on the term in \cref{eqn:constraint_difference}:
\begin{align*}
    \alpha_1 \norm{\nabla_y h_\mathcal{I}(x,y^*) (y_{\lambda^*,\boldsymbol{\alpha}}^* - y^*)} \leq O(\alpha_1^{1/2} \alpha_2^{-1/2} + \alpha_1^{-1}) ~.
\end{align*}

The term in \cref{eqn:h2-difference} can be bounded by:
\begin{align*}
    & \alpha_2 \norm{\nabla_x h_\mathcal{I}(x,y_{\lambda^*,\boldsymbol{\alpha}}^*)^\top h_\mathcal{I}(x,y_{\lambda^*,\boldsymbol{\alpha}}^*) - \nabla_x h_\mathcal{I}(x,y^*)^\top h_\mathcal{I}(x,y_{\lambda^*,\boldsymbol{\alpha}}^*)} \\
    = & \alpha_2 \norm{\nabla_x h_\mathcal{I}(x,y_{\lambda^*,\boldsymbol{\alpha}}^*) - \nabla_x h_\mathcal{I}(x,y^*)} O(\norm{h_\mathcal{I}(x,y^*_{\boldsymbol{\alpha},\lambda^*})})  \\
    = & \alpha_2\cdot O({\alpha_1^{-1}}) O({{\alpha_1^{-1/2} \alpha_2^{-1/2}}}) \\
    = & O(\alpha_1^{-3/2} \alpha_2^{1/2}) \numberthis\label{eq:bound_on_h_plus_grady_h}
\end{align*}

\noindent\textbf{Bounding \cref{eqn:dgdh-and-added-term}:}
This can be easily bounded by the smoothness of $g$ and $h$, and the bound on the dual solution $\norm{\lambda^*} \leq R$.
Thus \cref{eqn:dgdh-and-added-term} can be bounded by $R \cdot O({\alpha_1^{-1}}) = O({\alpha_1^{-1}})$.

\noindent\textbf{Bounding \cref{eqn:dh2-difference}:}
By the same argument in \cref{eq:bound_on_h_plus_grady_h}, we get:
\begin{align}
    & \alpha_2 \norm{\nabla_y h_\mathcal{I}(x,y_{\lambda^*,\boldsymbol{\alpha}}^*)^\top h_\mathcal{I}(x,y_{\lambda^*,\boldsymbol{\alpha}}^*) - \nabla_y h_\mathcal{I}(x,y^*)^\top h_\mathcal{I}(x,y_{\lambda^*,\boldsymbol{\alpha}}^*) } \nonumber \\
    \leq & \alpha_2 \norm{\nabla_y h_\mathcal{I}(x,y_{\lambda^*,\boldsymbol{\alpha}}^*) - \nabla_y h_\mathcal{I}(x,y^*)} \norm{h_\mathcal{I}(x,y_{\lambda^*,\boldsymbol{\alpha}}^*)} \nonumber \\
    = &  \alpha_2\cdot O({\alpha_1^{-1}}) O({{\alpha_1^{-1/2} \alpha_2^{-1/2}}}) \nonumber \\
    = & O(\alpha_1^{-3/2}\alpha_2^{1/2}) ~. \nonumber
\end{align}

Combining all upper bounds gives the claimed bound. 
\end{proof}

\subsection{Proof of the main result (\cref{thm:cost_of_computing_ystar_gammastar_inequality}): convergence and computation cost}\label{appendix:cost_of_computing_ystar_gammastar_inequality}
\computationCostInequality*
\begin{proof}
    First, given the bound in \cref{thm:diff_in_hypergrad_and_gradLagr}, we choose $\alpha_1 = \alpha^{-2}$ and $\alpha_2 = \alpha^{-4}$ to ensure the inexactness of the gradient oracle is bounded by $\alpha$. In the later analysis, we will still use $\alpha_1$ and $\alpha_2$ in the penalty function for clarity.

    Now we estimate the computation cost of the inexact gradient oracle:
    
    \noindent\textbf{Lower-level problem.}
    Given the oracle access to the optimal dual solution $\lambda^*(x)$, we can recover the primal solution $y^*(x)$ efficiently (e.g, by \cite{zhang2022solving}). Therefore, we can use the primal and dual solutions to construct the penalty function $\mathcal{L}_{\lambda^*, \alpha}(x,y)$ in ~\cref{eqn:penalty-lagrangian}.
    
    \noindent\textbf{Penalty function minimization problem.}
    The second main optimization problem is the penalty minimization problem in \cref{line:lagrangian-optimization} of \cref{alg:inexact-gradient-oracle}.
    Recall from \cref{eqn:penalty-lagrangian} that 
    \begin{align}\label{eqn:penalty-lagrangian-approximate}
        \mathcal{L}_{\lambda,\boldsymbol{\alpha}}(x,y) = f(x,y) + \alpha_1 \left( g(x,y) + \lambda^\top h(x,y) - g^*(x)  \right) + \frac{\alpha_2}{2} \norm{h_\mathcal{I}(x,y)}^2, 
    \end{align}
    where we use the approximate dual solution $\lambda$ as opposed to the optimal dual solution $\lambda^*$.
    Given \cref{eqn:penalty-lagrangian-approximate}, we solve the penalty minimization problem: 
    \begin{align*}
        y'_{\lambda, \alpha}(x) \coloneqq \arg\min_y \mathcal{L}_{\lambda,\boldsymbol{\alpha}}(x,y). 
    \end{align*}
    The penalty minimization is a unconstrained strongly convex optimization problem, which is known to have linear convergence rate. We further analyze its convexity and smoothness below to precisely estimate the computation cost:
    \begin{itemize}
        \item The strong convexity of $\mathcal{L}_{\lambda,\boldsymbol{\alpha}}(x,y)$ is lower bounded by $\frac{\alpha_1 \mu_g}{2} = O(\alpha_1)$.
        \item The smoothness of $\mathcal{L}_{\lambda,\boldsymbol{\alpha}}(x,y)$ is dominated by the smoothness of $\alpha_2 \norm{h_\mathcal{I}(x,y)}^2$ since $\alpha_2 \gg \alpha_1$. By \cref{thm:diff_in_hypergrad_and_gradLagr}, we know that the optimal solution must lie in an open ball $B(y^*, O(1/\alpha_1))$ with center $y^*$ (inner optimization primal solution) and a radius of the order of $O(\frac{1}{\alpha_1})$. This implies that we just need to search over a bounded feasible set of $y$, which we can bound $\norm{\nabla_y h(x,y)} \leq L_h$ and $h(x,y) \leq H$ within the bounded region $y \in B(y^*, O(1/\alpha_1))$. We can show that $h^2$ is smooth (gradient Lipschitz) within the bounded region by the following:
        \begin{align}
            \norm{\nabla^2_{yy} h^2} = \norm{h \nabla^2_{yy} h + \nabla_y h^\top \nabla_y h} \leq \norm{h \nabla^2_{yy} h} + \norm{\nabla_y h^\top \nabla_y h} \leq H C_h + L_h^2 \nonumber
        \end{align}
        which also implies $h^2_{\mathcal{I}}$ is also smooth (gradient Lipschitz).
        Therefore, $\alpha h^2_{\mathcal{I}}$ is $( H C_h + L_h^2) \alpha_2 = O(\alpha_2)$ smooth.
    \end{itemize}
     
    Choosing $\alpha_1 = \frac{1}{\alpha^2}$ and $\alpha_2 = \frac{1}{\alpha^4}$, the condition number of $\mathcal{L}_{\boldsymbol{\alpha}, \lambda}(x,y)$ becomes $\kappa = O(\alpha_2 / \alpha_1) = O(\frac{1}{\alpha^2})$.
    Therefore, by the linear convergence of gradient descent in strongly convex smooth optimization, the number of iterations needed to get to $\alpha$ accuracy is $O(\sqrt{\alpha^{-2}} \times \log (\frac{1}{\alpha})) = O(\frac{1}{\alpha} \log (\frac{1}{\alpha}))$.
    Therefore, we can get a near optimal solution $y'_{\lambda,\alpha}$ with inexactness $\alpha$ in $O(\frac{1}{\alpha})$ oracle calls.

\noindent\textbf{Computation cost and results.}
Overall, for the inner optimization, we can invoke the efficient optimal dual solution oracle to get the optimal dual solution $\lambda^*(x)$ and recover the optimal primal solution $y^*(x)$ from there.
For the penalty minimization, we need $O(\frac{1}{\alpha})$ oracle calls to solve an unconstrained strongly convex smooth optimization problem to get to $\alpha$ accuracy.
In conclusion, combining everything in \cref{appendix:cost_of_computing_ystar_gammastar_inequality}, we run $O(\frac{1}{\alpha})$ oracle calls to obtain an $\alpha$ accurate gradient oracle to approximate the hyperobjective gradient $\nabla_x F(x)$. This concludes the proof of \cref{thm:cost_of_computing_ystar_gammastar_inequality}.
\end{proof}

\begin{remark}
    The following analysis  quantifies how the error in the optimal dual solution propagates to the inexact gradient estimate. This is not needed if such a dual solution oracle exists. But in practice, the oracle may come with some error, for which we bound the error.
\end{remark}

\noindent\textbf{Bounding the error propagation in error in dual solution and the penalty minimization.}
First, if we do not get an exact optimal dual solution, the error in the dual solution $\lambda$ with $\norm{\lambda - \lambda^*} \leq \alpha$ will slightly impact the analysis in \cref{thm:diff_in_hypergrad_and_gradLagr}. Specifically, in \cref{appendix:proof-of-inexact-gradient}, the approximate $\lambda$ will impact the inexact gradient $\nabla_x \mathcal{L}_{\lambda,\alpha}(x,y'_{\lambda,\alpha})$ computation and the analysis in \cref{eqn:dgdh-and-added-term} and \cref{eqn:y_lambda_star_optimality}.
    In \cref{eqn:dgdh-and-added-term}, to change $\lambda$ to $\lambda^*$, we get an additional error:
    \begin{align*}\numberthis\label{eqn:dual-convergence-merit-function-1}
        & \alpha_1 \biggl( \nabla_x h(x,y')^\top (\lambda - \lambda^*) - \nabla_x h(x,y'_{\lambda,\alpha})^\top (\lambda - \lambda^*) \biggl) \nonumber \\
        = & \alpha_1 (\nabla_x h(x,y') - \nabla_x h(x,y'_{\lambda,\alpha}))^\top (\lambda - \lambda^*) \nonumber \\
        \leq & \alpha_1 C_h \norm{y' - y'_{\lambda,\alpha}} (\lambda - \lambda^*) \nonumber \\
        \leq & O(\alpha_1 \alpha_1^{-1} \alpha) = O(\alpha), 
    \end{align*}
    where the last inequality is due to $\norm{y' - y'_{\lambda,\alpha}} \leq O(\alpha_1^{-1})$ that is based on a similar analysis in \cref{thm:solution-bound} with a near-optimal $y'_{\lambda,\alpha}$ under $\alpha^2 = \alpha_1$ accuracy.
    
    Therefore, the error incurred by inexact $\lambda$ in \cref{eqn:dgdh-and-added-term} is at most $O(\alpha)$, which is of the same rate as the current gradient inexactness $O(\alpha)$.
    
    In \cref{eqn:y_lambda_star_optimality}, the optimality holds approximately for the approximate $\lambda$. Therefore, by the near optimality of $y_{\lambda,\boldsymbol{\alpha}}'$ (strongly convex optimization), we know that the following gradient is also $\alpha$-close to $0$, i.e.,
    \begin{align*}\numberthis\label{eqn:inexact-optimality-with-inexact-lambda}
        & \|\nabla_y f(x,y_{\lambda,\boldsymbol{\alpha}}') + \alpha_1 \left( \nabla_y g(x,y_{\lambda,\boldsymbol{\alpha}}') + \nabla_y h(x,y_{\lambda,\boldsymbol{\alpha}}')^\top \lambda \right) \\
        &\quad + \alpha_2 \nabla_y h_\mathcal{I}(x,y_{\lambda,\boldsymbol{\alpha}}')^\top h_\mathcal{I}(x,y_{\lambda,\boldsymbol{\alpha}}')\| \leq \alpha, 
    \end{align*}
    whose inexactness matches the inexactness of the gradient oracle $\alpha$, and thus we do not incur additional order of inexactness here.
    
    Moreover, there is an additional error because we need $\lambda^*$ as opposed to a near-optimal $\lambda$ to make the analysis in \cref{appendix:proof-of-inexact-gradient} work. The error between using $\lambda$ and $\lambda^*$ in \cref{eqn:inexact-optimality-with-inexact-lambda} can be bounded by:
    \begin{align}\label{eqn:dual-convergence-merit-function-2}
        \norm{\nabla_y h(x,y'_{\lambda,\alpha})^\top (\lambda - \lambda^*)} \leq L_h \alpha, 
    \end{align}
    where we use the local Lipschitzness of the function $h$ in an open ball near $y^*$.
    Therefore, the additional error is also  $O(\alpha)$, which matches the inexactness of the inexact gradient oracle. 

Therefore, we conclude that in order to bound the inexactness of the gradient oracle, we just need an efficient inexact dual solution with $\alpha$ accuracy.

\subsection{Practical oracle to optimal (approximate) dual solution}
Here we discuss how practical the assumption on the oracle access to the optimal dual solution is.

For linear inequality constraint $h(x,y) = Ax - By - b$, the LL problem is a constrained strongly convex smooth optimization problem. 

To show that we can compute an approximate solution to the optimal dual solution for linear inequality constraints, we apply the result from ~\cite{du2019linear}:
\begin{corollary}[Application of Corollary 3.1 in \cite{du2019linear}]\label{cor:linear-convergence-linear-inequality-LL}
    When $h(x,y) = Ax - By + b$ is linear in $y$, the primal and dual solutions can be written as:
    \begin{align}
        & y^*, \lambda^* = \arg\min_y \max_\lambda g(x,y) + (\lambda^*)^\top h(x,y) = g(x,y) - (\lambda^*)^\top B y + R(x) \nonumber \\
        \iff & y^*, \lambda^* = \arg\min_y \max_\lambda g(x,y) - (\lambda^*)^\top B y
    \end{align}
    where $g$ is strongly convex in $y$ and $B$ is of full rank by \cref{assumption:linEq_smoothness}. According to Corollary 3.1 from \cite{du2019linear}, the primal-dual gradient method guarantees a linear convergence. More precisely, in $t = O(\log \frac{1}{\alpha})$ iterations, we get:
    \begin{align}
        \norm{y^t - y^*} \leq \alpha \text{ and } \norm{\lambda^t - \lambda^*} \leq \alpha. 
    \end{align}
\end{corollary}
Given \cref{cor:linear-convergence-linear-inequality-LL}, we can efficiently approximate the primal and dual solutions up to high accuracy with $O(\log \frac{1}{\alpha})$ oracle calls when the inequality constraints are linear. This gives us an efficient approximate oracle access to the dual solution.

\begin{remark}
    Under the assumption of an optimal dual solution oracle, all the analyses mentioned in ~\cref{sec:inequality-bilevel} hold for the general convex inequality constraints. However, the main technical challenge is that the dual solution oracle for general convex inequality cannot be guaranteed in practice.
    In fact, to the best of our knowledge, there is no iterate convergence in the dual solution $\lambda$ for general convex inequality constraints. Most of the literature in strongly-convex-concave saddle point convergence only guarantees dual solution convergence in terms of its duality gap or some other merit functions. We are not aware of any successful bound on the dual solution iterate convergence, which is an important research question to answer by itself.
    This is the main technical bottleneck for general convex inequality constraints as well.
\end{remark}
\begin{remark}
    On the other hand, we need the dual solution iterate convergence with rate $O(1/\alpha)$ to ensure the error to be bounded. But this is not a necessary condition. To ensure a bound on the error propagation, we just need to bound some forms of merit functions (\cref{eqn:dual-convergence-merit-function-1} and \cref{eqn:dual-convergence-merit-function-2}) of the dual solutions, which we believe that this is much more tractable than the actual iterate dual solution convergence. We leave this as a future direction and this will generalize the analysis from linear inequality constraints to general convex inequality constraints.
\end{remark}

\subsection{The role of $\lambda^*(x)$ in the derivative of Equation~\cref{eqn:penalty-lagrangian}}
Notice that Equation~\cref{eqn:penalty-lagrangian}, we treat the dual solution $\lambda^*(x)$ as a constant to define the penalty function derivative. Yet, the dual solution $\lambda^*(x)$ is in fact also a function of $x$. Therefore, in theory, we should also compute its derivative with respect to $x$.

However, notice that the following:
\begin{align}\label{eqn:lambda_derivative}
    \nabla_x (\lambda^*(x))^\top h(x,y) & = \nabla_x h(x,y)^\top \lambda^* + \frac{d \lambda(x)}{dx}^\top h(x,y)
\end{align}
The later term in Equation~\cref{eqn:lambda_derivative} can be divided into two cases:
\begin{itemize}
    \item For active constraint $i \in \mathcal{I}$ with $h(x,y^*) = 0$, we know that $y^*_{\lambda,\alpha}$ is close to $y^*$ by \cref{thm:solution-bound}. Therefore, the derivative $\norm{\frac{d \lambda(x)}{dx}^\top h(x,y^*_{\lambda,\alpha})} \leq L_h L_\lambda \alpha_1 = O(\alpha_1) = O(\alpha^2)$ by the local smoothness of $h$ near $y^*$ and the Lipschitzness assumption of $\lambda^*$ in \cref{item:assumption_safe_constraints}.
    \item For inactive constraint $i \in \bar{\mathcal{I}}$ and $\lambda^*_i > 0$, we can solve the KKT conditions and get $\frac{d \lambda(x)}{dx} = 0$. Therefore, the second term becomes $0$.
    \item For inactive constraint $i \in \bar{\mathcal{I}}$ and $\lambda^*_i = 0$, the KKT system degenerates and we need to use subgradient. By solving the KKT system, we find that $\frac{d \lambda(x)}{dx} = 0$ is a valid subgradient. Therefore, by choosing this subgradient, the second term also vanishes.
\end{itemize}
Therefore, we do not need to compute the derivative of $\lambda^*$ as the terms involved its derivative is negligible compared to other major terms.

\section{Experimental setup}\label{appendix:experiment-setup}
All experiments were run on a computing cluster with Dual Intel Xeon Gold 6226 CPUs @ 2.7 GHz and DDR4-2933 MHz DRAM. No GPU was used, and we used 1 core with 8GB RAM per instance of the experiment. The cutoff time for running the algorithms is set to be 6 hours.
All experiments were run and averaged over 10 different random seeds.
All parameters in the constrained bilevel optimization in \cref{sec:exp},  including the objective parameters and the constrain parameters, are randomly generated from a normal distribution with $0$-mean and standard deviation $1$.

For our fully first-order algorithm, we implement the ``implementation-friendly'' \cref{alg: PIGD}, where the inexact gradient oracle subroutine is provided by implementing \cref{alg:inexact-gradient-oracle}.
All algorithms are implemented in PyTorch~\cite{paszke2019pytorch} to compute gradients, and using Cvxpy~\cite{diamond2016cvxpy} to solve the LL problem and the penalty minimization problem.
We implement our fully first-order method based on the solutions returned by Cvxpy with certain accuracy requirement, and use PyTorch to compute the inexact gradient discussed in \cref{sec:inequality-bilevel}.
We implement the non-fully first-order method using the CvxpyLayer~\cite{agrawal2019differentiable}, which is a Cvxpy compatible library that can differentiate through the LL convex optimization problem.

\end{document}